\documentclass[opre,dblanonrev]{informs4}
\usepackage{stmaryrd}

\RequirePackage{tgtermes}
\RequirePackage{newtxtext}
\RequirePackage{bm}
\RequirePackage{endnotes}

\OneAndAHalfSpacedXI

\usepackage{natbib}
 \bibpunct[, ]{(}{)}{,}{a}{}{,}%
 \def\bibfont{\small}%

\EquationsNumberedThrough

\TheoremsNumberedThrough
\ECRepeatTheorems

\usepackage{algorithm}
\usepackage{algpseudocode}
\usepackage{amsfonts,amsmath,amstext,amssymb,amsopn}
\usepackage{dsfont}
\usepackage{multirow}
\usepackage{bbm}
\usepackage{tikz}
\usepackage[margin=1in]{geometry}
\usepackage{scalerel,stackengine}
\usetikzlibrary{arrows.meta,
            decorations.pathreplacing,
                    calligraphy,%
                positioning}
\usepackage{booktabs}
\usepackage{placeins}

\usepackage{amssymb}

\newenvironment{myproof}[1][Proof]{%
  \par\noindent\textit{#1.}\quad%
}{%
  \par\noindent\mbox{} \hfill$\square$\par
}

\makeatletter
\def\BState{\State\hskip-\ALG@thistlm}
\makeatother

\def\h#1{{#1}^1}%

\def\n1{\h{\mathscr{N}}}

\def\j{\textsf{j}}

\usepackage{makecell}

\newcommand{\Reg}{\mathrm{Reg}}

\long\def\theARTICLETOPLEFT{}
\def\theARTICLETOPRIGHT{}

\begin{document}

\RUNTITLE{Value of Experiments in Assortment Optimization}

\TITLE{The Value of Experiments in Assortment Optimization}

\ARTICLEAUTHORS{%
\AUTHOR{Mohammadreza Ahmadnejadsaein}
\AFF{School of Operations Research and Information Engineering, Cornell University,
\EMAIL{ma932@cornell.edu}}

\AUTHOR{Amine Bennouna}
\AFF{Kellogg School of Management, Northwestern University,
\EMAIL{amine.bennouna@kellogg.northwestern.edu}}

\AUTHOR{Omar El Housni}
\AFF{School of Operations Research and Information Engineering, Cornell University,
\EMAIL{oe46@cornell.edu}}

}

\ABSTRACT{%
{%

Retailers can learn customer preferences by experimenting with different assortments.
However, the cost of each experiment and the operational burden of repeatedly changing product offerings can make extensive experimentation impractical.
We study how much a retailer can improve its final assortment decision using customer observations from one (or a very small number of) fixed experimental assortments.
We consider the multinomial logit model, where product revenues are known but customer preferences are unknown.
In this setting, we jointly design the experimental assortment and the rule for selecting a final assortment from observed customer choices to minimize worst-case regret: the expected revenue loss of the final assortment relative to an optimal assortment under the true preferences, maximized over possible customer preferences.
We first characterize exactly the minimax regret and an optimal assortment when no experimentation is possible.
We then show, surprisingly, that offering a single experimental assortment to just one customer can already substantially reduce regret, depending on how the products' revenues compare.
We derive the exact minimax regret and an optimal experimental policy in this setting.
When a single fixed experimental assortment is offered to many customers, we give a geometric characterization of the limiting minimax regret and construct an asymptotically optimal policy.
Depending on product revenues, this single experiment can yield no improvement over no experimentation or nearly complete regret reduction. Finally, we extend our analysis to multiple experimental assortments chosen in advance, characterizing the limiting minimax regret, constructing asymptotically optimal policies, and identifying which collections of assortments allow worst-case regret to vanish as the number of customer observations grows.

}

}%


\maketitle

\section{Introduction}\label{sec:intro}

Retailers and marketplaces must decide which products to offer their customers. Evaluating products in isolation can be misleading because a product’s demand depends on the alternatives offered alongside it. Adding a product may divert purchases from higher-revenue alternatives, while removing one may cause customers to leave without making a purchase. Choice models capture these substitution effects by specifying how customers choose among the offered products and the no-purchase option. Combining these choice probabilities with product revenues allows the retailer to evaluate the expected revenue of each assortment. Assortment optimization then seeks the set of products that maximizes this expected revenue \citep{talluri2004}.

Implementing this approach requires information about customer preferences, represented by the parameters of the choice model. These parameters are typically estimated from historical data on the assortments offered and the choices customers made \citep{train2009}. When a new set of products is introduced, however, the retailer may have little or no relevant historical data. How should it select an assortment when product revenues are known but customer preferences are not?

Experimentation offers a way to inform this decision. Before committing to a final assortment, the retailer can offer customers a subset of products as an experimental assortment and observe their choices. These observations provide information about customer preferences that can guide the final assortment decision.

Learning customer preferences through experimental assortments has been studied extensively. A substantial literature considers sequential learning, in which the offered assortment can be adapted in response to observed customer choices, either to maximize cumulative revenue while learning or to identify an optimal assortment \citep{saure2013,agrawal2017thompson,agrawal2019,lee2026design}.
However, frequent assortment changes can be difficult or impossible to implement in
physical stores, as replacing products and rearranging shelves require time and labor and can disrupt operations.
A report from the data science subsidiary of Kroger, a large U.S. grocery retailer, notes that ``testing countless
iterations'' of retail changes ``in stores is impractical,'' and indicates that a single physical-store
experiment can take multiple months.\footnote{\raggedright84.51$^\circ$, ``How virtual in-store testing can save time and money, and accelerate the innovation pipeline.'' Available at \url{https://www.8451.com/knowledge-hub/case-studies/how-virtual-in-store-testing-can-save-time-and-money-and-accelerate-the/} (published October 1, 2024).}

Similar constraints arise in online retail. Online retailers selling seasonal products may need to select their product range and place supplier orders before the selling season begins. Procurement lead times and inventory commitments can then limit their ability to revise the assortment as customer choices are observed \citep{kok2008assortment}. When opportunities for testing before these commitments are limited, the retailer must decide which products to carry based on only a small number of assortment experiments.

These operational constraints motivate approaches that reduce the burden of assortment experimentation. \citet{dong2020switching} develop sequential learning policies that achieve low cumulative revenue regret while limiting assortment switches. \citet{yu2026experimental} design a  collection of experimental assortments, chosen in advance, for choice-model estimation and nest identification. Practical constraints in physical stores may nevertheless limit retailers to one or very few experimental assortments. In this setting, the goal is to obtain information that improves the final assortment decision, which may require less information than estimating a full choice model. How much can such limited experimentation improve the retailer's assortment decision?

In this paper, we consider retailers that can run a pilot with a single or a very small number of experimental
assortments before committing to the final products they will offer to their customers. We study how much such a restricted
experiment can improve the final decision and when the improvement is substantial. For a given experimental budget, we determine which experimental assortments best inform a high-revenue final assortment. We also develop optimal (or near-optimal) decision rules for selecting the final assortment based on observed customer choices.

\subsection{Contributions and Main Insights}\label{sec:intro:contribution}

We study these questions under the widely used multinomial logit (MNL) model (see Section~\ref{sec:model:choice}). In this model, each product has a preference weight that captures its appeal to customers. The retailer knows the product revenues but not the preference weights. Before collecting data, it chooses one or a few experimental assortments to offer to customers. It then observes their choices, including no-purchase outcomes, and uses these observations to select a final assortment. The experimental assortments and the decision rule mapping observations to a final assortment together constitute a policy.  We evaluate a policy by its \emph{worst-case regret}: the expected revenue loss of the final assortment relative to an optimal assortment chosen with full knowledge of the preference weights, taking the worst case over all possible weights. Our goal is to jointly design the experimental assortments and the decision rule to minimize this regret.

\vspace{2mm}
\noindent\textbf{The benchmark without experimentation: Optimal assortment and regret.}
In Section~\ref{sec:G0}, we first establish a benchmark: what regret can the retailer guarantee without experimentation, and which assortment should it offer? When customer preferences are unknown, offering every product may seem like a safe choice, but it can lead to substantial worst-case revenue loss. The retailer faces a trade-off: excluding a product risks losing revenue when customers strongly prefer it, while including a low-revenue product can draw purchases away from higher-revenue alternatives. Without information about preferences, neither possibility can be ruled out. Theorem~\ref{thm:G0} characterizes the minimax regret, and Corollary~\ref{cor:G0half} yields a simple optimal policy: offer exactly the products whose revenues exceed $r_1/2$, where $r_1$ denotes the highest product revenue. This policy guarantees regret at most $r_1/2$. We use the minimax regret without experimentation as the benchmark for measuring how much experimentation improves the final decision.

\vspace{2mm}
\noindent\textbf{Even a single customer observation can meaningfully reduce regret.}
Surprisingly, even a single customer observation from one experiment can meaningfully improve the final assortment decision. In Section~\ref{sec:G1}, we consider a retailer that offers one experimental assortment to just one customer and observes the customer's choice, including the possibility of no purchase. We characterize exactly the minimax regret in this setting and, therefore, the best attainable benefit of a single observation (Theorem~\ref{thm:one:opt}). This benefit depends strongly on the product revenues: across the examples we exhibit, the reduction in regret ranges from $0\%$ (no benefit over the no-observation benchmark) to values approaching $87.5\%$.
For example, with uniformly spaced product revenues, the regret reduction approaches $30.4\%$ as the number of products grows. We also identify an optimal policy with an intuitive structure: experiment with products above a revenue threshold. If a purchase is observed, retain these products and add any additional products above a second threshold. If no purchase is observed, the assortment may expand further according to a third threshold. We determine all three thresholds exactly as functions of the product revenues. Remarkably, the decision rule uses only whether a purchase occurs, even though the observation also reveals which product was purchased.

\vspace{2mm}
\noindent\textbf{The value of a single experiment.}
How much can be gained by collecting more customer observations from the same experimental assortment? In Section~\ref{sec:rep}, we study a retailer that repeatedly offers one assortment and uses the observed choices to select its final assortment. The benefit again depends strongly on the product revenues: a single experiment can provide no improvement over the benchmark without experimentation or yield relative regret reductions arbitrarily close to $100\%$.

We characterize exactly the limiting minimax regret through a simple geometric rule (Theorem~\ref{thm:rep:main}). Place all product revenues, together with zero, on a line and partition them into three nonempty consecutive groups: high revenues, intermediate revenues, and low revenues including zero. The diameter of each group is the difference between its largest and smallest values. The smallest possible value of the largest group diameter is exactly the limiting minimax regret. In particular, this regret remains positive for every fixed set of at least three distinct positive revenues, regardless of how many customer observations are collected.

Theorem~\ref{thm:rep:main} also shows that a minimizing partition yields a simple asymptotically minimax-optimal policy: experiment with the high-revenue products, retain them as the final assortment if the observed purchase fraction meets a threshold determined by the product revenues, and otherwise add the intermediate-revenue products. As in the single-observation setting, the policy uses only whether customers purchase, without using the identities of the purchased products. With $m$ customer observations, its worst-case regret is within $r_1/(2\sqrt{m})$ of the limiting minimax regret. The geometric characterization also explains when the remaining regret is small: when the revenues, together with zero, form three tight clusters, a single experiment can bring worst-case regret close to zero.

\vspace{2mm}
\noindent\textbf{The value of multiple experiments.}
Finally, in Section~\ref{sec:limit}, we extend this geometric characterization to any fixed number of experiments. The retailer chooses up to $d$ distinct experimental assortments in advance, which can be tested in parallel (for example, across multiple stores within the same retail group), and uses the resulting observations to select a final assortment. We characterize the limiting minimax regret by placing the product revenues, together with zero, on a line and partitioning them into at most $d+2$ nonempty consecutive groups. The smallest possible value of the largest group diameter is exactly the limiting minimax regret (Theorem~\ref{thm:distinct}). A minimizing partition also yields an asymptotically minimax-optimal policy using at most $d$ nested revenue-ordered experimental assortments. Thus, the regret with no experimentation corresponds to a partition of two groups, with a single experiment to three, and each additional experiment permits one more group, potentially reducing the remaining worst-case revenue loss.

We further give necessary and sufficient conditions for a set of experimental assortments to permit vanishing worst-case regret as customer observations increase (Theorem~\ref{thm:distinct:vanishing}). In particular, with $n$ products, at least $n-1$ distinct experimental assortments are necessary, and $n-1$ appropriately chosen assortments suffice.

\subsection{Related Literature}\label{sec:intro:literature}

\noindent
\textbf{Assortment optimization.}
This literature studies which products to offer to maximize expected revenue. Algorithms have been developed under multinomial logit \citep{talluri2004}, nested logit \citep{davis2014nested}, Markov chain \citep{blanchet2016markov,feldman2017markov}, and nonparametric choice models \citep{aouad2018ranking}, with extensions addressing cardinality \citep{rusmevichientong2010}, capacity \citep{desir2020constrained}, and visibility constraints \citep{barre2026visibility}. Applications also include assortment personalization \citep{golrezaei2014real,elhousni2023customization} and two-sided markets \citep{ashlagi2022twosided,aouad2023twosided,elhousni2024twosided,ahmadnejadsaein2025adaptive}. Most of these works study assortment optimization under a specified choice model with known parameters.

When the parameters of the underlying choice model are unknown, observed customer choices can inform assortment decisions. These observations may come from historical transactions or from assortment experiments conducted to learn customer preferences. Existing approaches differ in how observations are collected and whether performance is evaluated during learning or through a final assortment decision. We next review existing approaches to assortment optimization when customer preferences are unknown.

\vspace{2mm}
\noindent
\textbf{Sequential assortment learning.}
Sequential assortment learning studies how a firm learns customer preferences while earning revenue through a sequence of assortment decisions \citep{caro2007,saure2013,agrawal2019, saha2025rankbreaking}. For example, \citet{agrawal2019} study the MNL-bandit problem, in which assortment offers balance learning customer preferences and earning revenue, and develop an upper-confidence-bound policy with nonasymptotic cumulative-regret guarantees. \citet{agrawal2017thompson} develop a Thompson-sampling policy for the same problem and establish worst-case cumulative-regret bounds. For uncapacitated MNL models, \citet{chen2021optimal} guide assortment decisions by directly estimating the revenues of revenue-ordered assortments, without separately estimating preference weights. To address the operational burden of frequent assortment changes, \citet{dong2020switching} develop policies with guarantees on both cumulative regret and the number of assortment changes. These approaches allow later assortments to depend on observed choices and evaluate performance through revenue earned throughout learning. We study experiments chosen in advance and evaluate their value through the revenue regret of the final assortment decision.

\vspace{2mm}
\noindent
\textbf{Offline assortment learning.}
Offline assortment learning selects a final assortment using historical offers and customer choices, with methods that account for estimation uncertainty \citep{han2025observational,zhang2026contextual}. \citet{dong2025pasta} develop PASTA, which maximizes the minimum expected revenue over a likelihood-based confidence set of choice models. Their guarantees require coverage of an optimal assortment in the historical offer distribution. Under an MNL model with bounded preference weights, \citet{han2025observational} use pessimistic rank-breaking to obtain guarantees requiring only sufficient coverage of each product in an optimal assortment. These guarantees depend on the bound on preference weights; our analysis considers worst-case regret uniformly over unbounded weights.

\vspace{2mm}
\noindent
\textbf{Experimental design in assortment optimization}.
Related work designs experiments to estimate choice models or inform a final assortment decision. \citet{yu2026experimental} design nonadaptive assortment experiments for choice-model estimation and establish nest-identification guarantees for nested logit models under identifiability and separation conditions. \citet{lee2026design} develop experimental designs for MNL models with linear utilities and apply them to sequential best-assortment identification with a prescribed confidence level. In a Bayesian framework, \citet{araman2022diffusion} study experimentation and stopping before a final decision, with an application to MNL assortment selection. Under a diffusion approximation, they show that an optimal experiment can be chosen once and repeated until stopping. Our formulation is decision-focused: we design experiments directly to optimize the final assortment decision rather than to estimate choice-model parameters. In particular, we jointly choose the experimental assortments and the final decision rule under a fixed budget of distinct assortments selected in advance, with the objective of minimizing worst-case revenue regret.

\vspace{2mm}
\noindent
Robust assortment optimization maximizes the minimum expected revenue over a set of possible choice models \citep{rusmevichientong2012,sturt2026expertise}.
\citet{lu2026robust} study robust assortment learning when future choice probabilities may differ from those of the model generating the data. \citet{chen2025prediction} select assortments through robust optimization using the partial preference orders revealed by historical transactions, without fitting a choice model. They also analyze experimental designs that sample single-product assortments. These robust objectives differ from our minimax-regret criterion, which minimizes the largest revenue loss relative to an optimal assortment under each possible preference vector.

\vspace{2mm}
\noindent
\textbf{The value of information in decision-making.}
Our work also relates to studies that measure information through its effect on decision quality. 
For linear optimization with an unknown cost vector, \citet{bennouna2025data, bennounaDataInformativeness2026} characterize which experiments, modeled as exact linear measurements, suffice to identify an optimal decision. Although MNL assortment optimization admits a linear programming formulation, the unknown parameters appear in the constraints, and experiments give stochastic rather than exact observations.
\citet{liu2023marginal} study the marginal value of acquiring customer preference data for assortment personalization and use this value to guide data collection. Their observations provide customer utility information, whereas we observe choices from offered assortments.  Related pricing studies evaluate decisions based on valuation samples \citep{allouah2022samples} or the purchase probability at a single historical price, known exactly or through an interval \citep{allouah2023single}, using worst-case revenue ratios. We quantify the value of stochastic choice observations through minimax revenue regret and determine how this value depends on the number and composition of experimental assortments.

\section{Model and Problem Formulation}\label{sec:model}

We formalize our setting in this section. A decision maker offers an {\em experimental assortment} to observe customer choices, then uses these observations to select a {\em final assortment}. 
Our goal is to design an experiment that helps the decision maker select a final assortment with high expected revenue.

\subsection{The Choice Model}\label{sec:model:choice}
Consider a universe of $n$ substitutable products, indexed by $[n]:=\{1,\dots,n\}$. Each customer purchases at most one product
from an offered set of products according to the multinomial logit (MNL) model. In this model, product $i$ has a
preference weight $v_i\ge0$, and the no-purchase option, indexed by $0$, has weight $v_0:=1$.
We allow $v_i=0$, in which case product $i$ is never purchased. For an assortment $S\subseteq[n]$, let
$S^+:=S\cup\{0\}$ and $V(S):=\sum_{i\in S}v_i$.  A customer offered $S$ selects $i\in S^+$ with probability
\[
\phi_v(i, S) \;:=\; \frac{v_i}{1 + V(S)}\textcolor{blue}{.}
\]
This defines a probability distribution over the products in $S$ and the no-purchase option. Each product $i$ has a known revenue $r_i$. We assume
that revenues are positive and distinct; without loss of generality, $r_1>r_2>\dots>r_n$. The expected revenue of $S$
is
\[
R_v(S) \;=\; \sum_{i \in S} r_i \,\phi_v(i, S) \;=\; \frac{\sum_{i \in S} r_i v_i}{1 + V(S)}.
\]
The optimal expected revenue under $v$ is $R_v^{{\star}}:=\max_{S\subseteq[n]}R_v(S)$. 
\citet{talluri2004} show that, under the MNL model, an optimal assortment is a revenue prefix, i.e., 
the set $[k]:=\{1,\dots,k\}$ of the $k$ highest-revenue products for some $k\in[n]$.
Throughout the paper, we use the convention that a maximum or minimum over an empty set equals zero.

\subsection{Experiment Design}\label{sec:model:design}

The revenue vector $r=(r_1,\dots,r_n)$ is fixed and known to the decision maker, while the preference vector $v=(v_1,\dots,v_n)$ is fixed but unknown.
We consider a setting where the decision maker proceeds in two phases: an experimentation phase, followed by a decision phase.
In the experimentation phase, the decision maker chooses an \emph{experimental assortment} $S\subseteq[n]$ to offer to $m$ customers
and observes their independent choices $O_t\sim\phi_v(\cdot,S)$, where $O_t\in S^+$ for $t=1,\dots,m$.
In the decision phase, the decision maker uses the observations $O=(O_1,\dots,O_m)$ to select a \emph{final assortment}
$\mathcal A(O)$, where the \emph{decision rule} $\mathcal A:([n]^+)^m\to2^{[n]}$ maps observations to an assortment.
Together, the experimental assortment and the decision rule constitute a \emph{policy} $(S,\mathcal A)$. We  restrict attention to deterministic experimental designs and deterministic decision rules; randomness arises only from customer choices.
Both $S$ and $\mathcal A$ may depend on the known revenue vector $r$; we suppress this dependence in the notation.
We also study the case $m=0$, in which the decision maker must select a final assortment without observing any customer choices.
Our goal is to find a policy that minimizes the worst-case regret, which we define next.

\subsection{Regret}\label{sec:model:gap}

We compare a policy's expected revenue with that of the optimal assortment in hindsight, that is, the assortment that maximizes expected revenue knowing the true preference vector $v$. For
a fixed $v$, the \emph{regret} of a policy $(S,\mathcal A)$ is
\[
\Reg_v(S,\mathcal A)\ :=\ {\max_{C\subseteq[n]}R_v(C)}-\mathbb E_{O{\sim\phi_v(\cdot,S)^{\otimes m}}}\big[R_v(\mathcal A(O))\big],
\]
where the expectation is taken over the observation vector $O$, where $O_t\sim\phi_v(\cdot,S)$ independently for $t=1,\dots,m$. We evaluate a policy by its \emph{worst-case regret},
\begin{equation}\label{eq:model:worst-gap}
\Reg(S,\mathcal A):=\sup_{v\in\mathbb R^n_+}\Reg_v(S,\mathcal A).
\end{equation}

To evaluate the fundamental value of a given experimental assortment, we consider its worst-case regret under the best possible decision rule using its observations. We then optimize over experimental assortments:
\[
\Reg_m(S):=\min_{\mathcal A}\Reg(S,\mathcal A),\qquad \Reg_m^{{\star}}:=\min_{S\subseteq[n]}\Reg_m(S),
\]
where the minimum over $\mathcal A$ is over all decision rules. We make the dependence on $m$ explicit, as these regrets depend on the number of customers involved in the experiment. We call $\Reg_m^{{\star}}$ the \emph{minimax regret} with $m$
observations; a policy is \emph{minimax optimal} if its worst-case regret as defined in~\eqref{eq:model:worst-gap} equals $\Reg_m^{{\star}}$.
Importantly, our criterion measures the value of an experiment by how much it improves the final assortment decision, rather than by how accurately it estimates the underlying choice model.

For $m=0$, these definitions reduce to $\Reg_0^{{\star}}=\min_{C\subseteq[n]}\sup_{v\in\mathbb R^n_+}\bigl(R_v^{{\star}}-R_v(C)\bigr)$,
the minimax regret without an experiment. Because a decision rule can ignore observations, $\Reg_m^{{\star}}$ cannot increase
with $m$: $\Reg_0^{{\star}}\ge \Reg_1^{{\star}}\ge \Reg_2^{{\star}}\ge\cdots$. Hence, $\Reg_0^{{\star}}$ is an upper bound for every $\Reg_m^{{\star}}$, and we use it as
a benchmark to measure how much an experiment reduces the minimax regret.

In Section~\ref{sec:limit}, we will extend these definitions to study the setting of multiple experiments run in parallel.

\section{Assortment Optimization with No Observation}\label{sec:G0}

In this section, we study the case $m=0$, in which the decision maker
selects a fixed assortment $C\subseteq[n]$ using only the known revenue
vector $r$, without observing any customer choices. The goal is to
choose $C$ to minimize regret in the worst case over all preference
vectors $v\in\mathbb R_+^n$. We characterize the minimax regret and
an optimal assortment by showing that the worst-case regret of $C$
is determined by two quantities:

{%
\[
\begin{aligned}
\text{\textit{Exclusion cost:}}\quad &\alpha(C):={\max\bigl(\{r_i:i\notin C\}\cup\{0\}\bigr)},\\[4pt]
\text{\textit{Inclusion cost:}}\quad &\beta(C):={\max\bigl(\{r_1-r_i:i\in C\}\cup\{0\}\bigr)}.
\end{aligned}
\]
}%
The exclusion cost $\alpha(C)$ is the largest revenue among the products
excluded from $C$, or zero if $C=[n]$. The inclusion cost $\beta(C)$ is
the difference between $r_1$ and the lowest revenue among the products
included in $C$, or zero if $C=\emptyset$. 
The following theorem characterizes the worst-case regret of any
assortment and shows that a revenue prefix is minimax optimal.

\begin{theorem}[{Minimax Regret with No Observation}]\label{thm:G0}
For every $C\subseteq[n]$,
\[
\sup_{v\in\mathbb R_+^n}\bigl(R_v^{{\star}}-R_v(C)\bigr)
=
\max\{\alpha(C),\beta(C)\}.
\]
With the convention $r_{n+1}:=0$,
\[
\Reg_0^{{\star}}
=
\min_{k\in[n]}\max\{r_{k+1},r_1-r_k\}.
\]
For every index $k$ attaining this minimum, the prefix $[k]$ is an optimal policy.
\end{theorem}

Theorem~\ref{thm:G0} shows that the worst-case regret of an assortment $C$ is the larger of its exclusion and inclusion costs. The exclusion cost captures the potential regret from a lost sale of a desirable product, while the inclusion cost captures the potential regret from including a low-revenue product that draws purchases away from a higher-revenue alternative.
This characterization shows that an optimal assortment can be chosen as a revenue prefix. Indeed, for any nonempty assortment $C$, adding all excluded products with revenues higher than the lowest revenue in $C$ does not increase $\alpha(C)$ and leaves $\beta(C)$ unchanged. For a prefix assortment $[k]$, the exclusion and inclusion costs are $r_{k+1}$ and $r_1-r_k$, respectively, so its worst-case regret is $\max\{r_{k+1},r_1-r_k\}$. The proof of Theorem~\ref{thm:G0} is deferred to Appendix~\ref{app:G0}.

The following corollary gives a simpler characterization of an optimal
policy through a threshold on product revenues.

{%
\begin{corollary}[{Optimal Policy}]\label{cor:G0half}
Let $k_0^{{\star}}:=\max\{i\in[n]:r_i>r_1/2\}$ be the largest index whose revenue exceeds $r_1/2$.
With no observation, the prefix $[k_0^{{\star}}]$ is an optimal policy, and $\Reg_0^{{\star}}\le r_1/2$.
\end{corollary}
}

We defer the proof of Corollary~\ref{cor:G0half} to Appendix~\ref{app:G0}.
The minimax regret and an optimal policy admit an intuitive geometric interpretation. Place the products' revenues on the real line and include zero. A prefix assortment $[k]$ partitions these points into two groups: the revenues of the excluded products together with zero, and those of the included products. The diameters of these groups are exactly the exclusion cost $r_{k+1}$ and the inclusion cost $r_1-r_k$, respectively. Thus, \emph{the minimax problem is to choose a cut} between consecutive revenue levels that minimizes the larger of the two group diameters. Corollary~\ref{cor:G0half} shows that a cut at $r_1/2$ is optimal. Figure~\ref{fig:G0-partition} illustrates this interpretation. This observation will help us contrast the no-observation benchmark with the results obtained through experimentation, which admit a closely related geometric interpretation.

\begingroup
\setlength{\intextsep}{8pt}
\begin{figure}[H]
\centering
{%
\begin{tikzpicture}[x=1cm,y=1cm,draw=black!75,text=black!85,>=Stealth]
  \fill[orange!12] (-0.25,-0.22) rectangle (4.05,0.22);
  \fill[teal!14] (6.75,-0.22) rectangle (9.65,0.22);
  \draw[->] (-0.4,0)--(10.1,0) node[right] {Revenue};
  \foreach \x in {0,0.65,0.95,3.05,3.8}
    \fill[orange!80!black] (\x,0) circle (2pt);
  \foreach \x in {7.0,7.3,8.85,9.4}
    \fill[teal!80!black] (\x,0) circle (2pt);
  \node[below=7pt] at (0,0) {$0$};
  \node[below=7pt] at (0.65,0) {$r_n$};
  \node[below=7pt] at (3.8,0) {$r_{k+1}$};
  \node[below=7pt] at (7.0,0) {$r_k$};
  \node[below=7pt] at (9.4,0) {$r_1$};
  \draw[dashed,black!55] (4.7,-0.5)--(4.7,0.7);
  \node[above,font=\small] at (4.7,0.7) {Cut at $r_1/2$};
  \draw[<->,orange!80!black] (0,0.55)--(3.8,0.55)
    node[midway,above,font=\small,text=black] {$\alpha([k])=r_{k+1}$};
  \draw[<->,teal!80!black] (7.0,0.55)--(9.4,0.55)
    node[midway,above,font=\small,text=black] {$\beta([k])=r_1-r_k$};
  \node[font=\small,align=center,text=black] at (1.9,-1.2) {excluded products and zero};
  \node[font=\small,align=center,text=black] at (8.2,-1.2) {included products $[k]$};
\end{tikzpicture}
}
\caption{A revenue prefix partitions the revenue levels into two groups. Its worst-case regret is the larger group diameter; an optimal prefix minimizes this maximum.}
\label{fig:G0-partition}
\end{figure}
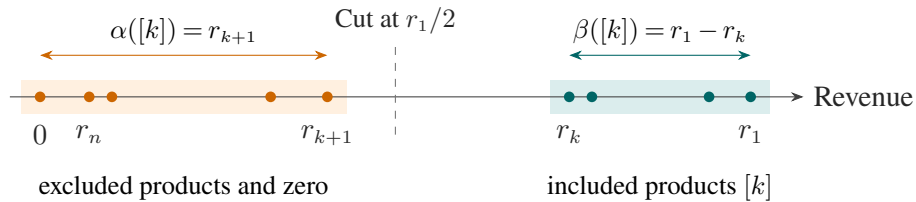
\endgroup

\noindent{\bf Examples.} We illustrate the minimax regret and optimal policies with four revenue structures. In Sections~\ref{sec:G1} and~\ref{sec:rep}, we use the same examples to compare regret reductions from a single observation and from multiple observations relative to the no-observation benchmark.

We consider the following four revenue profiles. Throughout this
example, assume that $n\ge20$ is even. Write $R:=r_1$ and set
$\epsilon_n:=R/n$. Revenues in each specified interval are distinct
and include both endpoints.
\begin{enumerate}
\item[(a)] \textbf{Uniformly spaced revenues:} $r_i=R-(i-1)\epsilon_n$ for $i\in[n]$.
\item[(b)] \textbf{One product at $R$; others near $R/2$:} One product has revenue $R$; the other $n-1$ have revenues in $[R/2-\epsilon_n,R/2]$.
\item[(c)] \textbf{Clusters near $R$ and $R/3$:} Half of the products have revenues in $[R-\epsilon_n,R]$ and half in $[R/3-\epsilon_n/2,R/3+\epsilon_n/2]$.
\item[(d)] \textbf{Clusters near $R$ and $R/8$:} Half of the products have revenues in $[R-\epsilon_n,R]$ and half in $[R/8-\epsilon_n/2,R/8+\epsilon_n/2]$.
\end{enumerate}

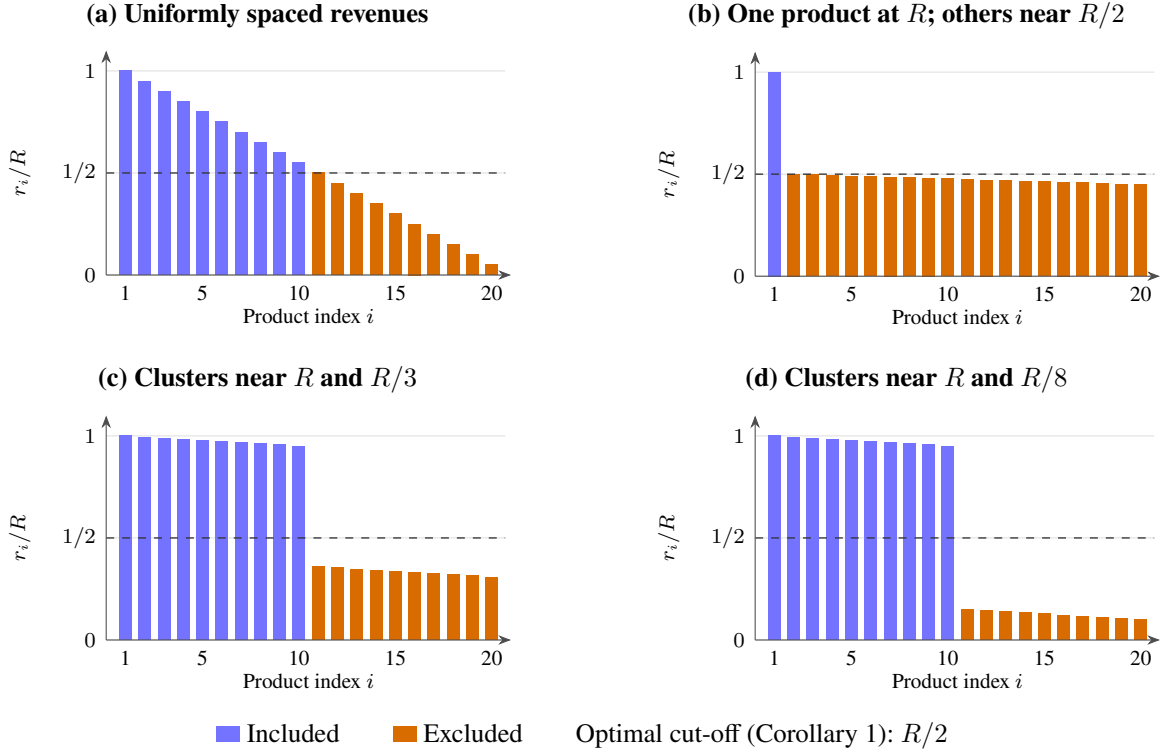
\begin{figure}[tbp]
\centering
\begin{minipage}[t]{0.48\linewidth}
\centering
{\small\textbf{(a) Uniformly spaced revenues}\par}
\vspace{2mm}
\begin{tikzpicture}[x=0.255cm,y=2.7cm,>=Stealth,draw=black!70,text=black,font=\scriptsize]
\path[use as bounding box] (-5.4,-0.29) rectangle (21.1,1.14);
\draw[black!12] (0,0)--(20.7,0);
\node[left] at (0,0) {$0$};
\draw[black!12] (0,0.5)--(20.7,0.5);
\node[left] at (0,0.5) {$1/2$};
\draw[black!12] (0,1)--(20.7,1);
\node[left] at (0,1) {$1$};
\fill[blue!55] (0.67,0) rectangle (1.33,1.00000000);
\fill[blue!55] (1.67,0) rectangle (2.33,0.95000000);
\fill[blue!55] (2.67,0) rectangle (3.33,0.90000000);
\fill[blue!55] (3.67,0) rectangle (4.33,0.85000000);
\fill[blue!55] (4.67,0) rectangle (5.33,0.80000000);
\fill[blue!55] (5.67,0) rectangle (6.33,0.75000000);
\fill[blue!55] (6.67,0) rectangle (7.33,0.70000000);
\fill[blue!55] (7.67,0) rectangle (8.33,0.65000000);
\fill[blue!55] (8.67,0) rectangle (9.33,0.60000000);
\fill[blue!55] (9.67,0) rectangle (10.33,0.55000000);
\fill[orange!85!black] (10.67,0) rectangle (11.33,0.50000000);
\fill[orange!85!black] (11.67,0) rectangle (12.33,0.45000000);
\fill[orange!85!black] (12.67,0) rectangle (13.33,0.40000000);
\fill[orange!85!black] (13.67,0) rectangle (14.33,0.35000000);
\fill[orange!85!black] (14.67,0) rectangle (15.33,0.30000000);
\fill[orange!85!black] (15.67,0) rectangle (16.33,0.25000000);
\fill[orange!85!black] (16.67,0) rectangle (17.33,0.20000000);
\fill[orange!85!black] (17.67,0) rectangle (18.33,0.15000000);
\fill[orange!85!black] (18.67,0) rectangle (19.33,0.10000000);
\fill[orange!85!black] (19.67,0) rectangle (20.33,0.05000000);
\draw[dashed,black!80,line width=0.5pt] (0,0.5)--(20.7,0.5);
\draw[->] (0,0)--(21,0);
\draw[->] (0,0)--(0,1.1);
\node[rotate=90,anchor=south] at (-3.4,0.5) {$r_i/R$};
\node[below] at (1,0) {1};
\node[below] at (5,0) {5};
\node[below] at (10,0) {10};
\node[below] at (15,0) {15};
\node[below] at (20,0) {20};
\node at (10.5,-0.20) {Product index $i$};
\end{tikzpicture}
\end{minipage}\hfill
\begin{minipage}[t]{0.48\linewidth}
\centering
{\small\textbf{(b) One product at $R$; others near $R/2$}\par}
\vspace{2mm}
\begin{tikzpicture}[x=0.255cm,y=2.7cm,>=Stealth,draw=black!70,text=black,font=\scriptsize]
\path[use as bounding box] (-5.4,-0.29) rectangle (21.1,1.14);
\draw[black!12] (0,0)--(20.7,0);
\node[left] at (0,0) {$0$};
\draw[black!12] (0,0.5)--(20.7,0.5);
\node[left] at (0,0.5) {$1/2$};
\draw[black!12] (0,1)--(20.7,1);
\node[left] at (0,1) {$1$};
\fill[blue!55] (0.67,0) rectangle (1.33,1.00000000);
\fill[orange!85!black] (1.67,0) rectangle (2.33,0.50000000);
\fill[orange!85!black] (2.67,0) rectangle (3.33,0.49722222);
\fill[orange!85!black] (3.67,0) rectangle (4.33,0.49444444);
\fill[orange!85!black] (4.67,0) rectangle (5.33,0.49166667);
\fill[orange!85!black] (5.67,0) rectangle (6.33,0.48888889);
\fill[orange!85!black] (6.67,0) rectangle (7.33,0.48611111);
\fill[orange!85!black] (7.67,0) rectangle (8.33,0.48333333);
\fill[orange!85!black] (8.67,0) rectangle (9.33,0.48055556);
\fill[orange!85!black] (9.67,0) rectangle (10.33,0.47777778);
\fill[orange!85!black] (10.67,0) rectangle (11.33,0.47500000);
\fill[orange!85!black] (11.67,0) rectangle (12.33,0.47222222);
\fill[orange!85!black] (12.67,0) rectangle (13.33,0.46944444);
\fill[orange!85!black] (13.67,0) rectangle (14.33,0.46666667);
\fill[orange!85!black] (14.67,0) rectangle (15.33,0.46388889);
\fill[orange!85!black] (15.67,0) rectangle (16.33,0.46111111);
\fill[orange!85!black] (16.67,0) rectangle (17.33,0.45833333);
\fill[orange!85!black] (17.67,0) rectangle (18.33,0.45555556);
\fill[orange!85!black] (18.67,0) rectangle (19.33,0.45277778);
\fill[orange!85!black] (19.67,0) rectangle (20.33,0.45000000);
\draw[dashed,black!80,line width=0.5pt] (0,0.5)--(20.7,0.5);
\draw[->] (0,0)--(21,0);
\draw[->] (0,0)--(0,1.1);
\node[rotate=90,anchor=south] at (-3.4,0.5) {$r_i/R$};
\node[below] at (1,0) {1};
\node[below] at (5,0) {5};
\node[below] at (10,0) {10};
\node[below] at (15,0) {15};
\node[below] at (20,0) {20};
\node at (10.5,-0.20) {Product index $i$};
\end{tikzpicture}
\end{minipage}
\par\vspace{4mm}
\begin{minipage}[t]{0.48\linewidth}
\centering
{\small\textbf{(c) Clusters near $R$ and $R/3$}\par}
\vspace{2mm}
\begin{tikzpicture}[x=0.255cm,y=2.7cm,>=Stealth,draw=black!70,text=black,font=\scriptsize]
\path[use as bounding box] (-5.4,-0.29) rectangle (21.1,1.14);
\draw[black!12] (0,0)--(20.7,0);
\node[left] at (0,0) {$0$};
\draw[black!12] (0,0.5)--(20.7,0.5);
\node[left] at (0,0.5) {$1/2$};
\draw[black!12] (0,1)--(20.7,1);
\node[left] at (0,1) {$1$};
\fill[blue!55] (0.67,0) rectangle (1.33,1.00000000);
\fill[blue!55] (1.67,0) rectangle (2.33,0.99444444);
\fill[blue!55] (2.67,0) rectangle (3.33,0.98888889);
\fill[blue!55] (3.67,0) rectangle (4.33,0.98333333);
\fill[blue!55] (4.67,0) rectangle (5.33,0.97777778);
\fill[blue!55] (5.67,0) rectangle (6.33,0.97222222);
\fill[blue!55] (6.67,0) rectangle (7.33,0.96666667);
\fill[blue!55] (7.67,0) rectangle (8.33,0.96111111);
\fill[blue!55] (8.67,0) rectangle (9.33,0.95555556);
\fill[blue!55] (9.67,0) rectangle (10.33,0.95000000);
\fill[orange!85!black] (10.67,0) rectangle (11.33,0.35833333);
\fill[orange!85!black] (11.67,0) rectangle (12.33,0.35277778);
\fill[orange!85!black] (12.67,0) rectangle (13.33,0.34722222);
\fill[orange!85!black] (13.67,0) rectangle (14.33,0.34166667);
\fill[orange!85!black] (14.67,0) rectangle (15.33,0.33611111);
\fill[orange!85!black] (15.67,0) rectangle (16.33,0.33055556);
\fill[orange!85!black] (16.67,0) rectangle (17.33,0.32500000);
\fill[orange!85!black] (17.67,0) rectangle (18.33,0.31944444);
\fill[orange!85!black] (18.67,0) rectangle (19.33,0.31388889);
\fill[orange!85!black] (19.67,0) rectangle (20.33,0.30833333);
\draw[dashed,black!80,line width=0.5pt] (0,0.5)--(20.7,0.5);
\draw[->] (0,0)--(21,0);
\draw[->] (0,0)--(0,1.1);
\node[rotate=90,anchor=south] at (-3.4,0.5) {$r_i/R$};
\node[below] at (1,0) {1};
\node[below] at (5,0) {5};
\node[below] at (10,0) {10};
\node[below] at (15,0) {15};
\node[below] at (20,0) {20};
\node at (10.5,-0.20) {Product index $i$};
\end{tikzpicture}
\end{minipage}\hfill
\begin{minipage}[t]{0.48\linewidth}
\centering
{\small\textbf{(d) Clusters near $R$ and $R/8$}\par}
\vspace{2mm}
\begin{tikzpicture}[x=0.255cm,y=2.7cm,>=Stealth,draw=black!70,text=black,font=\scriptsize]
\path[use as bounding box] (-5.4,-0.29) rectangle (21.1,1.14);
\draw[black!12] (0,0)--(20.7,0);
\node[left] at (0,0) {$0$};
\draw[black!12] (0,0.5)--(20.7,0.5);
\node[left] at (0,0.5) {$1/2$};
\draw[black!12] (0,1)--(20.7,1);
\node[left] at (0,1) {$1$};
\fill[blue!55] (0.67,0) rectangle (1.33,1.00000000);
\fill[blue!55] (1.67,0) rectangle (2.33,0.99444444);
\fill[blue!55] (2.67,0) rectangle (3.33,0.98888889);
\fill[blue!55] (3.67,0) rectangle (4.33,0.98333333);
\fill[blue!55] (4.67,0) rectangle (5.33,0.97777778);
\fill[blue!55] (5.67,0) rectangle (6.33,0.97222222);
\fill[blue!55] (6.67,0) rectangle (7.33,0.96666667);
\fill[blue!55] (7.67,0) rectangle (8.33,0.96111111);
\fill[blue!55] (8.67,0) rectangle (9.33,0.95555556);
\fill[blue!55] (9.67,0) rectangle (10.33,0.95000000);
\fill[orange!85!black] (10.67,0) rectangle (11.33,0.15000000);
\fill[orange!85!black] (11.67,0) rectangle (12.33,0.14444444);
\fill[orange!85!black] (12.67,0) rectangle (13.33,0.13888889);
\fill[orange!85!black] (13.67,0) rectangle (14.33,0.13333333);
\fill[orange!85!black] (14.67,0) rectangle (15.33,0.12777778);
\fill[orange!85!black] (15.67,0) rectangle (16.33,0.12222222);
\fill[orange!85!black] (16.67,0) rectangle (17.33,0.11666667);
\fill[orange!85!black] (17.67,0) rectangle (18.33,0.11111111);
\fill[orange!85!black] (18.67,0) rectangle (19.33,0.10555556);
\fill[orange!85!black] (19.67,0) rectangle (20.33,0.10000000);
\draw[dashed,black!80,line width=0.5pt] (0,0.5)--(20.7,0.5);
\draw[->] (0,0)--(21,0);
\draw[->] (0,0)--(0,1.1);
\node[rotate=90,anchor=south] at (-3.4,0.5) {$r_i/R$};
\node[below] at (1,0) {1};
\node[below] at (5,0) {5};
\node[below] at (10,0) {10};
\node[below] at (15,0) {15};
\node[below] at (20,0) {20};
\node at (10.5,-0.20) {Product index $i$};
\end{tikzpicture}
\end{minipage}
\par\vspace{3mm}
{\small\textcolor{blue!55}{\rule{9pt}{6pt}} Included\qquad
\textcolor{orange!85!black}{\rule{9pt}{6pt}} Excluded\qquad Optimal cut-off (Corollary~\ref{cor:G0half}): $R/2$}
\caption{Four illustrative revenue profiles with $n=20$ and $R=1$. A minimax policy includes exactly the products above $R/2$. Within each cluster, the plotted revenues are equally spaced.}
\label{fig:G0-examples}
\end{figure}

Figure~\ref{fig:G0-examples} illustrates the examples together with their minimax policies, and Table~\ref{tab:examples-G0} reports their optimal policies and minimax regret $\Reg_0^\star$.

\begin{table}[tbp]
\centering
\renewcommand{\arraystretch}{1.4}
\begin{tabular}{clcc}
\hline
& Example & An optimal policy & $\Reg_0^{{\star}}$ \\
\hline
(a) & Uniformly spaced revenues & $[\frac{n}{2}]$ & $\frac{R}{2}$ \\
(b) & One product at $R$; others near $R/2$ & $[1]$ & $\frac{R}{2}$ \\
(c) & Clusters near $R$ and $R/3$ & $[\frac{n}{2}]$ & $\frac{R}{3}+\frac{\epsilon_n}{2}$ \\
(d) & Clusters near $R$ and $R/8$ & $[\frac{n}{2}]$ & $\frac{R}{8}+\frac{\epsilon_n}{2}$ \\
\hline
\end{tabular}
\caption{Optimal policies and minimax regret with no observation for examples (a)--(d), with even $n\ge20$ and $\epsilon_n=R/n$.}
\label{tab:examples-G0}
\end{table}

{%

These examples show how the minimax-optimal policy and regret depend on the revenue structure. In examples (c) and (d), an optimal policy offers the expensive cluster, while in example (b), it offers only the highest-revenue product. Examples (b)--(d) illustrate that a larger revenue gap between the expensive and cheap products lowers the minimax regret. Examples (a) and (b) attain $R/2$, showing that the upper bound in Corollary~\ref{cor:G0half} is tight. At the other extreme, the regret can be arbitrarily close to zero: with two products of revenues $R$ and $\delta R$, where $0<\delta<1/2$, offering only the first product gives minimax regret $\delta R$, which tends to zero as $\delta\to0$. We derive the values in Table~\ref{tab:examples-G0} in Appendix~\ref{app:running-examples}.

}
\FloatBarrier

\section{Assortment Optimization with a Single Observation}\label{sec:G1}

{In this section, we study the setting where the decision maker first chooses an experimental assortment and observes a single customer's purchase realization, including the possibility of no purchase. Based on this observation, the decision maker then commits to a final assortment. Our goal is to identify the minimax-optimal experimental assortment and its associated optimal decision rule, and to quantify the value of this single experiment: can it meaningfully reduce the decision maker's regret relative to the no-observation benchmark?}

{The following theorem characterizes the minimax regret and an optimal policy with a single observation.}

\begin{theorem}[{{Minimax Optimal} Policy with a Single Observation}]\label{thm:one:opt}
With the convention $r_{n+1}:=0$, the minimax regret with one observation satisfies
{
\[
\Reg_1^\star
=\min_{1\le k\le p\le\ell\le n}
\max\left\{
 r_{\ell+1},\ r_1-r_p,\ r_{p+1}-r_\ell+\frac{r_\ell^2}{4r_k},\
 r_{k+1}-r_\ell,\ \frac{(r_1-r_\ell)^2}{4r_1}
\right\}.
\]
}%
{A triple $(k^{{\star}},p^{{\star}},\ell^{{\star}})$ that attains this minimum specifies a minimax-optimal policy as follows: the decision maker offers the experimental assortment $[k^{{\star}}]$. If the customer purchases an item,  the decision maker commits
to the final assortment $[p^\star]$. If the customer makes no
purchase,  the decision maker commits to $[\ell^\star]$.
}
\end{theorem}

We defer the proof of Theorem~\ref{thm:one:opt} to
Appendix~\ref{app:G1:lb}. We next provide intuition for the structure
of the optimal policy and interpret the terms in the minimax regret
expression. Note that $k^\star \le p^\star \le \ell^\star$. The optimal policy introduced in this theorem offers
the revenue prefix $[k^\star]$ during the experiment. For the final
assortment, it may expand to $[p^\star]$ after a purchase or further
to $[\ell^\star]$ after no purchase.

\vspace{2mm}

\noindent\textbf{Why the policy ignores product identity.} At first glance, it may seem surprising that the final assortment depends only on whether a purchase occurs, rather than on which product is chosen. One might instead try to use a maximum-likelihood estimator of the preference weights to select the final assortment. Given the experimental assortment $[k^\star]$ and a purchase of $i\in[k^\star]$, the likelihood approaches its supremum as the estimated weight of $i$ grows without bound and the other experimental weights are set to zero. This suggests choosing $[i]$ as the final assortment to maximize the estimated expected revenue. However, observing $i$ does not establish that its preference weight is much higher than those of the other experimental products. They could actually have equal preference weights, while a product held out of the experiment has a much larger weight without affecting the experimental observations. When revenues are sufficiently close, including such a product can increase the purchase probability with little reduction in revenue per purchase, giving higher expected revenue than $[i]$. Thus, selecting $[i]$ based on this fit can exclude products whose
inclusion would increase expected revenue under the true preference vector.

\noindent\textbf{Interpreting the regret terms.}
For intuition, consider a triple $(k,p,\ell)$ with $k\le p<\ell$.
The three linear terms reflect inclusion and exclusion losses,
as in Section~\ref{sec:G0}. The term $r_{\ell+1}$ captures the loss from
excluding a desirable product from both final assortments.
The term $r_1-r_p$ captures the loss when product $p$, included
after a purchase, draws demand away from the higher-revenue
product~$1$. The term $r_{k+1}-r_\ell$ captures the potential loss
from offering unobserved products together after no purchase.
Their relative preference weights remain unknown, and a
lower-revenue product can dominate purchases even though a
higher-revenue held-out product would yield greater expected
revenue on its own. The optimal policy balances these inclusion
and exclusion losses when choosing the experimental assortment
and the two final assortments.

A single observation introduces an additional source of loss.
A customer may make no purchase even when the experimental
products attract strong demand. This outcome leads the policy
to expand to $[\ell]$, although offering $[p]$ may yield higher
expected revenue. Conversely, a customer may purchase even when
experimental demand is weak, leading the policy to select $[p]$
although expanding to $[\ell]$ would yield higher expected revenue.
In either case, the additional loss relative to the better of
the two final assortments is the difference between their expected
revenues multiplied by the probability of selecting the one with lower
expected revenue. This loss is in addition to any regret that
remains because even the better final assortment may be
suboptimal. Both the observation probability and the revenue
difference depend on the preference weights. The two remaining terms, which contain quadratic contributions,
account for these effects, with their precise expressions
derived in Appendix~\ref{app:G1:lb}.

\noindent\textbf{A tight upper bound.}
While Theorem~\ref{thm:one:opt} gives an exact characterization
of the minimax regret, evaluating it requires optimizing over
three cutoff indices. The following lemma provides a simple
upper bound that depends only on the highest product revenue
$r_1$. Importantly, a single observation improves the tight
universal bound from $r_1/2$, established in Corollary~\ref{cor:G0half} for
the no-observation setting, to $8r_1/23$.

\begin{lemma} \label{lem:one:823}
For every $n\ge1$ and all revenues $r_1>\cdots>r_n>0$,
\[
    \operatorname{Reg}_1^\star \le \frac{8}{23}r_1.
\]
This bound is tight.
\end{lemma}
Tightness is attained already with four products: for the
revenue vector $(23,16,15,8)$, the minimax regret equals
$8=(8/23)r_1$. Thus, the coefficient $8/23$ cannot be reduced
in a bound valid for all revenue vectors. We defer the proof to Appendix~\ref{app:G1:823}.

A simple threshold-based policy guarantees worst-case regret
at most $8r_1/23$ without solving the minimization in
Theorem~\ref{thm:one:opt}. Define
\[
\begin{aligned}
\ell &:= \max\left\{i\in[n]:r_i>\frac{8}{23}r_1\right\}, \qquad p := \max\left\{i\in[n]:r_i\ge\frac{15}{23}r_1\right\},\\
k &:= \max\left\{i\in[n]:
r_i\ge\min\left\{r_1,r_\ell+\frac{8}{23}r_1\right\}\right\}.
\end{aligned}
\]
The policy experiments with $[k]$, commits to $[p]$ after
a purchase, and commits to $[\ell]$ after no purchase.
These thresholds satisfy $k\le p\le\ell$, and substituting
them into the expression in Theorem~\ref{thm:one:opt}
bounds each term by $8r_1/23$.

\noindent\textbf{Examples.}
We return to the four revenue profiles in
Examples in Section \ref{sec:G0}, to illustrate how the revenue
structure affects the value of a single observation.
Table~\ref{tab:examples-G1} summarizes the limiting minimax
regret and the relative regret reduction
$\frac{\Reg_0^\star-\Reg_1^\star}{\Reg_0^\star}$.

In example (a), with uniformly spaced revenues, the minimax
regret approaches the upper bound $8R/23$ as $n\to\infty$.
The relative regret reduction nevertheless approaches
$7/23\approx30.4\%$, illustrating the improvement over the
no-observation benchmark even when the one-observation bound
is asymptotically attained. The benefits can be substantially
larger for other revenue structures. In example (b), with one
product at $R$ and the others near $R/2$, the relative regret
reduction approaches $87.5\%$. In example (c), with clusters
near $R$ and $R/3$, it approaches approximately $66.6\%$.
Thus, even a single observation can substantially reduce
worst-case regret. In contrast, in example (d), with clusters
near $R$ and $R/8$, the reduction is exactly zero for every
even $n\ge20$. In this case, a single observation has no value
for reducing minimax regret relative to the no-observation
benchmark. Table~\ref{tab:examples-G1-finite} in
Appendix~\ref{app:running-examples} gives the finite-$n$ values.

\begin{table}[t]

\centering
\renewcommand{\arraystretch}{2}
\normalsize
\setlength{\tabcolsep}{4pt}
\begin{tabular}{@{}lccc|cl@{}}
\hline
Example & $\Reg_0^\star$ & $\Reg_1^\star$ & \makecell{Regret\\reduction} & Experiment & Decision rule \\
\hline
\makecell[l]{(a) Uniformly spaced\\\phantom{(a) }revenues} & $\mathchoice{\hbox{\small$\displaystyle\frac{R}{2}$}}{\hbox{\small$\displaystyle\frac{R}{2}$}}{\frac{R}{2}}{\frac{R}{2}}$ & $\mathchoice{\hbox{\small$\displaystyle\frac{8R}{23}$}}{\hbox{\small$\displaystyle\frac{8R}{23}$}}{\frac{8R}{23}}{\frac{8R}{23}}$ & $\frac{7}{23}\approx30.4\%$ & $\approx[7n/23]$ & \makecell[l]{Purchase: $\approx[8n/23]$\\No purchase: $\approx[15n/23]$} \\[2pt]
\makecell[l]{(b) One product at $R$;\\\phantom{(b) }others near $R/2$} & $\mathchoice{\hbox{\small$\displaystyle\frac{R}{2}$}}{\hbox{\small$\displaystyle\frac{R}{2}$}}{\frac{R}{2}}{\frac{R}{2}}$ & $\mathchoice{\hbox{\small$\displaystyle\frac{R}{16}$}}{\hbox{\small$\displaystyle\frac{R}{16}$}}{\frac{R}{16}}{\frac{R}{16}}$ & $\frac{7}{8}=87.5\%$ & $[1]$ & \makecell[l]{Purchase: $[1]$\\No purchase: $[n]$} \\[2pt]
\makecell[l]{(c) Clusters near\\\phantom{(c) }$R$ and $R/3$} & $\mathchoice{\hbox{\small$\displaystyle\frac{R}{3}$}}{\hbox{\small$\displaystyle\frac{R}{3}$}}{\frac{R}{3}}{\frac{R}{3}}$ & $\mathchoice{\hbox{\small$\displaystyle\frac{R}{9}$}}{\hbox{\small$\displaystyle\frac{R}{9}$}}{\frac{R}{9}}{\frac{R}{9}}$ & $\frac{2}{3}\approx66.6\%$ & $[n/2]$ & \makecell[l]{Purchase: $[n/2]$\\No purchase: $[n]$} \\[2pt]
\makecell[l]{(d) Clusters near\\\phantom{(d) }$R$ and $R/8$} & $\mathchoice{\hbox{\small$\displaystyle\frac{R}{8}$}}{\hbox{\small$\displaystyle\frac{R}{8}$}}{\frac{R}{8}}{\frac{R}{8}}$ & $\mathchoice{\hbox{\small$\displaystyle\frac{R}{8}$}}{\hbox{\small$\displaystyle\frac{R}{8}$}}{\frac{R}{8}}{\frac{R}{8}}$ & $0\%$ & $[n/2]$ & \makecell[l]{$[n/2]$ for\\either outcome} \\[2pt]
\hline
\end{tabular}
\caption{One-observation policies, regrets, and relative regret reductions for examples (a)--(d). Displayed policies are optimal for even $n\ge20$, except that the $\approx$ entries in example (a) show large-$n$ approximations; regret values and reductions are limits as $n\to\infty$ through even integers, with $R$ fixed. Regret reductions are relative to the no-observation minimax regret $\Reg_0^\star$. Percentage coefficients are rounded down to one decimal place.}
\label{tab:examples-G1}
\end{table}

Table~\ref{tab:examples-G1} also reveals the structure of the minimax-optimal policies. In examples (b) and (c), the decision maker experiments with expensive products. A purchase indicates that these products have some positive desirability, and retaining them is a minimax-optimal decision. A no-purchase outcome, by contrast, is consistent with these products being potentially undesirable. Offering only them can result in few or no purchases and hence high regret. The minimax-optimal response is to expand to the full assortment. In example (d), however, the cheap products generate so little revenue that a single no-purchase outcome still does not justify adding them. The decision maker instead retains the optimal no-observation assortment from Table~\ref{tab:examples-G0}. In the next section, we will see that further observations can support more elaborate decision rules and substantially greater regret reductions, even in this last example.

\section{Assortment Optimization with Multiple Observations}\label{sec:rep}

In this section, we study the setting where the decision maker
offers a single experimental assortment to $m$ customers and
uses their observed choices to select a final assortment.
We characterize the limiting minimax regret as $m\to\infty$
and identify an experimental assortment together with an
associated decision rule that attains this limit. We also
provide an explicit bound on the policy's regret for every
finite $m$.

\noindent\textbf{Is the full assortment an optimal experiment?} At first glance, it might seem that experimenting with the full assortment should be optimal: for every fixed preference vector $v$, the observed purchase frequencies allow arbitrarily accurate recovery of all preference weights with enough observations. Therefore, a single experimental assortment would identify an optimal final assortment asymptotically, yielding zero limiting regret. This conclusion is misleading for the minimax problem, however. As long as the observation budget $m$ is finite, the preference weights can scale with the budget and make different optimal assortments arbitrarily difficult to distinguish. 
Thus, convergence to zero for each fixed preference vector
does not imply that worst-case regret converges to zero.
For a fixed experimental assortment $S$, let $(\mathcal A_m)_{m\ge1}$ be a sequence of decision rules, where $\mathcal A_m$ uses $m$ observations from $S$. Taking the supremum before the limit can give
a different value from taking the limit first:
\[
\lim_{m\to\infty}\sup_v \Reg_v(S,\mathcal A_m)
\neq
\sup_v\lim_{m\to\infty}\Reg_v(S,\mathcal A_m).
\]

For example, consider two products with revenues $(1,1/2)$ and preference vectors $v^{(a)}=(H^{-1},H^2)$ and $v^{(b)}=(H,H^2)$, where $H>1$. Their unique optimal assortments are $\{1,2\}$ and $\{1\}$, respectively. Yet, when experimenting with the full assortment on any fixed number $m$ of customers, the decision maker observes only purchases of product~$2$ with probability arbitrarily close to one as $H\to\infty$, under either vector. Hence, both preference vectors produce the same observation with high probability, but no single decision rule can achieve zero regret under both, as they have different optimal assortments. Choosing either $\{1\}$ or $\{1,2\}$ therefore gives regret approaching $1/2$ under one of the two preference vectors as $H\to\infty$.

Hence, under any finite observation budget, an experiment must account for preference vectors that can be arbitrarily difficult to distinguish. Our interest is precisely this finite-budget setting: we study $\Reg_m^\star$, which takes the worst case over preferences at each budget, and then its limit as $m\to\infty$.

If $n=1$, committing to the sole product gives zero regret for every $m$. We therefore assume $n\ge2$ throughout this section. In the following theorem, we prove that the limiting minimax regret is
\begin{equation}\label{eq:Ginf-explicit-intro}
\Reg_\infty^\star
:=\min_{{1\le k<\ell\le n}}
\max\{r_1-r_k,\ {r_{k+1}-r_\ell},\ r_{\ell+1}\},
\end{equation}
where $r_{n+1}:=0$.

\begin{theorem}[Asymptotically Minimax Optimal Policy with Multiple Observations] \label{thm:rep:main}
Let $(k^\star,\ell^\star)$ be a pair attaining the minimum in~\eqref{eq:Ginf-explicit-intro}.
For $m\ge1$, the decision maker offers $[k^\star]$ to $m$ customers. Let
$\hat z$ denote the observed fraction of these customers who purchase
a product from $[k^\star]$. Let $\mathcal A_m$ be the decision rule that
selects $[k^\star]$ if ${\hat z\ge r_{\ell^\star}/r_{k^\star}}$,
and $[\ell^\star]$ otherwise. Then
\begin{equation} \label{eq:sandwich}
\Reg_\infty^\star
\le \Reg_m^\star
\le \Reg_m([k^\star])
\le \Reg([k^\star],\mathcal A_m)
\le \Reg_\infty^\star+\frac{r_1}{2\sqrt m}.
\end{equation}
Consequently, both $\Reg_m^\star$ and $\Reg_m([k^\star])$
converge to $\Reg_\infty^\star$ as $m\to\infty$, and the
policy is asymptotically minimax optimal.
\end{theorem}

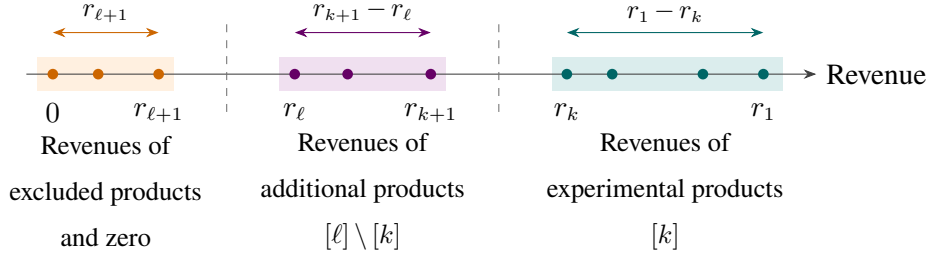
\begin{figure}[t]
\centering
\begin{tikzpicture}[x=1cm,y=1cm,draw=black!75,text=black,>=Stealth]
  \fill[orange!12] (-0.2,-0.22) rectangle (1.6,0.22);
  \fill[violet!12] (3.0,-0.22) rectangle (5.2,0.22);
  \fill[teal!14] (6.6,-0.22) rectangle (9.65,0.22);
  \draw[->] (-0.35,0)--(10.1,0) node[right] {Revenue};
  \foreach \x in {0,0.6,1.4}
    \fill[orange!80!black] (\x,0) circle (2pt);
  \foreach \x in {3.2,3.9,5.0}
    \fill[violet!80!black] (\x,0) circle (2pt);
  \foreach \x in {6.8,7.4,8.6,9.4}
    \fill[teal!80!black] (\x,0) circle (2pt);
  \node[below=7pt] at (0,0) {$0$};
  \node[below=7pt] at (1.4,0) {$r_{\ell+1}$};
  \node[below=7pt] at (3.2,0) {$r_\ell$};
  \node[below=7pt] at (5.0,0) {$r_{k+1}$};
  \node[below=7pt] at (6.8,0) {$r_k$};
  \node[below=7pt] at (9.4,0) {$r_1$};
  \draw[dashed,black!55] (2.3,-0.5)--(2.3,0.7);
  \draw[dashed,black!55] (5.9,-0.5)--(5.9,0.7);
  \draw[<->,orange!80!black] (0,0.55)--(1.4,0.55)
    node[midway,above,font=\small,text=black] {$r_{\ell+1}$};
  \draw[<->,violet!80!black] (3.2,0.55)--(5.0,0.55)
    node[midway,above,font=\small,text=black] {$r_{k+1}-r_\ell$};
  \draw[<->,teal!80!black] (6.8,0.55)--(9.4,0.55)
    node[midway,above,font=\small,text=black] {$r_1-r_k$};
  \node[font=\small,align=center,text=black] at (0.7,-1.55) {Revenues of\\excluded products\\and zero};
  \node[font=\small,align=center,text=black] at (4.1,-1.55) {Revenues of\\additional products\\$[\ell]\setminus[k]$};
  \node[font=\small,align=center,text=black] at (8.1,-1.55) {Revenues of\\experimental products\\$[k]$};
\end{tikzpicture}
\caption{Two cutoffs $k<\ell$ partition the revenues, including zero, into three nonempty consecutive groups. The limiting minimax regret is the smallest possible largest group diameter. At a minimizing pair, the policy retains the experimental group after sufficiently frequent purchases and otherwise adds the middle group.}
\label{fig:Ginf-partition}
\end{figure}

We defer the proof of Theorem~\ref{thm:rep:main}
to Appendix~\ref{sec:rep:proof}. We next interpret the limiting regret and
the structure of the associated policy. As in the
no-observation setting, place the revenues, including zero,
on the real line. The cutoffs $k<\ell$ partition them into
three nonempty consecutive groups: the revenues of the
experimental products $[k]$, those of the additional products
$[\ell]\setminus[k]$, and the remaining revenues together
with zero. Their diameters are $r_1-r_k$,
$r_{k+1}-r_\ell$, and $r_{\ell+1}$, respectively.
Thus, minimizing the largest group diameter gives both
the limiting minimax regret and the cutoffs defining an
asymptotically minimax optimal experiment and decision rule.
Figure~\ref{fig:Ginf-partition} illustrates this interpretation, extending the
two-group representation of the no-observation problem
in Figure \ref{fig:G0-partition}.

\noindent\textbf{Interpreting the regret terms.}
The three terms reflect inclusion and exclusion losses.
The term $r_1-r_{k^\star}$ captures the loss from retaining
all experimental products in both final assortments.
For any finite observation budget, product $k^\star$ can
dominate purchases even though offering product~1 alone
would yield expected revenue close to $r_1$. Retaining
product $k^\star$ can then reduce final expected revenue
to nearly $r_{k^\star}$.

The term $r_{k^\star+1}-r_{\ell^\star}$ captures a similar
loss among the products added when the observed purchase
frequency is below the threshold. Their preference weights
are unobserved, so a lower-revenue product can dominate
purchases even though a higher-revenue product in this group
would yield greater expected revenue on its own.
Finally, $r_{\ell^\star+1}$ captures the loss from excluding
a desirable product from both final assortments. If the
included products attract almost no demand, final expected
revenue can be close to zero, while offering the excluded
product alone can yield expected revenue close to its
revenue. This explains why the last group includes zero;
if no products are excluded, its diameter is zero.

Although the experiment records individual product choices,
this asymptotically minimax optimal decision rule uses only
the overall purchase frequency. Every product in $[k^\star]$
has revenue at least $r_{k^\star}$, so $r_{k^\star}\hat z$
estimates a lower bound on this assortment's expected revenue.
The policy compares this estimate with $r_{\ell^\star}$,
the lowest revenue among the additional products: it retains
$[k^\star]$ when the estimate is at least $r_{\ell^\star}$
and otherwise expands to $[\ell^\star]$. Additional observations
estimate the purchase probability more accurately, and
Theorem~\ref{thm:rep:main} bounds the resulting additional
regret by $r_1/(2\sqrt m)$.

\medskip

\noindent\textbf{A tight upper bound and a simple policy.}
The upper bound in Theorem~\ref{thm:rep:main} approaches
$\Reg_\infty^\star$ at rate $m^{-1/2}$. For small $m$,
ignoring all but the first observation can give a sharper
guarantee. Section~\ref{sec:model:gap} gives
$\Reg_m^\star\le\Reg_1^\star$, and
Lemma~\ref{lem:one:823} gives
$\Reg_1^\star\le 8r_1/23$. Together with
Theorem~\ref{thm:rep:main}, these bounds yield
\[
\Reg_m^\star
\le
\min\left\{
    \Reg_1^\star,\,
    \Reg_\infty^\star+\frac{r_1}{2\sqrt m}
\right\}
\le
\min\left\{
    \frac{8r_1}{23},\,
    \Reg_\infty^\star+\frac{r_1}{2\sqrt m}
\right\}.
\]

As with $r_1/2$ for $\Reg_0^\star$
(Corollary~\ref{cor:G0half}) and $8r_1/23$ for
$\Reg_1^\star$ (Lemma~\ref{lem:one:823}), the limiting
minimax regret has a tight universal bound.
The three-group interpretation suggests dividing the revenue
range $[0,r_1]$ into three intervals of length $r_1/3$.
Accordingly, we experiment with the products whose revenues
are at least $2r_1/3$: define
\begin{equation}\label{eq:rep:two-thirds-index}
k:=\max\left\{i\in[n]:r_i\ge\frac{2r_1}{3}\right\}.
\end{equation}
The next lemma specifies the cutoff for expanding the final
assortment and establishes $\Reg_\infty^\star\le r_1/3$.

\begin{lemma}\label{lem:rep:one-third}
For every $n\ge2$ and all revenues $r_1>\cdots>r_n>0$,
let $k$ be as in~\eqref{eq:rep:two-thirds-index} and let
$\ell:=\min\{i\in[n]:r_{i+1}\le r_1/3\}$, where $r_{n+1}:=0$. Then
\[
\max\{r_1-r_k,\;r_{k+1}-r_\ell,\;r_{\ell+1}\}
\le \frac{r_1}{3},
\]
so $\Reg_\infty^\star\le r_1/3$.
This bound on $\Reg_\infty^\star$ is tight.
\end{lemma}

Tightness is attained already with three products:
for the revenue vector $(3,2,1)$, the limiting minimax
regret equals $1=r_1/3$. Thus, the coefficient $1/3$
cannot be reduced in a bound valid for all revenue vectors.
The same example shows that repeating a single experiment can yield no improvement over no experimentation: Theorem~\ref{thm:G0} gives $\Reg_0^\star=1$, attained by the prefix $[2]$, so $\Reg_m^\star=1$ for every $m$ because $\Reg_\infty^\star\le\Reg_m^\star\le\Reg_0^\star$.
We defer the proof of Lemma~\ref{lem:rep:one-third}
to Appendix~\ref{app:rep:one-third}. These cutoffs give a simple policy without solving the
minimization in~\eqref{eq:Ginf-explicit-intro}.
The decision maker offers $[k]$ to $m$ customers and selects
$[k]$ if the observed purchase fraction satisfies
${\hat z\ge r_\ell/r_k}$, and $[\ell]$ otherwise.
The policy analysis underlying
Theorem~\ref{thm:rep:main}, together with
Lemma~\ref{lem:rep:one-third}, guarantees worst-case regret
at most $r_1/3+r_1/(2\sqrt m)$.

\noindent\textbf{Examples.}
We return to the four revenue profiles introduced in
Section~\ref{sec:G0}. Table~\ref{tab:examples-Ginf}
summarizes the regrets and relative regret reductions
$(\Reg_0^\star-\Reg_j^\star)/\Reg_0^\star$, for
$j=1,m,\infty$, as $n\to\infty$ through even integers.
For $\Reg_\infty^\star$, the limit $m\to\infty$ is taken
before $n\to\infty$. The finite-$m$ column combines the
one-observation guarantee with
Theorem~\ref{thm:rep:main} to give lower bounds on the
limit inferior of the relative regret reduction as
$n\to\infty$, for each fixed $m$.

\begin{table}[t]
\centering
\small
\renewcommand{\arraystretch}{1.8}
\setlength{\tabcolsep}{4pt}
\resizebox{\linewidth}{!}{%
\begin{tabular}{lccc|ccc}
\hline
Example & $\Reg_0^\star$ & $\Reg_1^\star$ & $\Reg_\infty^\star$
& \makecell{Regret reduction\\with 1 observation}
& \makecell{Regret reduction\\with $m$ observations\\(lower bound)}
& \makecell{Regret reduction\\($m\to\infty$)} \\
\hline
\makecell[l]{(a) Uniformly spaced\\\phantom{(a) }revenues} & $\frac{R}{2}$ & $\frac{8R}{23}$ & $\frac{R}{3}$ & $\frac{7}{23}\approx30.4\%$ & $\max\{30.4,33.3-\frac{100}{\sqrt m}\}\%$ & $\frac{1}{3}\approx33.3\%$ \\[2pt]
\makecell[l]{(b) One product at $R$;\\\phantom{(b) }others near $R/2$} & $\frac{R}{2}$ & $\frac{R}{16}$ & $0$ & $\frac{7}{8}=87.5\%$ & $\max\{87.5,100-\frac{100}{\sqrt m}\}\%$ & $100\%$ \\[2pt]
\makecell[l]{(c) Clusters near\\\phantom{(c) }$R$ and $R/3$} & $\frac{R}{3}$ & $\frac{R}{9}$ & $0$ & $\frac{2}{3}\approx{66.6}\%$ & $\max\{{66.6},100-\frac{150}{\sqrt m}\}\%$ & $100\%$ \\[2pt]
\makecell[l]{(d) Clusters near\\\phantom{(d) }$R$ and $R/8$} & $\frac{R}{8}$ & $\frac{R}{8}$ & $0$ & $0\%$ & $\max\{0,100-\frac{400}{\sqrt m}\}\%$ & $100\%$ \\[2pt]
\hline
\end{tabular}}
\caption{Regrets and relative regret reductions for examples (a)--(d) as $n\to\infty$ through even integers, with $R$ fixed. Regret reductions are relative to the no-observation minimax regret $\Reg_0^\star$. Percentage coefficients are rounded down to one decimal place.}
\label{tab:examples-Ginf}
\end{table}

For uniformly spaced revenues, additional observations
provide only a modest improvement beyond one observation:
the relative regret reduction increases from approximately
$30.4\%$ to $33.3\%$. Thus, a single observation already
achieves most of the possible regret reduction.
In examples~(b) and~(c), additional observations have
substantially greater value, with the relative regret
reduction approaching $100\%$ in the stated limits.
Example~(d) gives a particularly striking contrast:
one observation yields exactly zero regret reduction,
whereas multiple observations allow the reduction to
approach $100\%$. For each fixed $n$ in examples~(b)--(d),
however, $\Reg_\infty^\star=\epsilon_n=R/n$; this remaining
regret vanishes only as $n\to\infty$.
Table~\ref{tab:examples-Ginf-finite} in
Appendix~\ref{app:running-examples} gives the finite-$n$
expressions.

Table~\ref{tab:examples-policies} illustrates the structure
of the corresponding experiments and decision rules,
with large-$n$ approximations indicated in the caption.
The exact finite-$n$ policies follow from
Theorems~\ref{thm:G0}--\ref{thm:rep:main} and
Appendix~\ref{app:running-examples}. They are minimax
optimal with zero or one observation and asymptotically
minimax optimal as $m\to\infty$ with multiple observations.

In examples~(b)--(d), the same experimental assortment
supports an optimal policy with one observation and an
asymptotically minimax optimal policy with multiple
observations. The additional observations allow the
decision maker to use the purchase frequency to decide
whether to retain the expensive products alone or expand
to the full assortment. This distinction is especially
important in example~(d), where the one-observation policy
ignores the observation entirely. In example~(a), the
displayed policies instead use different experimental
assortments as well as different decision rules.

\begin{table}[t]
\centering
\normalsize
\renewcommand{\arraystretch}{1.6}
\setlength{\tabcolsep}{4pt}
\begin{tabular}{lc|cl|cl}
\hline
& \multicolumn{1}{c|}{No observation} & \multicolumn{2}{c|}{One observation} & \multicolumn{2}{c}{\makecell{Multiple observations\\($m\to\infty$)}} \\
\cline{2-6}
Example & Decision rule & Experiment & Decision rule & Experiment & Decision rule \\
\hline
(a) & $[n/2]$ & $\approx[7n/23]$ & \Gape[4pt][0pt]{\makecell[l]{Purchase: $\approx[8n/23]$\\No purchase: $\approx[15n/23]$}} & $\approx[n/3]$ & \makecell[l]{${\hat z\ge1/2}$: $\approx[n/3]$\\Otherwise: $\approx[2n/3]$} \\[4pt]
(b) & $[1]$ & $[1]$ & \makecell[l]{Purchase: $[1]$\\No purchase: $[n]$} & $[1]$ & \makecell[l]{$\hat z\ge1/2$: $[1]$\\Otherwise: $[n]$} \\[4pt]
(c) & $[n/2]$ & $[n/2]$ & \makecell[l]{Purchase: $[n/2]$\\No purchase: $[n]$} & $[n/2]$ & \makecell[l]{$\hat z\ge1/3$: $[n/2]$\\Otherwise: $[n]$} \\[4pt]
(d) & $[n/2]$ & $[n/2]$ & \makecell[l]{$[n/2]$ for\\either outcome} & $[n/2]$ & \makecell[l]{$\hat z\ge1/8$: $[n/2]$\\Otherwise: $[n]$} \\[4pt]
\hline
\end{tabular}
\caption{Policies for examples (a)--(d), with even $n\ge20$. Here $\hat z$ is the observed purchase fraction. The $\approx$ entries and thresholds ${1/2}$, $1/3$, and $1/8$ show large-$n$ approximations.}
\label{tab:examples-policies}
\end{table}

Table~\ref{tab:examples-Ginf} shows that, for uniformly
spaced revenues, the limiting minimax regret approaches
the universal bound $r_1/3$ in
Lemma~\ref{lem:rep:one-third} as $n\to\infty$.
More generally, for every fixed revenue vector with at
least three distinct positive revenues,
$\Reg_\infty^\star>0$: repeating a single experimental
assortment cannot eliminate worst-case regret.
Section~\ref{sec:limit} studies whether using several
experimental assortments, chosen in advance, can
eliminate this remaining regret.

\section{Assortment Optimization with Multiple Experiments}\label{sec:limit}

{%
In this section, we generalize the setting of Section~\ref{sec:rep}. We study the value of multiple distinct experimental assortments that the decision maker commits to before observing any customer choices. For a fixed budget on the number of distinct assortments, we characterize the limiting minimax regret and construct a policy whose worst-case regret converges to this limit as the number of customer observations grows. We also characterize which sets of distinct experimental assortments suffice to guarantee expected revenue arbitrarily close to optimal.

\noindent\textbf{Policies with multiple distinct experiments.} We start by extending the setting of Section~\ref{sec:model} to allow multiple distinct experimental assortments. We fix in advance a sequence of possibly different experimental assortments $S_1,\dots,S_m\subseteq[n]$. In experimental round $t$, we offer $S_t$ and observe a choice $O_t\sim\phi_v(\cdot,S_t)$, independently across rounds. After observing all $m$ choices, a decision rule $\mathcal A$ selects the final assortment $\mathcal A(O_1,\dots,O_m)$. Together, the experimental sequence and the decision rule form a policy. The minimum worst-case regret for this fixed sequence is
\[
\Reg_m^{\mathrm{seq}}(S_1,\dots,S_m)
:=\min_{\mathcal A}\sup_{v\in\mathbb R_+^n}
\mathbb E_v\!\left[R_v^\star-R_v\!\left(\mathcal A(O_1,\dots,O_m)\right)\right].
\]
For integers $m\geq d\geq1$, we define the minimax regret with $m$ observations and at most $d$ distinct experimental assortments by minimizing over the experimental sequence:
\[
\Reg_{m,d}^\star
:=\min_{\substack{S_1,\dots,S_m\subseteq[n]\\
|\{S_1,\dots,S_m\}|\le d}}
\Reg_m^{\mathrm{seq}}(S_1,\dots,S_m).
\]
When $d=1$, we offer the same experimental assortment in all $m$ rounds, so we recover the setting of Section~\ref{sec:rep} and have $\Reg_{m,1}^\star=\Reg_m^\star$. We extend the notation to $m=d=0$ by setting $\Reg_{0,0}^\star:=\Reg_0^\star$, the no-observation minimax regret of Section~\ref{sec:G0}.

\subsection{Limiting Minimax Regret with Multiple Experimental Assortments}\label{sec:limit:minimax}

In this section, we characterize the limit of $\Reg_{m,d}^{\star}$ as the observation budget $m$ grows, for any fixed $d\ge1$. We then construct a policy that uses at most $d$ distinct nested revenue prefixes as experimental assortments and whose worst-case regret converges to this limit.

We prove that this limit is characterized by a geometric construction naturally extending that of Sections~\ref{sec:G0} and~\ref{sec:rep}. Place the revenues, including zero, on the real line and consider partitions $\mathcal P$ of $\{r_1,\ldots,r_n,0\}$ into at most $d+2$ nonempty consecutive groups. The diameter of a group $B\in\mathcal P$ is $\max B-\min B$, the difference between its largest and smallest revenues. This diameter is zero for a singleton. We define the minimum possible largest group diameter by
\begin{equation}\label{eq:distinct:partition}
\Reg_{\infty,d}^\star
:=\min_{\mathcal P}\max_{B\in\mathcal P}(\max B-\min B).
\end{equation}
Splitting a group cannot increase the largest group diameter, so we may choose a minimizing partition with exactly $\min\{d+2,n+1\}$ groups. For $d=0,1$, the partition formula recovers the regret values in Theorems~\ref{thm:G0} and~\ref{thm:rep:main}, respectively (Figures~\ref{fig:G0-partition} and~\ref{fig:Ginf-partition}).

To illustrate the formula, consider the case $d=2$ and $n>3$. Choose four consecutive groups, ordered by decreasing revenue, and let $a<b<c$ denote the last product indices in the first three groups. The fourth contains the remaining revenues and zero. Minimizing the largest group diameter gives
\[
\Reg_{\infty,2}^\star
=\min_{1\le a<b<c\le n}
\max\big\{r_1-r_a,\ r_{a+1}-r_b,\ r_{b+1}-r_c,\ r_{c+1}\big\}{,}
\]
where $r_{n+1}:=0$. Figure~\ref{fig:distinct:partition} illustrates this partition.

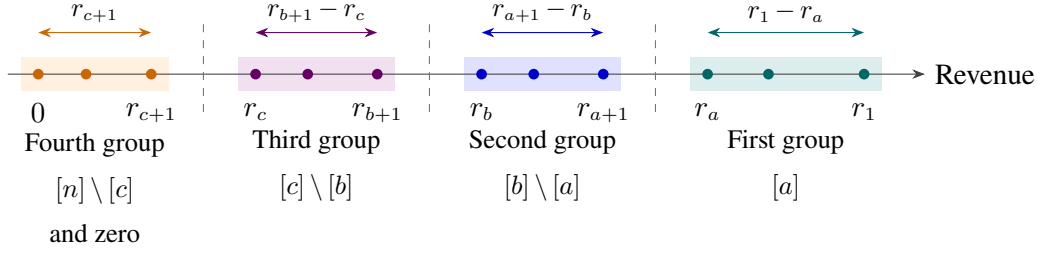
\begin{figure}[ht]
\centering
\begin{tikzpicture}[x=1.15cm,y=1cm,draw=black!75,text=black,>=Stealth]
  \fill[orange!12] (-0.2,-0.22) rectangle (1.5,0.22);
  \fill[violet!12] (2.3,-0.22) rectangle (4.1,0.22);
  \fill[blue!12] (4.9,-0.22) rectangle (6.7,0.22);
  \fill[teal!14] (7.5,-0.22) rectangle (9.7,0.22);
  \draw[->] (-0.35,0)--(10.2,0) node[right] {Revenue};
  \foreach \x in {0,0.55,1.3}
    \fill[orange!80!black] (\x,0) circle (2pt);
  \foreach \x in {2.5,3.1,3.9}
    \fill[violet!80!black] (\x,0) circle (2pt);
  \foreach \x in {5.1,5.7,6.5}
    \fill[blue!80!black] (\x,0) circle (2pt);
  \foreach \x in {7.7,8.4,9.5}
    \fill[teal!80!black] (\x,0) circle (2pt);
  \node[below=7pt] at (0,0) {$0$};
  \node[below=7pt] at (1.3,0) {$r_{c+1}$};
  \node[below=7pt] at (2.5,0) {$r_c$};
  \node[below=7pt] at (3.9,0) {$r_{b+1}$};
  \node[below=7pt] at (5.1,0) {$r_b$};
  \node[below=7pt] at (6.5,0) {$r_{a+1}$};
  \node[below=7pt] at (7.7,0) {$r_a$};
  \node[below=7pt] at (9.5,0) {$r_1$};
  \foreach \x in {1.9,4.5,7.1}
    \draw[dashed,black!55] (\x,-0.5)--(\x,0.7);
  \draw[<->,orange!80!black] (0,0.55)--(1.3,0.55)
    node[midway,above,font=\small,text=black] {$r_{c+1}$};
  \draw[<->,violet!80!black] (2.5,0.55)--(3.9,0.55)
    node[midway,above,font=\small,text=black] {$r_{b+1}-r_c$};
  \draw[<->,blue!80!black] (5.1,0.55)--(6.5,0.55)
    node[midway,above,font=\small,text=black] {$r_{a+1}-r_b$};
  \draw[<->,teal!80!black] (7.7,0.55)--(9.5,0.55)
    node[midway,above,font=\small,text=black] {$r_1-r_a$};
  \node[font=\small,align=center,text=black] at (0.65,-1.5) {Fourth group\\$[n]\setminus[c]$\\and zero};
  \node[font=\small,align=center] at (3.2,-1.2) {Third group\\$[c]\setminus[b]$};
  \node[font=\small,align=center] at (5.8,-1.2) {Second group\\$[b]\setminus[a]$};
  \node[font=\small,align=center] at (8.6,-1.2) {First group\\$[a]$};
\end{tikzpicture}
\caption{Four nonempty consecutive revenue groups for two distinct experimental assortments, with $a<b<c$. The value $\Reg_{\infty,2}^{\star}$ is the smallest possible largest group diameter.}
\label{fig:distinct:partition}
\end{figure}

The next theorem shows that, for fixed $d\ge1$, $\Reg_{\infty,d}^{\star}$ gives the limiting minimax regret: $\Reg_{m,d}^{\star}$ approaches this value as the observation budget $m$ grows. Thus, a policy using at most $d$ distinct experimental assortments can achieve worst-case regret as close as desired to $\Reg_{\infty,d}^{\star}$ with a sufficiently large but finite number of repetitions.

\begin{theorem}[Limiting Minimax Regret with Multiple Experiments]\label{thm:vanish}\label{thm:distinct}
{For all integers $m\ge d\ge1$,}
\begin{equation}\label{eq:distinct:sandwich}
\Reg_{\infty,d}^\star
\le\Reg_{m,d}^\star
\le\Reg_{\infty,d}^\star+
{\frac{r_1}{2\sqrt{\lfloor m/d\rfloor}}
\min\!\left\{\sqrt d,\sqrt{2\ln(d+1)}\right\}}.
\end{equation}

Consequently, for every fixed $d\ge1$,
\[
\lim_{m\to\infty}\Reg_{m,d}^\star=\Reg_{\infty,d}^\star.
\]

Moreover, there exists a policy that chooses its entire experimental sequence in advance, offers at most $d$ distinct nested revenue-ordered assortments over $m$ rounds, and has worst-case regret at most the right-hand side of~\eqref{eq:distinct:sandwich}.
\end{theorem}

{We defer the proof of Theorem~\ref{thm:vanish} to Appendix~\ref{app:limit:minimax}. We first present a policy that achieves the upper bound and then explain its experimental design and decision rule.}

{%
\noindent\textbf{Revenue prefix policy with rounding.} Fix a partition whose largest group diameter equals $\Reg_{\infty,d}^{\star}$. For each product $i$, let $\widetilde r_i$ be the minimum revenue in its group, and set $\widetilde r_0=0$. Order the $h\le d+1$ groups with positive minima by decreasing revenue. For $j=1,\ldots,h$, let $k_j$ be the largest product index in group $j$, corresponding to the product with the lowest original revenue.

If $h\ge2$, offer each of $[k_1],\ldots,[k_{h-1}]$ in $\tau:=\lfloor m/d\rfloor$ rounds. Offer $[k_1]$ in the remaining rounds and ignore those additional observations. Let $O_{j,t}$ be the $t$th retained choice from $[k_j]$, and estimate its expected rounded revenue by
\[
\widehat R_j
:=\frac{1}{\tau}\sum_{t=1}^{\tau}\widetilde r_{O_{j,t}},
\qquad
\forall\,j\in\{1,\ldots,h-1\}.
\]
Choose the smallest index $\widehat j\in\{1,\ldots,h-1\}$ maximizing $\widehat R_j$, and commit to
\begin{equation}\label{eq:distinct:policy}
\mathcal A(O):=
\begin{cases}
[k_{\widehat j}],
&\widehat R_{\widehat j}\ge\widetilde r_{k_h},\\[2pt]
[k_h],
&\widehat R_{\widehat j}<\widetilde r_{k_h}.
\end{cases}
\end{equation}
If $h\le1$, offer $[1]$ throughout, ignore the observations, and commit to $\{i\in[n]:\widetilde r_i>0\}$, which is empty when $h=0$.

\FloatBarrier

}

Figure~\ref{fig:distinct:policy} illustrates the policy for the $d=2$ partition in Figure~\ref{fig:distinct:partition} and even $m\ge2$. We offer each of $[a]$ and $[b]$ in $m/2$ rounds; their sample average rounded revenues, labeled ${\widehat R_1}$ and ${\widehat R_2}$, determine the final assortment through comparison with $r_c$.

\begin{figure}[H]
\centering
\begin{tikzpicture}[x=1.15cm,y=1cm,draw=black!75,text=black,>=Stealth]
  \fill[orange!12] (-0.2,-0.22) rectangle (1.5,0.22);
  \fill[violet!12] (2.3,-0.22) rectangle (4.1,0.22);
  \fill[blue!12] (4.9,-0.22) rectangle (6.7,0.22);
  \fill[teal!14] (7.5,-0.22) rectangle (9.7,0.22);
  \draw[->] (-0.35,0)--(10.2,0) node[right] {Revenue};
  \foreach \x in {0,0.55,1.3}
    \fill[orange!80!black] (\x,0) circle (2pt);
  \foreach \x in {2.5,3.1,3.9}
    \fill[violet!80!black] (\x,0) circle (2pt);
  \foreach \x in {5.1,5.7,6.5}
    \fill[blue!80!black] (\x,0) circle (2pt);
  \foreach \x in {7.7,8.4,9.5}
    \fill[teal!80!black] (\x,0) circle (2pt);
  \node[below=7pt] at (0,0) {$0$};
  \node[below=7pt] at (1.3,0) {$r_{c+1}$};
  \node[below=7pt] at (2.5,0) {$r_c$};
  \node[below=7pt] at (3.9,0) {$r_{b+1}$};
  \node[below=7pt] at (5.1,0) {$r_b$};
  \node[below=7pt] at (6.5,0) {$r_{a+1}$};
  \node[below=7pt] at (7.7,0) {$r_a$};
  \node[below=7pt] at (9.5,0) {$r_1$};
  \foreach \x in {1.9,4.5,7.1}
    \draw[dashed,black!55] (\x,-0.5)--(\x,0.7);

  \draw[decorate,decoration={brace,amplitude=4pt},thick,draw=teal!80!black]
    (7.7,0.8)--(9.5,0.8)
    node[midway,above=7pt,font=\small,align=center,text=black]
    {\textbf{Experiment $[a]$:} with $m/2$ customers to estimate ${\widehat R_1}$.};
  \draw[decorate,decoration={brace,amplitude=4pt},thick,draw=blue!80!black]
    (5.1,1.9)--(9.5,1.9)
    node[midway,above=7pt,font=\small,align=center,text=black]
    {\textbf{Experiment $[b]$:} with $m/2$ customers to estimate ${\widehat R_2}$.};

  \foreach \x in {2.5,5.1,7.7}
    \draw[thick] (\x,0) circle (4pt);
  \draw[->,thick] (2.5,-0.6)--(2.5,-0.9)--(1.8,-1.3);
  \draw[->,thick] (5.1,-0.6)--(5.1,-0.9)--(5.4,-1.3);
  \draw[->,thick] (7.7,-0.6)--(7.7,-0.9)--(9,-1.3);
  \node[anchor=north,align=center,text width=4cm,font=\small] at (1.8,-1.5)
    {\textbf{Offer $[c]$}\\[5pt]
     if $\max\{{\widehat R_1},{\widehat R_2}\}<r_c$};
  \node[anchor=north,align=center,text width=4cm,font=\small] at (5.4,-1.5)
    {\textbf{Offer $[b]$}\\[5pt]
     if ${\widehat R_2}>{\widehat R_1}$ and ${\widehat R_2}\ge r_c$};
  \node[anchor=north,align=center,text width=4cm,font=\small] at (9,-1.5)
    {\textbf{Offer $[a]$}\\[5pt]
     if ${\widehat R_1}\ge\max\{{\widehat R_2},r_c\}$};
\end{tikzpicture}
\caption{The policy for $d=2$ and even $m$, using a minimizing revenue partition of the form shown in Figure~\ref{fig:distinct:partition}. The experiments are $[a]$ and $[b]$; the final assortment is $[a]$, $[b]$, or $[c]$.}
\label{fig:distinct:policy}
\end{figure}
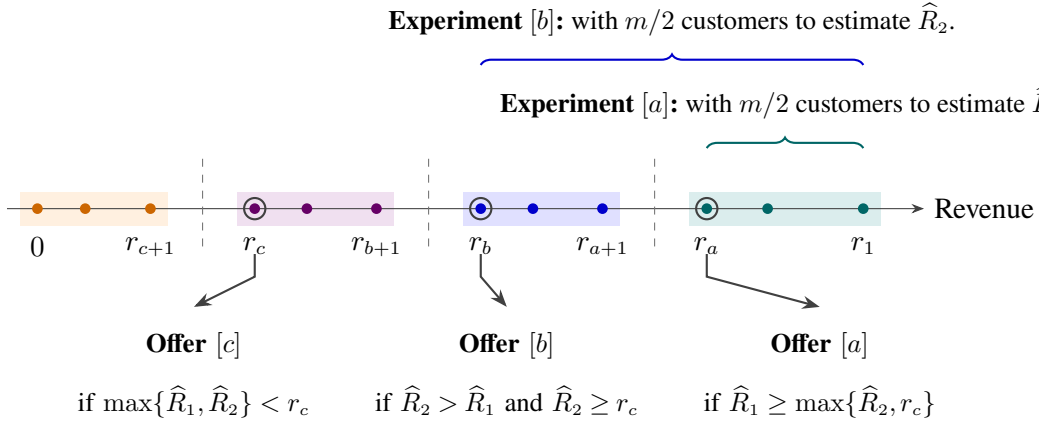
\FloatBarrier

{We now interpret the policy and explain why it achieves the upper bound in Theorem~\ref{thm:vanish}. The policy chooses its final assortment from $[k_1],\ldots,[k_h]$. These candidates suffice under the rounded revenues: the revenue-ordering property of the MNL model allows an optimal assortment to include or exclude products with equal revenue together. Since products within each group share the same rounded revenue, we can choose an optimal assortment consisting of whole groups in decreasing revenue order and excluding the group rounded to zero. Such an assortment is one of the policy's $h\le d+1$ candidate prefixes.}

{The experiments reveal whether adding lower-revenue groups increases expected rounded revenue. Adding such a group can turn no purchase into a purchase, but it can also divert purchases from higher-revenue products. For each offered prefix, averaging the chosen products' rounded revenues, with zero for no purchase, estimates its expected rounded revenue. Comparing these estimates guides the choice of how many groups to include in the final assortment. Under MNL, the policy's threshold comparison lets us leave the longest prefix untested, so $h-1\le d$ distinct experimental assortments suffice.}

{The policy's regret under the original revenues has two sources: rounding and estimation. Rounding leaves choice probabilities unchanged and lowers each product's revenue by at most the largest group diameter, $\Reg_{\infty,d}^{\star}$, so any assortment's expected revenue falls by at most this amount. In particular, the original optimal assortment loses at most this amount. An optimal assortment under the rounded revenues does at least as well as the original optimal assortment in the rounded problem, and its revenue can only increase under the original revenues. Under the original revenues, this assortment's expected revenue is therefore within $\Reg_{\infty,d}^{\star}$ of optimal. Since observed rounded revenues lie in $[0,r_1]$, the additional expected loss from using estimates vanishes uniformly over preference weights as repetitions increase, for fixed $d$. Together, rounding and estimation give the two terms in Theorem~\ref{thm:vanish}'s upper bound.}

{%
The partition formula also gives an upper bound that depends only on $r_1$, $m$, and $d$. Partitioning the revenue points according to $d+2$ equal intervals of $[0,r_1]$, and omitting empty groups, gives $\Reg_{\infty,d}^{\star}\le r_1/(d+2)$. Applying Theorem~\ref{thm:vanish} then yields the following corollary.

\begin{corollary}\label{cor:distinct:universal}
For all integers $m\ge d\ge1$,
\[
\Reg_{m,d}^{\star}
\le r_1\left[
\frac{1}{d+2}
+\frac{\min\!\left\{\sqrt d,\sqrt{2\ln(d+1)}\right\}}
{2\sqrt{\lfloor m/d\rfloor}}
\right].
\]
\end{corollary}

We defer the proof of Corollary~\ref{cor:distinct:universal} to Appendix~\ref{app:limit:corollary}. For a prescribed $d$, this corollary shows that an observation budget of order $d^3\ln(d+1)$ suffices to guarantee worst-case regret $O(r_1/(d+2))$.
}

{%
\subsection{Necessary and Sufficient Conditions for Vanishing Regret}\label{sec:diff}

In this section, we characterize when a given set of distinct experimental assortments permits vanishing worst-case regret through repetition. Theorem~\ref{thm:vanish} implies that vanishing regret requires $\Reg_{\infty,d}^{\star}=0$. By~\eqref{eq:distinct:partition}, each of the $n+1$ distinct revenue points, including zero, must therefore lie in its own group. Since at most $d+2$ groups are allowed, at least $n-1$ distinct experimental assortments are necessary. Thus, guaranteeing revenue arbitrarily close to optimal requires nearly as many distinct experimental assortments as there are products, even when customer observations are abundant.

Having sufficiently many distinct experimental assortments is necessary but not sufficient for vanishing regret. The following theorem characterizes exactly which sets of experimental assortments permit vanishing worst-case regret.

\begin{theorem}[Necessary and Sufficient Conditions for Vanishing Regret]\label{thm:distinct:vanishing}
Assume $n\ge2$, and fix a set $\mathcal F$ of $d\ge1$ distinct experimental assortments. For each integer $N\ge1$, consider policies using only assortments in $\mathcal F$, with each assortment offered at most $N$ times. The minimum worst-case regret over these policies converges to zero as $N\to\infty$ if and only if
\[
\forall\,i\in\{1,\ldots,n-1\},\qquad
\exists\,S\in\mathcal F\ \text{such that}\ i\in S\subseteq[i].
\]
\end{theorem}

We defer the proof of Theorem~\ref{thm:distinct:vanishing} to Appendix~\ref{app:limit:vanishing}. We now give intuition for why the condition is necessary and sufficient.

For necessity, suppose the condition fails for some product $j<n$. Then product $j$ is either never offered or always offered with at least one lower-revenue product. In the first case, its preference weight cannot be learned from the observations, while that weight can affect the optimal final assortment. In the second case, the lower-revenue products can attract almost every purchase, masking the preference of product $j$. For any finite observation budget, their weights can be large enough to prevent us from reliably deciding whether to offer $j$ alone or include a lower-revenue product in the final assortment. The resulting worst-case regret remains bounded away from zero as customer observations increase.

For sufficiency, suppose the condition holds. Without loss of generality, assume that for every product $j<n$, its corresponding experiment $S_j\in\mathcal F$ satisfies $j\in S_j\subseteq[j]$. After $N$ repetitions of $S_j$, we record how often customers choose product $j$ or make no purchase. When the observed no-purchase frequency is positive, we estimate $v_j$ by dividing the observed frequency of purchases of product $j$ by the observed no-purchase frequency. If the no-purchase probability in $S_j$ is not too small, a sufficiently large finite $N$ gives a small estimation error with high probability. If the no-purchase probability in $S_j$ is too small, we may be unable to estimate $v_j$ with small error. However, interestingly, we can then restrict attention to assortments $C\subseteq[j]$. To see why, every purchase from $S_j$ earns at least $r_j$, so its expected revenue exceeds $r_{j+1}$ once no-purchases are sufficiently rare. For any $C\subseteq[j]$, the expected revenue after adding products outside $[j]$ is a convex combination of $R_v(C)$ and the added products' revenues. Those revenues are at most $r_{j+1}$, so these additions cannot improve the best expected revenue over $C\subseteq[j]$. Thus, if some product $j<n$ has a sufficiently small no-purchase probability in its corresponding experiment $S_j$, an optimal assortment can be chosen as a subset of $[j]$.

Now suppose $i<n$ is the first product in decreasing revenue order whose experiment $S_i$ has a small observed no-purchase frequency. For sufficiently many repetitions, with high probability, we can restrict attention to $C\subseteq[i]$ and accurately estimate all weights $v_j$ for $j<i$. Under MNL, an optimal assortment within $[i]$ can be chosen independently of the lowest-revenue product's preference weight $v_i$. For this restricted problem, we use these sufficiently accurate MNL parameter estimates to find a revenue prefix $C\subseteq[i]$ that is optimal under the estimated parameters. The selected assortment is then nearly optimal under the true parameters. If no such $i<n$ exists, the same approach applies to the full problem using estimates of $v_1,\ldots,v_{n-1}$.

For each product $i<n$, the experiment required by Theorem~\ref{thm:distinct:vanishing} may contain any subset of the higher-revenue products. Including all of them gives a revenue prefix, while including none gives a singleton. Including only some is also allowed, so vanishing regret does not require every experimental assortment to have either form.

\noindent\textbf{Examples.} For $n=5$, consider the following three sets of experimental assortments:
\[
\begin{aligned}
{\mathcal F_1:=}\;&
\bigl\{\{1\},\{1,2\},\{1,2,3\},\{1,2,3,4\}\bigr\},\\
{\mathcal F_2:=}\;&
\bigl\{\{1\},\{2\},\{3\},\{4\}\bigr\},\\
{\mathcal F_3:=}\;&
\bigl\{\{1\},\{1,2\},\{2,3\},\{4\}\bigr\}.
\end{aligned}
\]
In ${\mathcal F_3}$, product $3$ is offered with product $2$ but without product $1$. Thus, $\{2,3\}$ is neither a revenue prefix nor a singleton, yet it satisfies the requirement for product $3$. In all three sets, the four displayed assortments satisfy the requirement for products $1,2,3,4$, respectively, since each contains its corresponding product and no lower-revenue product. Each set uses $d=4$ distinct experimental assortments, the minimum number required. Offering each assortment $N$ times gives $m=4N$ observations. For each set, Theorem~\ref{thm:distinct:vanishing} guarantees suitable decision rules whose worst-case regret tends to zero as $N\to\infty$.

The characterization also holds when the assortments are offered different numbers of times, provided the number of observations from each assortment tends to infinity. For sufficiency, let $N$ be the smallest repetition count, use $N$ observations from each assortment, and ignore the extras. The necessity argument also applies to unequal repetition counts.

}

}

{%
\section{Conclusion}\label{sec:conclusion}
We take a decision-focused approach to limited, non-adaptive experimentation, designing experiments to inform which assortment to offer rather than to estimate customer preferences. We show that, surprisingly, even a single experiment can substantially reduce the retailer's regret, even with very limited customer observations. We characterize a minimax-optimal experimentation policy with a single customer observation and construct near-optimal policies with multiple observations, providing regret guarantees in each setting. We also extend our analysis to multiple experiments run in parallel. For any fixed budget of distinct experimental assortments, our geometric characterization connects limiting minimax regret to partitions of product revenues, yielding asymptotically optimal policies and precise conditions for vanishing regret. While our analysis focuses on MNL assortment optimization, we hope this decision-focused experimentation framework motivates research on experimentation under other choice models and in related assortment settings.

{%
\addcontentsline{toc}{section}{Bibliography}
\bibliographystyle{abbrvnat}

\bibliography{ref}

\begin{thebibliography}{38}
\providecommand{\natexlab}[1]{#1}
\providecommand{\url}[1]{\texttt{#1}}
\expandafter\ifx\csname urlstyle\endcsname\relax
  \providecommand{\doi}[1]{doi: #1}\else
  \providecommand{\doi}{doi: \begingroup \urlstyle{rm}\Url}\fi

\bibitem[Agrawal et~al.(2017)Agrawal, Avadhanula, Goyal, and Zeevi]{agrawal2017thompson}
S.~Agrawal, V.~Avadhanula, V.~Goyal, and A.~Zeevi.
\newblock Thompson sampling for the {MNL}-bandit.
\newblock In \emph{Proceedings of the 2017 Conference on Learning Theory}, volume~65 of \emph{Proceedings of Machine Learning Research}, pages 76--78. PMLR, 2017.

\bibitem[Agrawal et~al.(2019)Agrawal, Avadhanula, Goyal, and Zeevi]{agrawal2019}
S.~Agrawal, V.~Avadhanula, V.~Goyal, and A.~Zeevi.
\newblock {MNL}-bandit: A dynamic learning approach to assortment selection.
\newblock \emph{Operations Research}, 67\penalty0 (5):\penalty0 1453--1485, 2019.
\newblock \doi{10.1287/opre.2018.1832}.

\bibitem[Ahmadnejadsaein and El~Housni(2025)]{ahmadnejadsaein2025adaptive}
M.~Ahmadnejadsaein and O.~El~Housni.
\newblock Adaptive two-sided assortment optimization: Revenue maximization.
\newblock arXiv preprint arXiv:2507.04156, 2025.
\newblock Revised September 2026.

\bibitem[Allouah et~al.(2022)Allouah, Bahamou, and Besbes]{allouah2022samples}
A.~Allouah, A.~Bahamou, and O.~Besbes.
\newblock Pricing with samples.
\newblock \emph{Operations Research}, 70\penalty0 (2):\penalty0 1088--1104, 2022.
\newblock \doi{10.1287/opre.2021.2200}.

\bibitem[Allouah et~al.(2023)Allouah, Bahamou, and Besbes]{allouah2023single}
A.~Allouah, A.~Bahamou, and O.~Besbes.
\newblock Optimal pricing with a single point.
\newblock \emph{Management Science}, 69\penalty0 (10):\penalty0 5866--5882, 2023.
\newblock \doi{10.1287/mnsc.2023.4683}.

\bibitem[Aouad and Saban(2023)]{aouad2023twosided}
A.~Aouad and D.~Saban.
\newblock Online assortment optimization for two-sided matching platforms.
\newblock \emph{Management Science}, 69\penalty0 (4):\penalty0 2069--2087, 2023.
\newblock \doi{10.1287/mnsc.2022.4464}.

\bibitem[Aouad et~al.(2018)Aouad, Farias, Levi, and Segev]{aouad2018ranking}
A.~Aouad, V.~Farias, R.~Levi, and D.~Segev.
\newblock The approximability of assortment optimization under ranking preferences.
\newblock \emph{Operations Research}, 66\penalty0 (6):\penalty0 1661--1669, 2018.
\newblock \doi{10.1287/opre.2018.1754}.

\bibitem[Araman and Caldentey(2022)]{araman2022diffusion}
V.~F. Araman and R.~A. Caldentey.
\newblock Diffusion approximations for a class of sequential experimentation problems.
\newblock \emph{Management Science}, 68\penalty0 (8):\penalty0 5958--5979, 2022.
\newblock \doi{10.1287/mnsc.2021.4195}.

\bibitem[Ashlagi et~al.(2022)Ashlagi, Krishnaswamy, Makhijani, Saban, and Shiragur]{ashlagi2022twosided}
I.~Ashlagi, A.~K. Krishnaswamy, R.~Makhijani, D.~Saban, and K.~Shiragur.
\newblock Technical note---assortment planning for two-sided sequential matching markets.
\newblock \emph{Operations Research}, 70\penalty0 (5):\penalty0 2784--2803, 2022.
\newblock \doi{10.1287/opre.2022.2327}.

\bibitem[Barr{\'e} et~al.(2026)Barr{\'e}, El~Housni, Ibn~Brahim, Lodi, and Segev]{barre2026visibility}
T.~Barr{\'e}, O.~El~Housni, M.~Ibn~Brahim, A.~Lodi, and D.~Segev.
\newblock Assortment optimization with visibility constraints.
\newblock \emph{Mathematical Programming}, 216\penalty0 (1--2):\penalty0 177--220, 2026.
\newblock \doi{10.1007/s10107-025-02238-9}.

\bibitem[Bennouna et~al.(2025)Bennouna, Bennouna, Amin, and Ozdaglar]{bennouna2025data}
O.~Bennouna, A.~Bennouna, S.~Amin, and A.~Ozdaglar.
\newblock What data enables optimal decisions? {An} exact characterization for linear optimization.
\newblock In \emph{Advances in Neural Information Processing Systems}, volume~38, 2025.
\newblock \doi{10.52202/085713-5534}.

\bibitem[Bennouna et~al.(2026)Bennouna, Bennouna, Amin, and Ozdaglar]{bennounaDataInformativeness2026}
O.~Bennouna, A.~Bennouna, S.~Amin, and A.~Ozdaglar.
\newblock Data informativeness in linear optimization under uncertainty.
\newblock arXiv preprint arXiv:2602.15365, 2026.

\bibitem[Blanchet et~al.(2016)Blanchet, Gallego, and Goyal]{blanchet2016markov}
J.~Blanchet, G.~Gallego, and V.~Goyal.
\newblock A {Markov} chain approximation to choice modeling.
\newblock \emph{Operations Research}, 64\penalty0 (4):\penalty0 886--905, 2016.
\newblock \doi{10.1287/opre.2016.1505}.

\bibitem[Caro and Gallien(2007)]{caro2007}
F.~Caro and J.~Gallien.
\newblock Dynamic assortment with demand learning for seasonal consumer goods.
\newblock \emph{Management Science}, 53\penalty0 (2):\penalty0 276--292, 2007.
\newblock \doi{10.1287/mnsc.1060.0613}.

\bibitem[Chen et~al.(2025)Chen, Cire, Gao, and Wang]{chen2025prediction}
N.~Chen, A.~Cire, P.~Gao, and S.~Wang.
\newblock Assortment optimization without prediction: An end-to-end framework with transaction data.
\newblock Working paper, SSRN 5280529, 2025.

\bibitem[Chen et~al.(2021)Chen, Wang, and Zhou]{chen2021optimal}
X.~Chen, Y.~Wang, and Y.~Zhou.
\newblock Optimal policy for dynamic assortment planning under multinomial logit models.
\newblock \emph{Mathematics of Operations Research}, 46\penalty0 (4):\penalty0 1639--1657, 2021.
\newblock \doi{10.1287/moor.2021.1133}.

\bibitem[Davis et~al.(2014)Davis, Gallego, and Topaloglu]{davis2014nested}
J.~M. Davis, G.~Gallego, and H.~Topaloglu.
\newblock Assortment optimization under variants of the nested logit model.
\newblock \emph{Operations Research}, 62\penalty0 (2):\penalty0 250--273, 2014.
\newblock \doi{10.1287/opre.2014.1256}.

\bibitem[D{\'e}sir et~al.(2020)D{\'e}sir, Goyal, Segev, and Ye]{desir2020constrained}
A.~D{\'e}sir, V.~Goyal, D.~Segev, and C.~Ye.
\newblock Constrained assortment optimization under the {Markov} chain--based choice model.
\newblock \emph{Management Science}, 66\penalty0 (2):\penalty0 698--721, 2020.
\newblock \doi{10.1287/mnsc.2018.3230}.

\bibitem[Dong et~al.(2025)Dong, Mo, Qi, Shi, Fang, and Tarokh]{dong2025pasta}
J.~Dong, W.~Mo, Z.~Qi, C.~Shi, E.~X. Fang, and V.~Tarokh.
\newblock {PASTA}: A unified framework for offline assortment learning.
\newblock arXiv preprint arXiv:2510.01693, 2025.

\bibitem[Dong et~al.(2020)Dong, Li, Zhang, and Zhou]{dong2020switching}
K.~Dong, Y.~Li, Q.~Zhang, and Y.~Zhou.
\newblock Multinomial logit bandit with low switching cost.
\newblock In \emph{Proceedings of the 37th International Conference on Machine Learning}, volume 119 of \emph{Proceedings of Machine Learning Research}, pages 2607--2615. PMLR, 2020.

\bibitem[El~Housni and Topaloglu(2023)]{elhousni2023customization}
O.~El~Housni and H.~Topaloglu.
\newblock Joint assortment optimization and customization under a mixture of multinomial logit models: On the value of personalized assortments.
\newblock \emph{Operations Research}, 71\penalty0 (4):\penalty0 1197--1215, 2023.
\newblock \doi{10.1287/opre.2022.2384}.

\bibitem[El~Housni et~al.(2024)El~Housni, Hennebelle, and Torrico]{elhousni2024twosided}
O.~El~Housni, U.~Hennebelle, and A.~Torrico.
\newblock Two-sided assortment optimization: Adaptivity gaps and approximation algorithms.
\newblock arXiv preprint arXiv:2403.08929, 2024.
\newblock Revised May 2026.

\bibitem[Feldman and Topaloglu(2017)]{feldman2017markov}
J.~B. Feldman and H.~Topaloglu.
\newblock Revenue management under the {Markov} chain choice model.
\newblock \emph{Operations Research}, 65\penalty0 (5):\penalty0 1322--1342, 2017.
\newblock \doi{10.1287/opre.2017.1628}.

\bibitem[Golrezaei et~al.(2014)Golrezaei, Nazerzadeh, and Rusmevichientong]{golrezaei2014real}
N.~Golrezaei, H.~Nazerzadeh, and P.~Rusmevichientong.
\newblock Real-time optimization of personalized assortments.
\newblock \emph{Management Science}, 60\penalty0 (6):\penalty0 1532--1551, 2014.
\newblock \doi{10.1287/mnsc.2014.1939}.

\bibitem[Han et~al.(2025)Han, Zhong, Lu, Blanchet, and Zhou]{han2025observational}
Y.~Han, H.~Zhong, M.~Lu, J.~Blanchet, and Z.~Zhou.
\newblock Learning an optimal assortment policy under observational data.
\newblock arXiv preprint arXiv:2502.06777, 2025.

\bibitem[K{\"o}k et~al.(2008)K{\"o}k, Fisher, and Vaidyanathan]{kok2008assortment}
A.~G. K{\"o}k, M.~L. Fisher, and R.~Vaidyanathan.
\newblock Assortment planning: Review of literature and industry practice.
\newblock In N.~Agrawal and S.~A. Smith, editors, \emph{Retail Supply Chain Management: Quantitative Models and Empirical Studies}, pages 99--153. Springer, 2008.
\newblock \doi{10.1007/978-0-387-78902-6_6}.

\bibitem[Lee and Oh(2026)]{lee2026design}
J.~Lee and M.-H. Oh.
\newblock Optimal design for multinomial logit model with applications to best assortment identification.
\newblock In \emph{Proceedings of the 43rd International Conference on Machine Learning}, volume 306 of \emph{Proceedings of Machine Learning Research}, pages 64899--64932. PMLR, 2026.

\bibitem[Liu et~al.(2023)Liu, Cao, and Shen]{liu2023marginal}
M.~Liu, J.~Cao, and Z.-J.~M. Shen.
\newblock Marginal value of one data point in assortment personalization.
\newblock Working paper, SSRN 4487888, 2023.
\newblock Revised August 2025.

\bibitem[Lu et~al.(2026)Lu, Han, Zhong, Zhou, and Blanchet]{lu2026robust}
M.~Lu, Y.~Han, H.~Zhong, Z.~Zhou, and J.~Blanchet.
\newblock Robust assortment optimization from observational data.
\newblock arXiv preprint arXiv:2602.10696, 2026.

\bibitem[Rusmevichientong and Topaloglu(2012)]{rusmevichientong2012}
P.~Rusmevichientong and H.~Topaloglu.
\newblock Robust assortment optimization in revenue management under the multinomial logit choice model.
\newblock \emph{Operations Research}, 60\penalty0 (4):\penalty0 865--882, 2012.
\newblock \doi{10.1287/opre.1120.1063}.

\bibitem[Rusmevichientong et~al.(2010)Rusmevichientong, Shen, and Shmoys]{rusmevichientong2010}
P.~Rusmevichientong, Z.-J.~M. Shen, and D.~B. Shmoys.
\newblock Dynamic assortment optimization with a multinomial logit choice model and capacity constraint.
\newblock \emph{Operations Research}, 58\penalty0 (6):\penalty0 1666--1680, 2010.
\newblock \doi{10.1287/opre.1100.0866}.

\bibitem[Saha and Gaillard(2025)]{saha2025rankbreaking}
A.~Saha and P.~Gaillard.
\newblock Finally rank-breaking conquers {MNL} bandits: Optimal and efficient algorithms for {MNL} assortment.
\newblock In \emph{Proceedings of the Thirteenth International Conference on Learning Representations}, 2025.

\bibitem[Saur{\'e} and Zeevi(2013)]{saure2013}
D.~Saur{\'e} and A.~Zeevi.
\newblock Optimal dynamic assortment planning with demand learning.
\newblock \emph{Manufacturing \& Service Operations Management}, 15\penalty0 (3):\penalty0 387--404, 2013.
\newblock \doi{10.1287/msom.2013.0429}.

\bibitem[Sturt(2026)]{sturt2026expertise}
B.~Sturt.
\newblock The value of human expertise.
\newblock arXiv preprint arXiv:2608.26051, 2026.

\bibitem[Talluri and van Ryzin(2004)]{talluri2004}
K.~Talluri and G.~van Ryzin.
\newblock Revenue management under a general discrete choice model of consumer behavior.
\newblock \emph{Management Science}, 50\penalty0 (1):\penalty0 15--33, 2004.
\newblock \doi{10.1287/mnsc.1030.0147}.

\bibitem[Train(2009)]{train2009}
K.~E. Train.
\newblock \emph{Discrete Choice Methods with Simulation}.
\newblock Cambridge University Press, second edition, 2009.

\bibitem[Yu et~al.(2026)Yu, Ma, and Zhao]{yu2026experimental}
X.~Yu, W.~Ma, and M.~Zhao.
\newblock Experimental assortments for choice estimation and nest identification.
\newblock arXiv preprint arXiv:2602.16137, 2026.

\bibitem[Zhang et~al.(2026)Zhang, Han, Shan, Blanchet, and Zhou]{zhang2026contextual}
Y.~Zhang, Y.~Han, H.~Shan, J.~Blanchet, and Z.~Zhou.
\newblock Optimal online and offline algorithms for contextual {MNL} with applications to assortment and pricing.
\newblock arXiv preprint arXiv:2604.19008, 2026.

\end{thebibliography}
}

\newpage

\begin{APPENDICES}
\section{Proofs for Section~\ref{sec:G0}}\label{app:G0}\label{sec:G0:proof}

\begin{myproof}[Proof of Theorem~\ref{thm:G0}]
We first determine the worst-case regret of a fixed assortment, then show that an optimal assortment can be chosen as a revenue prefix.
Fix an assortment $C\subseteq[n]$.

\noindent\textbf{Lower bound.} If $\alpha(C)>0$, choose $i\notin C$ with $r_i=\alpha(C)$.
Set $v_i=H$ and all remaining product weights to zero. Then
\[
R_{v}^\star-R_{v}(C)
=\frac{r_iH}{1+H}\longrightarrow\alpha(C)
\qquad\text{as }H\to\infty.
\]
If $\beta(C)>0$, let $j\in C$ have the smallest revenue in $C$;
in particular, $j\ne1$. Set $v_j=H^2$, $v_1=H$,
and all remaining product weights to zero. Since
\[
\frac{r_1H}{1+H}\le R_{v}^\star\le r_1,
\]
we have $R_{v}^\star\to r_1$, while
$R_{v}(C)\to r_j$. Thus, the regret converges to
$r_1-r_j=\beta(C)$. When either cost is zero, its lower bound
follows from the nonnegativity of regret. Therefore,
\[
\sup_{v\in\mathbb R_+^n}\bigl(R_v^\star-R_v(C)\bigr)
\ge\max\{\alpha(C),\beta(C)\}.
\]

\noindent\textbf{Upper bound.}
Fix $v\in\mathbb R_+^n$, and let $S^\star$ maximize $R_v(S)$
over all $S\subseteq[n]$, so that $R_v(S^\star)=R_v^\star$.
Set $A:=C\cap S^\star$ and $B:=C\setminus S^\star$.
Since every product in $S^\star\setminus C$ has revenue at most
$\alpha(C)$ and $R_v(A)\ge0$,
\[
R_v^\star-R_v(A)
=\frac{\sum_{i\in S^\star\setminus C}
       \bigl(r_i-R_v(A)\bigr)v_i}{1+V(S^\star)}
\le\alpha(C).
\]
If $V(B)=0$, then $R_v(C)=R_v(A)$, and the desired bound follows.
Otherwise, define
\[
\bar R_v(B):=\frac{\sum_{i\in B}r_iv_i}{V(B)}
\ge r_1-\beta(C).
\]
Using $R_v^\star\le r_1$, we obtain
\begin{align*}
R_v^\star-R_v(C)
&=\frac{1+V(A)}{1+V(C)}\bigl(R_v^\star-R_v(A)\bigr)
 +\frac{V(B)}{1+V(C)}\bigl(R_v^\star-\bar R_v(B)\bigr)\\
&\le\frac{1+V(A)}{1+V(C)}\alpha(C)
 +\frac{V(B)}{1+V(C)}\beta(C)\\
&\le\max\{\alpha(C),\beta(C)\},
\end{align*}
where the last inequality holds because the two coefficients are
nonnegative and sum to one. Together with the lower bound, this
proves the first identity.

\noindent\textbf{Optimal prefix.}
For nonempty $C$, let $k:=\max C$ be the largest product index in $C$.
Since revenues are decreasing, $r_k$ is the lowest revenue in $C$.
Moreover, $C\subseteq[k]$ and $k\in C$, so
\[
\alpha([k])=r_{k+1}\le\alpha(C),
\qquad
\beta([k])=r_1-r_k=\beta(C).
\]
Thus, the first identity implies that $[k]$ has no larger worst-case
regret than $C$. Moreover, the empty assortment has worst-case regret
$r_1$, whereas $\{1\}$ has worst-case regret $r_2<r_1$, with $r_2=0$
when $n=1$. Hence, an optimal assortment can be chosen among the
nonempty prefixes, giving
\[
\Reg_0^\star
=\min_{k\in[n]}\max\{r_{k+1},r_1-r_k\}.
\]
Every index attaining this minimum yields an optimal policy $[k]$.
\end{myproof}

\begin{myproof}[Proof of Corollary~\ref{cor:G0half}]
Recall that $k_0^\star=\max\{i\in[n]:r_i>r_1/2\}$ is the largest product index whose revenue exceeds $r_1/2$. Then
\[
r_{k_0^\star}>\frac{r_1}{2}
\qquad\text{and}\qquad
r_{k_0^\star+1}\le\frac{r_1}{2},
\]
where $r_{n+1}=0$. By Theorem~\ref{thm:G0}, the worst-case regret of a prefix $[j]$ is $\max\{r_{j+1},r_1-r_j\}$, so
it suffices to show that $k_0^\star$ minimizes this maximum over $j\in[n]$. For $j=k_0^\star$, both terms are at most $r_1/2$. If
$j<k_0^\star$, the worst-case regret of $[j]$ is at least $r_{j+1}\ge r_{k_0^\star}>r_1/2$. If $j>k_0^\star$, its worst-case regret is at least
\[
r_1-r_j\ge r_1-r_{k_0^\star+1}\ge\frac{r_1}{2}.
\]
Therefore, $[k_0^\star]$ is an optimal policy and $\Reg_0^{\star}\le r_1/2$.
\end{myproof}

\section{Proofs for Section~\ref{sec:G1}}\label{app:G1}

\subsection{Proof of Theorem~\ref{thm:one:opt}}\label{sec:G1:proof}\label{app:G1:lb}

{
For $1\le k\le p\le\ell\le n$, let $M(k,p,\ell)$ denote the maximum appearing in Theorem~\ref{thm:one:opt}:
\begin{equation}\label{eq:M}
M(k,p,\ell):=\max\left\{
 r_{\ell+1},\ r_1-r_p,\ r_{p+1}-r_\ell+\frac{r_\ell^2}{4r_k},\
 r_{k+1}-r_\ell,\ \frac{(r_1-r_\ell)^2}{4r_1}
\right\}.
\end{equation}
}
{
Let
\[
M^\star:=\min_{1\le k\le p\le\ell\le n}M(k,p,\ell).
\]
We first compare $M^\star$ with the no-observation minimax regret. For every $p\in[n]$, the third and fourth terms of $M(p,p,p)$ are at most $r_{p+1}$, and the fifth is at most $r_1-r_p$. Hence, Theorem~\ref{thm:G0} gives
\[
M^\star
\le \min_{1\le p\le n}M(p,p,p)
=\min_{1\le p\le n}\max\{r_{p+1},\,r_1-r_p\}
=\Reg_0^\star.
\]

To prove Theorem~\ref{thm:one:opt}, we establish that every policy has worst-case regret at least $M^\star$, then show that the nested policy associated with a minimizing triple achieves this bound. For the lower bound, any policy with worst-case regret at least $\Reg_0^\star$ already has regret at least $M^\star$. It therefore remains to consider policies with worst-case regret below $\Reg_0^\star$.
}

\vspace{2mm}
\noindent
{\bf Step 1: Structural lower bound.} {We first identify properties that any policy must satisfy to improve on the no-observation minimax regret, then use these properties to derive the five lower bounds in $M(k,p,\ell)$.}

\begin{lemma}\label{lem:one:lb}%
With one observation, if a policy $(S,\mathcal A)$ has worst-case regret below $\Reg_0^{{\star}}$, then there exist indices
$1\le k\le p<\ell\le n$ such that
\[
\Reg(S,\mathcal A)\ge M(k,p,\ell).
\]
\end{lemma}

{
We prove the lemma after establishing two intermediate claims. With one observation, the regret of a policy $(S,\mathcal A)$ for a fixed preference vector $v$ is
\begin{equation}\label{eq:one-gap-decomposition}
\Reg_v(S,\mathcal A)
=
\sum_{i\in S^+}\phi_v(i,S)
\bigl(R_v^\star-R_v(\mathcal A(i))\bigr).
\end{equation}
Each term in~\eqref{eq:one-gap-decomposition} is nonnegative, so keeping only some selected terms gives a lower
bound on $\Reg_v(S,\mathcal A)$.} {
To choose the three indices in Lemma~\ref{lem:one:lb}, let $U:=\bigcup_{i\in S}\mathcal A(i)$ be the union of the final assortments chosen
after a purchase. When $S$, $U$, and $\mathcal A(0)\setminus S$ are nonempty, define their largest indices by
\begin{equation}\label{eq:one-boundary-indices}
k:=\max S,
\qquad
p:=\max U,
\qquad
\ell:=\max\bigl(\mathcal A(0)\setminus S\bigr).
\end{equation}
Since revenues decrease with the product index, $r_k$, $r_p$, and $r_\ell$ are the lowest revenues in $S$, $U$, and $\mathcal A(0)\setminus S$, respectively.
The next claim shows that these indices are well defined and satisfy $k\le p<\ell$ whenever the policy's worst-case
regret is below $\Reg_0^\star$.}

\begin{claim}\label{clm:one:structure}
With one observation, suppose a policy $(S,\mathcal A)$ has worst-case regret below $\Reg_0^{{\star}}$. Then
\begin{equation}\label{eq:one-canonical-properties}
i\in\mathcal A(i)\quad\text{for every }i\in S,
\qquad
\mathcal A(0)\setminus S\ne\emptyset,
\qquad
1\in S.
\end{equation}
{
The indices in~\eqref{eq:one-boundary-indices} are well defined and satisfy
\[
1\le k\le p<\ell\le n.
\]
}
\end{claim}

\begin{myproof}[Proof of Claim~\ref{clm:one:structure}]
{
If $S=\emptyset$, the policy always commits to $\mathcal A(0)$, so its worst-case regret is at least $\Reg_0^\star$. Thus $S$ is nonempty. Suppose $j\notin\mathcal A(j)$ for some $j\in S$, and set $C:=\mathcal A(j)$.

Fix any $w\in\mathbb R_+^n$. Increasing only the preference weight of $j$ cannot decrease optimal expected revenue. To see this, choose an optimal assortment $B$ under $w$, omitting any products with zero preference weight. If $j\notin B$, increasing its weight leaves the expected revenue of $B$ unchanged. If $j\in B$, optimality implies $r_j\ge R_w(B)=R_w^\star$. In this case, the expected revenue of $B$ after the increase is a weighted average of $R_w^\star$ and $r_j$, and is therefore at least $R_w^\star$.

Set $v_j=H\ge w_j$ and $v_i=w_i$ for $i\ne j$. The preceding argument gives $R_v^\star\ge R_w^\star$, while $j\notin C$ gives $R_v(C)=R_w(C)$. As $H\to\infty$, the probability of observing $j$ satisfies
\[
\phi_v(j,S)
=\frac{H}{1+H+\sum_{i\in S\setminus\{j\}}w_i}
\longrightarrow1.
\]
Keeping only the term for observation $j$ in~\eqref{eq:one-gap-decomposition} gives
\[
\Reg(S,\mathcal A)
\ge \phi_v(j,S)\bigl(R_w^\star-R_w(C)\bigr).
\]
Letting $H\to\infty$, then taking the supremum over $w$, Theorem~\ref{thm:G0} gives
\begin{equation}\label{eq:one-canonical-own}
\Reg(S,\mathcal A)
\ge \sup_{w\in\mathbb R_+^n}\bigl(R_w^\star-R_w(C)\bigr)
=\max\{\alpha(C),\beta(C)\}
\ge\Reg_0^\star.
\end{equation}
The bound in~\eqref{eq:one-canonical-own} contradicts the assumed regret below $\Reg_0^\star$. Therefore, $i\in\mathcal A(i)$ for every $i\in S$. Hence $S\subseteq U$, so $k$ and $p$ are well defined with $k\le p$.
}

{
We next show that $\mathcal A(0)$ contains a product with index greater than $p$, and therefore with revenue below $r_p$, the lowest revenue in $U$. Let $q:=\max(U\cup\mathcal A(0))$, the largest index appearing in any final assortment, so $q\ge p$. We rule out $q=p$ using the following two bounds:
\begin{equation}\label{eq:one-canonical-boundary}
\Reg(S,\mathcal A)\ge r_{q+1},
\qquad
\Reg(S,\mathcal A)\ge r_1-r_p.
\end{equation}

For the first bound, if $q<n$, set $v_{q+1}=H$ and all other weights to zero. Since $S\subseteq U$, product $q+1$ belongs to neither the experimental assortment nor any final assortment. The observation is therefore no purchase, and the final assortment has zero expected revenue. Offering product $q+1$ alone gives expected revenue converging to $r_{q+1}$ as $H\to\infty$. If $q=n$, the bound follows from $r_{n+1}=0$.

For the second bound, the case $p=1$ is immediate. Otherwise, choose the largest index $j\in S$ such that $p\in\mathcal A(j)$. Such an index exists because $p\in U$, and $j\le p$ because $S\subseteq U$. If $p>j$, then $p\notin S$; otherwise, $p\in\mathcal A(p)$ would contradict the choice of $j$.

Starting from zero weights, give product $p$ weight $H^3$, add $H^2$ to the weight of product $j$, and add $H$ to the weight of product $1$. Coincident indices receive the sum of these weights. The largest experimental weight is on $j$, and the largest weight in $\mathcal A(j)$ is on $p$. Moreover, offering product $1$ alone gives expected revenue converging to $r_1$. Thus, as $H\to\infty$,
\[
\phi_v(j,S)\longrightarrow1,
\qquad
R_v(\mathcal A(j))\longrightarrow r_p,
\qquad
R_v^\star\longrightarrow r_1.
\]
Keeping only the term for observation $j$ in~\eqref{eq:one-gap-decomposition} proves the second bound.

If $q=p$, these bounds and Theorem~\ref{thm:G0} give
\[
\Reg(S,\mathcal A)
\ge\max\{r_{p+1},r_1-r_p\}
\ge\Reg_0^\star,
\]
a contradiction. Therefore $q>p$. Since $q>p=\max U$, product $q$ belongs to $\mathcal A(0)\setminus S$ and is its largest index. Hence $\ell=q>p$.
}

{
Finally, we show that $1\in S$. For every product $i\in[\ell]\setminus S$, we have
\begin{equation}\label{eq:one-canonical-unoffered}
\Reg(S,\mathcal A)\ge r_i-r_\ell.
\end{equation}
If $i=\ell$, the bound follows from nonnegative regret. If $i<\ell$, set $v_i=H$, $v_\ell=H^2$, and all other weights to zero. Both products lie outside $S$, so no purchase is observed with probability one. Since $\ell\in\mathcal A(0)$ has the largest weight, as $H\to\infty$, $R_v(\mathcal A(0))\to r_\ell$ and $R_v(\{i\})\to r_i$. The optimal expected revenue is at least $R_v(\{i\})$, so these limits prove~\eqref{eq:one-canonical-unoffered}.

If $1\notin S$, applying~\eqref{eq:one-canonical-unoffered} with $i=1$ and combining it with~\eqref{eq:one-canonical-boundary}, where $q=\ell$, gives
\[
\Reg(S,\mathcal A)
\ge\max\{r_{\ell+1},r_1-r_\ell\}
\ge\Reg_0^\star.
\]
The last inequality follows from Theorem~\ref{thm:G0}, contradicting the assumed regret below $\Reg_0^\star$. Therefore, $1\in S$.
}
\end{myproof}

{We next derive the five lower bounds in $M(k,p,\ell)$ from these structural properties. The following claim does not require worst-case regret below $\Reg_0^\star$, so we can also use it to obtain the lower bound for the nested policies with $p<\ell$ in Step~2.}

\begin{claim}\label{clm:one:bounds}
With one observation, suppose a policy $(S,\mathcal A)$ satisfies the properties in~{\eqref{eq:one-canonical-properties}},
and let $k,p,\ell$ be as in~{\eqref{eq:one-boundary-indices}}. If $p<\ell$, then
\[
\Reg(S,\mathcal A)\ge M(k,p,\ell).
\]
\end{claim}

\begin{myproof}[Proof of Claim~\ref{clm:one:bounds}]
{
By the assumed properties in~\eqref{eq:one-canonical-properties}, $S\subseteq U$, so $k\le p$. The definition of $\ell$ gives $\mathcal A(0)\subseteq S\cup[\ell]$, while $p=\max U<\ell$. Thus, $\ell$ is the largest index appearing in any final assortment.
The constructions establishing~\eqref{eq:one-canonical-boundary} and~\eqref{eq:one-canonical-unoffered} in the proof of Claim~\ref{clm:one:structure} use only the properties in~\eqref{eq:one-canonical-properties} and these index relations, so their bounds apply here. The two bounds in~\eqref{eq:one-canonical-boundary} therefore give the first two inequalities below. Since $k=\max S<\ell$, product $k+1$ belongs to $[\ell]\setminus S$, and~\eqref{eq:one-canonical-unoffered} gives the third:
\[
\Reg(S,\mathcal A)\ge r_{\ell+1},
\qquad
\Reg(S,\mathcal A)\ge r_1-r_p,
\qquad
\Reg(S,\mathcal A)\ge r_{k+1}-r_\ell.
\]
}
{
It remains to show the two quadratic lower bounds in~\eqref{eq:M}. Fix $z\in(0,1)$, and let $v$ denote the preference weights specified in each construction, with all unspecified weights zero.

We first show that $\Reg(S,\mathcal A)\ge r_{p+1}-r_\ell+r_\ell^2/(4r_k)$. Since $p=\max U<\ell$ and $S\subseteq U$, the products with indices $p+1$ and $\ell$ are held out of the experiment and excluded from every final assortment after a purchase. Set $v_k=z/(1-z)$. When $p+1<\ell$, set $v_{p+1}=H$ and $v_\ell=H^2$. When $p+1=\ell$, set $v_\ell=H+H^2$.

During the experiment, product $k$ is purchased with probability $z$, and no purchase is observed with probability $1-z$. The final assortment $\mathcal A(k)$ contains $k$ but neither $p+1$ nor $\ell$, while $\mathcal A(0)$ contains $\ell$. Thus,
\[
R_{v}(\mathcal A(k))=r_kz,
\qquad
\lim_{H\to\infty}R_{v}(\mathcal A(0))=r_\ell.
\]
Offering product $p+1$ alone gives expected revenue converging to $r_{p+1}$ as $H\to\infty$. Since $R_v^\star\ge R_v(\{p+1\})$, substituting the preceding observation probabilities and final-assortment revenues into~\eqref{eq:one-gap-decomposition} gives
\[
\begin{aligned}
\Reg(S,\mathcal A)
\ge\liminf_{H\to\infty}\Reg_{v}(S,\mathcal A)
\ge r_{p+1}-(1-z)r_\ell-z^2r_k
=r_{p+1}-r_\ell+zr_\ell-z^2r_k.
\end{aligned}
\]
Since $0<r_\ell<r_k$, we may choose $z=r_\ell/(2r_k)\in(0,1)$. Substituting this value gives
\[
\Reg(S,\mathcal A)
\ge r_{p+1}-r_\ell+\frac{r_\ell^2}{4r_k}.
\]
}
{
We next show that $\Reg(S,\mathcal A)\ge(r_1-r_\ell)^2/(4r_1)$. For a new construction, fix any $z\in(0,1)$ and set
\[
v_1=\frac{z}{1-z},
\qquad
v_\ell=H,
\]
with all other weights zero. Since $1\in S$ and $\ell\notin S$, product $1$ is purchased with probability $z$ during the experiment, and no purchase is observed with probability $1-z$. Because $\mathcal A(1)$ contains $1$ but not $\ell$, we have $R_{v}(\mathcal A(1))=r_1z$. For every finite $H$,~\eqref{eq:one-gap-decomposition} gives
\[
\begin{aligned}
\Reg_{v}(S,\mathcal A)
&=R_{v}^\star-r_1z^2-(1-z)R_{v}(\mathcal A(0))\\
&\ge(1-z)\bigl(r_1z-R_{v}(\mathcal A(0))\bigr),
\end{aligned}
\]
where the inequality follows from $R_{v}^\star\ge R_{v}(\{1\})=r_1z$. Since $\ell\in\mathcal A(0)$, its weight $H$ dominates the only other possibly positive weight, $v_1=z/(1-z)$, which stays fixed as $H$ grows. Thus $R_{v}(\mathcal A(0))\to r_\ell$, and
\[
\Reg(S,\mathcal A)
\ge\liminf_{H\to\infty}\Reg_{v}(S,\mathcal A)
\ge(1-z)(r_1z-r_\ell).
\]
}
The right-hand side is maximized at $z=(r_1+r_\ell)/(2r_1)$, and hence
\[
\Reg(S,\mathcal A)\ge\frac{(r_1-r_\ell)^2}{4r_1}.
\]
Together with the three bounds above, these two bounds prove the claim.
\end{myproof}

{We now proceed with the proof of Lemma~\ref{lem:one:lb}.}

\begin{myproof}[Proof of Lemma~\ref{lem:one:lb}]
By Claim~\ref{clm:one:structure}, the policy satisfies~{\eqref{eq:one-canonical-properties}}, and the indices
in~{\eqref{eq:one-boundary-indices}} satisfy $1\le k\le p<\ell\le n$. Claim~\ref{clm:one:bounds} then gives
$\Reg(S,\mathcal A)\ge M(k,p,\ell)$.
\end{myproof}

\vspace{2mm}
\noindent
{
{\bf Step 2: Matching upper bound.}
For arbitrary $1\le k\le p\le\ell\le n$, consider the policy with experimental assortment $[k]$ and final assortment $[p]$ after a purchase or $[\ell]$ after no purchase. We show that its worst-case regret is $M(k,p,\ell)$. We first show how expanding a revenue prefix changes its expected revenue. Let $[0]:=\emptyset$, so that $R_v([0])=V([0])=0$. For every $s\in[n]$,
\begin{equation}\label{eq:one-prefix-step}
R_v([s])
=
\frac{(1+V([s-1]))R_v([s-1])+r_sv_s}{1+V([s])}.
\end{equation}
Thus adding product $s$ makes the expected revenue a convex combination of $R_v([s-1])$ and $r_s$. When expanding $[a]$ to $[b]$, every added product has revenue at least $r_b$. Repeated application of~\eqref{eq:one-prefix-step} therefore gives
\[
R_v([b])\ge\min\{R_v([a]),r_b\},
\qquad 1\le a\le b\le n.
\]
The following claim gives monotonicity properties and bounds on the expected revenue of revenue prefixes. We use these properties to compare the expected revenues of the policy's final assortments with the optimal expected revenue.
}

\begin{claim}\label{lem:one:prefix}
Fix $v\in\mathbb R_{++}^n$, and let $s_v^\star:=\max\{i\in[n]:r_i\ge R_v^{{\star}}\}$, so that $[s_v^\star]$ is an optimal
assortment under $v$ \citep{talluri2004}. Then:
\begin{enumerate}
\item $R_v([s])$ is nondecreasing for $s\le s_v^\star$ and nonincreasing for $s\ge s_v^\star$.
\item $R_v([s_v^\star])\le r_{s_v^\star}$. Moreover, $R_v([s])\ge r_{s+1}$ for $s\ge s_v^\star$ and
$R_v([s])\ge r_s$ for $s\ge s_v^\star+1$.
\item For every $s\in[n]$,
\[
r_s\frac{V([s])}{1+V([s])}
\le R_v([s])
\le r_1\frac{V([s])}{1+V([s])}.
\]
\end{enumerate}
\end{claim}

\begin{myproof}[Proof of Claim~\ref{lem:one:prefix}]
By~\eqref{eq:one-prefix-step},
\[
R_v([s])-R_v([s-1])
=
\frac{v_s}{1+V([s])}\bigl(r_s-R_v([s-1])\bigr).
\]
For $s\le s_v^\star$, we have $r_s\ge R_v^{{\star}}\ge R_v([s-1])$, so the prefix revenue is nondecreasing. If $s_v^\star<n$,
then $r_{s_v^\star+1}<R_v^{{\star}}=R_v([s_v^\star])$, so it decreases at $s_v^\star+1$. {If $R_v([s])<R_v([s-1])$ and $s<n$, the convex combination in~\eqref{eq:one-prefix-step} gives $r_s<R_v([s])$. Since $r_{s+1}<r_s$, the difference identity then implies $R_v([s+1])<R_v([s])$. This proves the first statement.}

{By definition, $R_v([s_v^\star])=R_v^{{\star}}\le r_{s_v^\star}$. For $s_v^\star\le s<n$, the first statement gives $R_v([s+1])\le R_v([s])$. Since $v_{s+1}>0$, the difference identity at $s+1$ then implies $R_v([s])\ge r_{s+1}$. The case $s=n$ follows from $r_{n+1}=0$. For $s\ge s_v^\star+1$, applying the preceding bound at $s-1$ gives $R_v([s-1])\ge r_s$. The convex combination in~\eqref{eq:one-prefix-step} therefore gives $R_v([s])\ge r_s$. Finally,
\[
r_sV([s])\le\sum_{i=1}^s r_iv_i\le r_1V([s]),
\]
which proves the last statement after division by $1+V([s])$.}
\end{myproof}

We now evaluate the nested policy for every triple.

\begin{lemma}\label{lem:one:ub}%
For $1\le k\le p\le\ell\le n$, let $(S,\mathcal A)$ be the nested policy with one observation given by $S=[k]$,
$\mathcal A(i)=[p]$ for $i\in[k]$, and $\mathcal A(0)=[\ell]$. Its worst-case regret is
\[
\Reg(S,\mathcal A)=M(k,p,\ell).
\]
\end{lemma}

\begin{myproof}[Proof of Lemma~\ref{lem:one:ub}]
We prove the regret formula first for $p<\ell$ and then for $p=\ell$.
First suppose $p<\ell$. The policy satisfies~{\eqref{eq:one-canonical-properties}}. Because $U=[p]$, its indices
in~{\eqref{eq:one-boundary-indices}} are $k$, $p$ and $\ell$. Hence Claim~\ref{clm:one:bounds} gives
\begin{equation}\label{eq:one-nested-lower}
\Reg(S,\mathcal A)\ge M(k,p,\ell).
\end{equation}

{For the reverse inequality, we show that $\Reg_v(S,\mathcal A)\le M(k,p,\ell)$ for every $v\in\mathbb R_{++}^n$. Continuity of $\Reg_v(S,\mathcal A)$ in $v$ then extends this bound to $\mathbb R_+^n$. Fix such a positive vector $v$, and let $[s_v^\star]$ denote an optimal assortment under $v$ as specified in Claim~\ref{lem:one:prefix}. Let $z_v:=\frac{V([k])}{1+V([k])}$ be the purchase probability from $[k]$. By~\eqref{eq:one-gap-decomposition},}
\begin{equation}\label{eq:one-nested-gap}
\Reg_v(S,\mathcal A)
=
R_v([s_v^\star])-(1-z_v)R_v([\ell])-z_vR_v([p]).
\end{equation}
{We bound this expression on four regions of $\mathbb R_{++}^n$, defined by the position of $s_v^\star$.}

\medskip\noindent
\textbf{Region 1: $\{v\in\mathbb R_{++}^n:s_v^\star>\ell\}$.} Claim~\ref{lem:one:prefix} gives
\[
\Reg_v(S,\mathcal A)\le R_v([s_v^\star])\le r_{s_v^\star}\le r_{\ell+1}.
\]

\medskip\noindent
{
\textbf{Region 2: $\{v\in\mathbb R_{++}^n:p<s_v^\star\le\ell\}$.} Claim~\ref{lem:one:prefix} gives
\[
R_v^\star\le r_{p+1},
\qquad
R_v([p])\ge R_v([k])\ge r_kz_v.
\]
Expanding the optimal prefix $[s_v^\star]$ to $[\ell]$ gives
\[
R_v([\ell])\ge\min\{R_v^\star,r_\ell\}.
\]
Substituting these bounds into~\eqref{eq:one-nested-gap} gives
\[
\Reg_v(S,\mathcal A)
\le R_v^\star-(1-z_v)\min\{R_v^\star,r_\ell\}-r_kz_v^2.
\]
For fixed $z_v\in[0,1]$, the right-hand side is nondecreasing in $R_v^\star$. Replacing $R_v^\star$ by its upper bound $r_{p+1}$ and using $r_{p+1}\ge r_\ell$ gives
\[
\begin{aligned}
\Reg_v(S,\mathcal A)
&\le r_{p+1}-r_\ell+r_\ell z_v-r_kz_v^2\\
&\le r_{p+1}-r_\ell+\frac{r_\ell^2}{4r_k}.
\end{aligned}
\]
The last inequality follows by completing the square. This argument also includes $s_v^\star=\ell$.
}

\medskip\noindent
{
\textbf{Region 3: $\{v\in\mathbb R_{++}^n:k<s_v^\star\le p\}$.} Claim~\ref{lem:one:prefix} gives $R_v^\star\le r_{k+1}$, $R_v([\ell])\ge r_\ell$, and $R_v([p])\ge r_{p+1}$. Substituting these bounds into~\eqref{eq:one-nested-gap} gives
\[
\Reg_v(S,\mathcal A)
\le r_{k+1}-(1-z_v)r_\ell-z_vr_{p+1}
\le r_{k+1}-r_\ell,
\]
where the last inequality uses $r_{p+1}\ge r_\ell$. When $k=p$, this region is empty. In that case,
\[
r_{k+1}-r_\ell
=r_{p+1}-r_\ell
\le r_{p+1}-r_\ell+\frac{r_\ell^2}{4r_k}.
\]
Thus, the term associated with Region 3 is no larger than the bound from Region 2 and does not affect the maximum defining $M(k,p,\ell)$.
}

\medskip\noindent
{
\textbf{Region 4: $\{v\in\mathbb R_{++}^n:s_v^\star\le k\}$.} Since $[s_v^\star]\subseteq[k]$, the probability of a purchase when $[s_v^\star]$ is offered is at most $z_v$. Claim~\ref{lem:one:prefix} therefore gives $R_v^\star\le r_1z_v$. Since $\ell>p\ge s_v^\star$, it also gives $R_v([\ell])\ge r_\ell$. Expanding the optimal prefix $[s_v^\star]$ to $[p]$ gives
\[
R_v([p])\ge\min\{R_v^\star,r_p\}.
\]
Therefore,~\eqref{eq:one-nested-gap} gives
\[
\begin{aligned}
\Reg_v(S,\mathcal A)
&\le R_v^\star-z_v\min\{R_v^\star,r_p\}-(1-z_v)r_\ell\\
&=\max\left\{(1-z_v)(R_v^\star-r_\ell),\,
R_v^\star-z_vr_p-(1-z_v)r_\ell\right\}\\
&\le\max\left\{(1-z_v)(r_1z_v-r_\ell),\,
z_v(r_1-r_p)-(1-z_v)r_\ell\right\}\\
&\le\max\left\{\frac{(r_1-r_\ell)^2}{4r_1},\,r_1-r_p\right\}.
\end{aligned}
\]
The last inequality follows by maximizing the first expression over $z_v\in[0,1]$ and bounding the second by $r_1-r_p$.

These four regions cover $\mathbb R_{++}^n$, and their bounds are at most $M(k,p,\ell)$. Hence $\Reg(S,\mathcal A)\le M(k,p,\ell)$. Together with~\eqref{eq:one-nested-lower}, this proves the result when $p<\ell$.
}

If $p=\ell$, the policy commits to $[p]$ after every observation. Theorem~\ref{thm:G0} gives
\[
\Reg(S,\mathcal A)
=\max\{r_{p+1},r_1-r_p\}.
\]
The remaining three terms of $M(k,p,p)$ satisfy
\[
r_{p+1}-r_p+\frac{r_p^2}{4r_k}\le r_{p+1},
\qquad
r_{k+1}-r_p\le r_1-r_p,
\qquad
\frac{(r_1-r_p)^2}{4r_1}\le r_1-r_p.
\]
The first inequality uses $r_p\le r_k$, while the other two use $r_{k+1}\le r_1$ and
$0\le r_1-r_p\le r_1$, respectively.
Therefore,
\[
M(k,p,p)
=\max\{r_{p+1},r_1-r_p\}
=\Reg(S,\mathcal A),
\]
which proves the result when $p=\ell$.
\end{myproof}

\vspace{2mm}
\noindent
{
{\bf Step 3: Concluding the proof of Theorem~\ref{thm:one:opt}.}
For any policy $(S,\mathcal A)$ with $\Reg(S,\mathcal A)\ge\Reg_0^\star$, the comparison $M^\star\le\Reg_0^\star$ gives $\Reg(S,\mathcal A)\ge M^\star$. Otherwise, Lemma~\ref{lem:one:lb} gives
\[
\Reg(S,\mathcal A)\ge M(k,p,\ell)\ge M^\star
\]
for some $1\le k\le p<\ell\le n$. Thus every policy has worst-case regret at least $M^\star$, so $\Reg_1^\star\ge M^\star$.

By Lemma~\ref{lem:one:ub}, every triple attaining $M^\star$ defines a nested policy with worst-case regret $M^\star$, so $\Reg_1^\star\le M^\star$. Therefore, $\Reg_1^\star=M^\star$, and every minimizing triple defines an optimal nested policy, proving Theorem~\ref{thm:one:opt}.
}

{
\subsection{{Proof of Lemma~\ref{lem:one:823}}}\label{app:G1:823}

{
We first choose revenue cutoffs that give regret at most $8r_1/23$, then construct a four-product instance attaining this bound.

\noindent\textbf{Upper bound.} Let $\theta:=8/23$. By Theorem~\ref{thm:one:opt}, it suffices to find
$1\le k\le p\le\ell\le n$ such that $M(k,p,\ell)\le\theta r_1$. Define
\begin{align*}
\ell&:=\max\{i\in[n]:r_i>\theta r_1\},
\qquad
p:=\max\{i\in[n]:r_i\ge(1-\theta)r_1\},\\
k&:=\max\{i\in[n]:r_i\ge\min\{r_1,r_\ell+\theta r_1\}\}.
\end{align*}
Each defining set contains product $1$. Since $r_\ell>\theta r_1$ and $2\theta>1-\theta$,
\[
\min\{r_1,r_\ell+\theta r_1\}>(1-\theta)r_1>\theta r_1,
\]
so $k\le p\le\ell$.
}

The definitions of $p$ and $\ell$ give
\[
r_{\ell+1}\le \theta r_1,
\qquad
r_1-r_p\le \theta r_1,
\]
where $r_{n+1}=0$. By the maximality of $k$,
\[
r_{k+1}<\min\{r_1,r_\ell+\theta r_1\}\le r_\ell+\theta r_1,
\]
and hence $r_{k+1}-r_\ell<\theta r_1$. Since $r_\ell>\theta r_1$, we also have
\[
\frac{(r_1-r_\ell)^2}{4r_1}
<
\frac{(1-\theta)^2}{4}r_1
<\theta r_1.
\]

{
It remains to bound the third term. Write $x:=r_\ell$. Since $x>\theta r_1$ and $\theta<1/2$, the definition of $k$ gives
\[
r_k\ge\min\{r_1,x+\theta r_1\}>2\theta r_1.
\]
By the maximality of $p$,
\begin{equation}\label{eq:app-823-mixed-bound}
r_{p+1}-x+\frac{x^2}{4r_k}
<
(1-\theta)r_1-x+\frac{x^2}{8\theta r_1}
=:h(x).
\end{equation}
For $\theta r_1<x\le r_1$, we have $h'(x)=-1+x/(4\theta r_1)\le-9/32<0$.
Consequently, by~\eqref{eq:app-823-mixed-bound}, the third term is less than
\[
h(\theta r_1)=\left(1-2\theta+\frac{\theta}{8}\right)r_1=\theta r_1.
\]
}

Every term in $M(k,p,\ell)$ is therefore at most $\theta r_1=8r_1/23$, so Theorem~\ref{thm:one:opt} gives
$\Reg_1^{{\star}}\le8r_1/23$.

\noindent\textbf{Tightness.} Take $n=4$ and
\[
(r_1,r_2,r_3,r_4)=(23,16,15,8).
\]
We show that $M(k,p,\ell)\ge 8$ for every triple $1\le k\le p\le\ell\le4$. If $\ell\le3$,
then $r_{\ell+1}\ge r_4=8$. If $p\ge3$, then $r_1-r_p\ge r_1-r_3=8$.
It remains to consider $\ell=4$ and $p\le2$. If $p=1$, then $k=1$ and
\[
r_{p+1}-r_\ell+\frac{r_\ell^2}{4r_k}
=16-8+\frac{64}{4\cdot23}>8.
\]
If $p=2$ and $k=1$, then $r_{k+1}-r_\ell=16-8=8$. If $p=2$ and $k=2$,
then
\[
r_{p+1}-r_\ell+\frac{r_\ell^2}{4r_k}
=15-8+\frac{64}{4\cdot16}=8.
\]
Hence Theorem~\ref{thm:one:opt} gives $\Reg_1^{{\star}}\ge8$, and the upper bound gives $\Reg_1^{{\star}}\le8r_1/23=8$.

}

\section{Proofs for Section~\ref{sec:rep}}\label{app:rep}

\subsection{Proof of Theorem~\ref{thm:rep:main}}\label{sec:rep:proof}

{%
{%
We prove the lower bound and the policy guarantee in~\eqref{eq:sandwich} in three steps. First, we obtain two lower bounds for experimental assortments whose regret is below $\Reg_0^\star$. Second, we use these bounds and Theorem~\ref{thm:G0} to show that every experimental assortment has regret at least $\Reg_\infty^\star$ in~\eqref{eq:Ginf-explicit-intro}. Third, we bound the regret of the purchase-frequency policy introduced in Theorem~\ref{thm:rep:main}.

For each $W\subseteq[n]$, let $\Reg_0^{\star,W}$ denote the minimax regret with no observation when only the products in $W$ are available. In particular, $\Reg_0^{\star,[n]}=\Reg_0^\star$. We use the convention that a maximum or minimum over an empty set equals zero, so $\Reg_0^{\star,\emptyset}=0$.
\par}

\vspace{2mm}
{%
\noindent\textbf{Step 1: A lower bound for every experimental assortment.}\label{sec:rep:lb} Fix $m\ge1$ and an experimental assortment $S$. Suppose all $m$ customers purchase the same product $i\in S$. The next lemma shows that a policy which excludes $i$ from the final assortment after these observations has worst-case regret at least $\Reg_0^\star$. Thus, an optimal decision rule achieving regret below $\Reg_0^\star$ must include $i$ in the final assortment whenever all $m$ customers purchase $i$ during the experiment. We will apply this property to the lowest-revenue product in $S$ to prove
\[
\Reg_m(S)\ge r_1-\min_{j\in S}r_j
\qquad\text{when }\Reg_m(S)<\Reg_0^\star.
\]
\par}

{%
\begin{lemma}\label{lem:rep:contains}
Suppose $\Reg_m(S)<\Reg_0^{\star}$, and let $\mathcal A^\star$ be a decision rule attaining $\Reg_m(S)$. Then $i\in\mathcal A^\star(i,\ldots,i)$ for every $i\in S$. In other words, whenever all $m$ customers purchase the same product $i$ during the experiment, $\mathcal A^\star$ must include $i$ in the final assortment.
\end{lemma}
}

{%
Lemma~\ref{lem:rep:contains} generalizes the first property in Claim~\ref{clm:one:structure}, $i\in\mathcal A(i)$ for every $i\in S$, from one observation to $m$ observations.
\par}

\begin{myproof}[Proof of Lemma~\ref{lem:rep:contains}]
{%
Suppose that $i\notin\mathcal A^\star(i,\dots,i)$ for some $i\in S$, and let $C:=\mathcal A^\star(i,\dots,i)$. We show that $\Reg_m(S)\ge\alpha(C)$ and $\Reg_m(S)\ge\beta(C)$. Theorem~\ref{thm:G0} then gives $\Reg_m(S)\ge\Reg_0^\star$, contradicting the hypothesis.

Let $H>1$. In each construction, all weights not specified are zero. As $H\to\infty$, the weight of $i$ dominates every other weight in $S$, so the probability of observing $(i,\dots,i)$ converges to one. Since the loss is nonnegative,
\[
\Reg_m(S)\ge \mathbb P_v[O_t=i\text{ for every }t\in[m]]\big(R_v^{\star}-R_v(C)\big).
\]
To prove $\Reg_m(S)\ge\alpha(C)$, let $q:=\min([n]\setminus C)$, the index of the highest-revenue product excluded from $C$. Then $q\le i$ and $r_q=\alpha(C)$.
}
{%
Set $v_q=H$, and if $q<i$, also set $v_i=H^2$. Since $i,q\notin C$, all products in $C$ have zero preference weight, so $R_v(C)=0$. Offering product $q$ alone gives expected revenue $r_qH/(1+H)$. Taking $H\to\infty$ in the preceding regret bound therefore gives
\[
\Reg_m(S)\ge r_q=\alpha(C).
\]
\par}
{%
Suppose next that $C$ is nonempty, and let $p:=\max C$, the index of its lowest-revenue product. If $p=1$, then $\beta(C)=0$. For $p>1$, set $v_i=H^3$ and $v_p=H^2$; if $i>1$, also set $v_1=H$. Since $i\notin C$, every product in $C$ other than $p$ has weight at most $H$. Thus, as $H\to\infty$, $R_v(C)\to r_p$, and offering product~$1$ alone gives expected revenue tending to $r_1$. Taking this limit in the regret bound for observing $(i,\ldots,i)$ gives
\[
\Reg_m(S)\ge r_1-r_p=\beta(C).
\]
\par}
{%
If $C=\emptyset$, then $\Reg_m(S)\ge0=\beta(C)$. Combining the two bounds and applying Theorem~\ref{thm:G0} gives
\[
\Reg_m(S)\ge\max\{\alpha(C),\beta(C)\}\ge\Reg_0^\star,
\]
contradicting $\Reg_m(S)<\Reg_0^\star$. Hence $i\in\mathcal A^\star(i,\ldots,i)$ for every $i\in S$.
\par}
\end{myproof}

{%
We now apply Lemma~\ref{lem:rep:contains} to the lowest-revenue product in $S$ to obtain a regret lower bound when $\Reg_m(S)<\Reg_0^\star$. The next lemma also shows that $\Reg_m(S)$ is at least the no-observation minimax regret on the products held out of the experiment.
\par}

\begin{lemma}\label{lem:rep:lb}
For every $S\subseteq[n]$ and every $m\ge1$:
\begin{enumerate}
\item[\textup{(i)}] $\Reg_m(S)\ \ge\ \Reg_0^{\star,[n]\setminus S}$;
\item[\textup{(ii)}] if $\Reg_m(S)<\Reg_0^{\star}$, then $\Reg_m(S)\ \ge\ r_1-\min_{i\in S}r_i$.
\end{enumerate}
\end{lemma}

\begin{myproof}[Proof of Lemma~\ref{lem:rep:lb}]
We prove \textup{(i)} by setting the preference weights of all experimental products to zero, then prove \textup{(ii)} by applying Lemma~\ref{lem:rep:contains}.

{%
\textup{(i)} Fix a decision rule $\mathcal A$ with $m$ observations. Consider preference vectors satisfying $v_i=0$ for every $i\in S$, leaving the remaining preference weights unrestricted. All $m$ customers then make no purchase, so the final assortment is always $\mathcal A(0,\ldots,0)$. Removing the zero-weight products in $S$ from any assortment leaves its expected revenue unchanged. After removing these products, the rule therefore chooses the fixed assortment $\mathcal A(0,\ldots,0)\cap([n]\setminus S)$ in the no-observation problem on $[n]\setminus S$. Taking the worst case over the unrestricted preference weights and using the definition of no-observation minimax regret gives $\Reg(S,\mathcal A)\ge\Reg_0^{\star,[n]\setminus S}$. Minimizing over $\mathcal A$ proves \textup{(i)}.
\par}

{%
\textup{(ii)} Suppose $\Reg_m(S)<\Reg_0^\star$. Since $\Reg_m(\emptyset)=\Reg_0^\star$, $S$ is nonempty. Let $j:=\max S$, the index of its lowest-revenue product. If $j=1$, the required bound is $r_1-r_j=0$. Otherwise, let $\mathcal A^\star$ be a decision rule attaining $\Reg_m(S)$. For $H>0$, set $v_1=H$, $v_j=H^2$, and $v_i=0$ otherwise. Then
\[
\mathbb P_v[O_t=j\text{ for every }t\in[m]]
\geq
\left(\frac{H^2}{1+H+H^2}\right)^m.
\]
\par}
{%
Lemma~\ref{lem:rep:contains} gives $j\in\mathcal A^\star(j,\ldots,j)$ whenever all $m$ customers purchase $j$ during the experiment. Since the loss is nonnegative for every observation vector, the preceding probability bound gives, for every $H>0$,
\[
\Reg_m(S)\ge
\left(\frac{H^2}{1+H+H^2}\right)^m
\bigl(R_v^\star-R_v(\mathcal A^\star(j,\ldots,j))\bigr).
\]
As $H\to\infty$, the probability term converges to one. The final assortment contains product $j$ with weight $H^2$, while the only other positive weight is $v_1=H$, so its expected revenue tends to $r_j$. Moreover, $r_1H/(1+H)\le R_v^\star\le r_1$, so $R_v^\star\to r_1$. Taking the limit gives
\[
\Reg_m(S)\ge r_1-r_j=r_1-\min_{i\in S}r_i.
\]
\par}
\end{myproof}

\vspace{2mm}
{%
\noindent\textbf{Step 2: Minimizing the lower bound.}\label{sec:rep:minimize} We prove that $\Reg_m(S)\ge\Reg_\infty^\star$ for every experimental assortment~$S$.

We now show that removing products cannot increase the no-observation minimax regret. Let $W'\subseteq W\subseteq[n]$, and set the preference weights of products in $W\setminus W'$ to zero, leaving the weights in $W'$ unrestricted. Any fixed assortment on $W$ then earns the same expected revenue as its intersection with $W'$, and the optimal expected revenue is unchanged. Taking the worst case over these weights and minimizing over fixed assortments gives
\[
\Reg_0^{\star,W'}\le\Reg_0^{\star,W}.
\]

We next connect the no-observation regret to the definition of $\Reg_\infty^\star$ in~\eqref{eq:Ginf-explicit-intro}. For $1\le k<n$, consider the experimental assortment $[k]=\{1,\ldots,k\}$. The products held out of the experiment are $\{k+1,\ldots,n\}$. Apply Theorem~\ref{thm:G0} to the no-observation problem on these products. It gives a minimax optimal fixed assortment $\{k+1,\ldots,\ell\}$, for some $k<\ell\le n$, satisfying
\[
\Reg_0^{\star,[n]\setminus[k]}
=\max\{r_{k+1}-r_\ell,\ r_{\ell+1}\}.
\]
Since $1\le k<\ell\le n$, the minimum in~\eqref{eq:Ginf-explicit-intro} is no greater than its value at this pair. Therefore,
\[
\begin{aligned}
\Reg_\infty^\star
&\le\max\{r_1-r_k,\ r_{k+1}-r_\ell,\ r_{\ell+1}\}\\
&=\max\bigl\{r_1-r_k,\ \max\{r_{k+1}-r_\ell,\ r_{\ell+1}\}\bigr\}\\
&=\max\{r_1-r_k,\ \Reg_0^{\star,[n]\setminus[k]}\}.
\end{aligned}
\]

We now show that $\Reg_\infty^\star\le\Reg_0^\star$. This will establish $\Reg_m(S)\ge\Reg_\infty^\star$ whenever $\Reg_m(S)\ge\Reg_0^\star$. Setting $k=1$ makes $r_1-r_k=0$, so the preceding comparison and the product-removal inequality give
\[
\Reg_\infty^\star
\le\Reg_0^{\star,[n]\setminus[1]}
\le\Reg_0^\star.
\]
Consequently, for every $S$ with $\Reg_m(S)\ge\Reg_0^\star$,
\[
\Reg_m(S)\ge\Reg_0^\star\ge\Reg_\infty^\star.
\]
This includes $S=\emptyset$, since $\Reg_m(\emptyset)=\Reg_0^\star$.

It remains to consider $\Reg_m(S)<\Reg_0^\star$. Then $S$ is nonempty. Let $k:=\max S$, the index of its lowest-revenue product, so $r_k=\min_{i\in S}r_i$. If $k<n$, every product with index greater than $k$ is held out of the experiment. Thus,
\[
[n]\setminus[k]\subseteq[n]\setminus S.
\]
Lemma~\ref{lem:rep:lb} gives the first inequality below, product-removal monotonicity gives the second, and the comparison with~\eqref{eq:Ginf-explicit-intro} gives the third:
\[
\begin{aligned}
\Reg_m(S)
&\ge\max\{\Reg_0^{\star,[n]\setminus S},\ r_1-r_k\}\\
&\ge\max\{\Reg_0^{\star,[n]\setminus[k]},\ r_1-r_k\}\\
&\ge\Reg_\infty^\star.
\end{aligned}
\]

When $k=n$, Lemma~\ref{lem:rep:lb} gives $\Reg_m(S)\ge r_1-r_n$. Evaluating~\eqref{eq:Ginf-explicit-intro} at the pair $(n-1,n)$ gives $\Reg_\infty^\star\le r_1-r_{n-1}$. Hence,
\[
\Reg_m(S)\ge r_1-r_n
\ge r_1-r_{n-1}
\ge\Reg_\infty^\star.
\]
We have therefore proved $\Reg_m(S)\ge\Reg_\infty^\star$ for every $S$. Taking the minimum over $S$ yields $\Reg_m^\star\ge\Reg_\infty^\star$, the first inequality in~\eqref{eq:sandwich}.
\par}

\vspace{2mm}
{%
\noindent\textbf{Step 3: A matching policy.}\label{sec:rep:policy} We now use the purchase-frequency policy introduced in Theorem~\ref{thm:rep:main} to prove the upper bound in~\eqref{eq:sandwich}. We analyze this policy for any pair of indices $1\le k\le\ell\le n$, then choose a minimizing pair in~\eqref{eq:Ginf-explicit-intro}. With $m\ge1$ observations, the policy offers $[k]$ to $m$ customers and lets $\hat z$ be the fraction who purchase a product from this experimental assortment. Its decision rule is
\begin{equation}\label{eq:rep:threshold-policy}
\mathcal A(O_1,\ldots,O_m)=
\begin{cases}
[k],& {r_k\hat z\ge r_\ell},\\
[\ell],& {r_k\hat z<r_\ell}.
\end{cases}
\end{equation}
\par}

\begin{lemma}\label{lem:rep:binary-z-threshold}
For every $1\le k\le\ell\le n$ and $m\ge1$, the policy $([k],\mathcal A)$ in~\eqref{eq:rep:threshold-policy} satisfies
\[
\Reg_m([k])\le\Reg([k],\mathcal A)
\le\max\{r_1-r_k,\,r_{k+1}-r_\ell,\,r_{\ell+1}\}
+\frac{r_1}{2\sqrt m}.
\]
\end{lemma}

\begin{myproof}[Proof of Lemma~\ref{lem:rep:binary-z-threshold}]
{%
Fix $v\in\mathbb R_+^n$, and let
\[
z_v:=\frac{V([k])}{1+V([k])}
\]
be the probability of a purchase from the experimental assortment $[k]$. Let $\epsilon:=r_k(\hat z-z_v)$ be the signed estimation error in $r_k\hat z$. We bound the loss from each possible final assortment in terms of $|\epsilon|$. Since every product in $[k]$ has revenue between $r_k$ and $r_1$,
\[
r_kz_v\le R_v([k])\le r_1z_v.
\]

To bound $R_v^\star$, let $[s_v^\star]$ be an optimal revenue prefix under $v$, as guaranteed by Section~\ref{sec:model:choice}. If $s_v^\star\le k$, its purchase probability is at most $z_v$, so its expected revenue is at most $r_1z_v$. Otherwise, it is obtained from $[k]$ by adding products with revenues at most $r_{k+1}$. By~\eqref{eq:one-prefix-step}, adding a product makes the new expected revenue a convex combination of the current assortment's expected revenue and the added product's revenue. Since every added product has revenue at most $r_{k+1}$, repeating this procedure bounds the larger prefix's expected revenue by $\max\{R_v([k]),r_{k+1}\}$. Using $R_v([k])\le r_1z_v$, we obtain in both cases
\[
R_v^\star\le\max\{r_1z_v,\,r_{k+1}\}.
\]

We next bound the loss when the policy chooses $[\ell]$. If $s_v^\star\le\ell$, expanding it to $[\ell]$ adds products with revenues at least $r_\ell$. By~\eqref{eq:one-prefix-step}, $R_v([\ell])\ge\min\{R_v^\star,r_\ell\}$, so the loss is at most $\max\{R_v^\star-r_\ell,0\}$. If $s_v^\star>\ell$, expanding $[\ell]$ to $[s_v^\star]$ adds products with revenues at most $r_{\ell+1}$. Also,~\eqref{eq:one-prefix-step} gives $R_v^\star\le\max\{R_v([\ell]),r_{\ell+1}\}$. The loss is therefore at most $r_{\ell+1}$, because $R_v([\ell])\ge0$. Since $r_{\ell+1}\ge0$, both cases give
\[
R_v^\star-R_v([\ell])
\le\max\{R_v^\star-r_\ell,\,r_{\ell+1}\}.
\]

If $r_k\hat z\ge r_\ell$, the policy chooses $[k]$. Using $R_v([k])\ge r_kz_v$, the bound on $R_v^\star$, and $r_kz_v=r_k\hat z-\epsilon$, we obtain
\[
\begin{aligned}
R_v^\star-R_v([k])
&\le\max\{r_1-r_k,\,r_{k+1}-r_\ell+\epsilon\}\\
&\le\max\{r_1-r_k,\,r_{k+1}-r_\ell,\,r_{\ell+1}\}
+\max\{0,\epsilon\}.
\end{aligned}
\]
The first inequality also uses $0\le z_v\le1$.

If $r_k\hat z<r_\ell$, the policy chooses $[\ell]$. In this case,
\[
r_1z_v-r_\ell
=(r_1-r_k)z_v+r_k\hat z-r_\ell-\epsilon
\le r_1-r_k-\epsilon.
\]
Combining this inequality with the bounds on $R_v^\star$ and $R_v^\star-R_v([\ell])$ gives
\[
\begin{aligned}
R_v^\star-R_v([\ell])
&\le\max\{r_1-r_k-\epsilon,\,r_{k+1}-r_\ell,\,r_{\ell+1}\}\\
&\le\max\{r_1-r_k,\,r_{k+1}-r_\ell,\,r_{\ell+1}\}
+\max\{0,-\epsilon\}.
\end{aligned}
\]
Combining the two decision cases, we obtain, for every observation sequence,
\[
R_v^\star-R_v\bigl(\mathcal A(O_1,\ldots,O_m)\bigr)
\le\max\{r_1-r_k,\,r_{k+1}-r_\ell,\,r_{\ell+1}\}
+|\epsilon|.
\]

Since $m\hat z$ is binomial with parameters $m$ and $z_v$, we have $\mathbb E_v[\epsilon]=0$ and $\mathbb E_v[\epsilon^2]=r_k^2z_v(1-z_v)/m$. The Cauchy--Schwarz inequality gives
\[
\mathbb E_v[|\epsilon|]
\le\sqrt{\mathbb E_v[\epsilon^2]}
=r_k\sqrt{\frac{z_v(1-z_v)}{m}}
\le\frac{r_1}{2\sqrt m}.
\]

Taking expectations and using $\mathbb E_v[|\epsilon|]\le r_1/(2\sqrt m)$ gives a bound independent of $v$. Taking the supremum over $v$ therefore proves the stated upper bound on $\Reg([k],\mathcal A)$. Finally, $\Reg_m([k])\le\Reg([k],\mathcal A)$, since $\Reg_m([k])$ minimizes worst-case regret over decision rules using these observations.
\par}
\end{myproof}

{%
Now take any minimizing pair $(k^\star,\ell^\star)$ in~\eqref{eq:Ginf-explicit-intro}. For this pair, the maximum in Lemma~\ref{lem:rep:binary-z-threshold} equals $\Reg_\infty^\star$. Since $\Reg_m^\star$ is the minimum of $\Reg_m(S)$ over experimental assortments $S$, applying Lemma~\ref{lem:rep:binary-z-threshold} gives
\[
\Reg_m^\star
\le\Reg_m([k^\star])
\le\Reg([k^\star],\mathcal A)
\le\Reg_\infty^\star+\frac{r_1}{2\sqrt m}.
\]
\par}
{%
Together with the lower bound $\Reg_m^\star\ge\Reg_\infty^\star$ from Step~2, this proves~\eqref{eq:sandwich}. The upper and lower bounds in~\eqref{eq:sandwich} differ by $r_1/(2\sqrt m)$, which converges to zero as $m\to\infty$. Hence
\[
\lim_{m\to\infty}\Reg_m^\star
=\lim_{m\to\infty}\Reg_m([k^\star])
=\Reg_\infty^\star,
\]
which completes the proof of Theorem~\ref{thm:rep:main}.
\par}
}

\subsection{Proof of Lemma~\ref{lem:rep:one-third}}\label{app:rep:one-third}

We first verify the upper bound for the proposed cutoffs, then give a three-product instance attaining it.

\noindent\textbf{Upper bound.} Recall from~\eqref{eq:rep:two-thirds-index} that $k=\max\{i\in[n]:r_i\ge 2r_1/3\}$, and set
\[
\ell:=\min\left\{i\in[n]:r_{i+1}\le \frac{r_1}{3}\right\}.
\]
We show that each of $r_{\ell+1}$, $r_{k+1}-r_\ell$ and $r_1-r_k$ is at most $r_1/3$. Both indices are well defined, since $r_1\ge\frac{2r_1}{3}$ and $r_{n+1}=0\le\frac{r_1}{3}$. Also $k\le\ell$: otherwise
$k\ge\ell+1$, so $r_k\le r_{\ell+1}\le\frac{r_1}{3}$, contradicting $r_k\ge\frac{2r_1}{3}$. The definitions of $k$ and $\ell$ give
\[
r_1-r_k\le \frac{r_1}{3},\qquad r_{\ell+1}\le\frac{r_1}{3}.
\]
Maximality of $k$ gives $r_{k+1}<\frac{2r_1}{3}$, while $r_\ell>\frac{r_1}{3}$ follows directly if $\ell=1$ and
from minimality of $\ell$ otherwise. Therefore
\[
r_{k+1}-r_\ell<\frac{2r_1}{3}-\frac{r_1}{3}=\frac{r_1}{3}.
\]
{%
The preceding inequalities give
\[
\max\{r_1-r_k,\;r_{k+1}-r_\ell,\;r_{\ell+1}\}
\le \frac{r_1}{3}.
\]
When $k<\ell$, the definition in~\eqref{eq:Ginf-explicit-intro} gives $\Reg_\infty^\star\le r_1/3$.

When $k=\ell$, \eqref{eq:Ginf-explicit-intro} requires a pair with distinct indices. If $k<n$, use $(k,k+1)$: its maximum is $\max\{r_1-r_k,r_{k+2}\}\le r_1/3$. If $k=n$, use $(n-1,n)$: its maximum is $r_1-r_{n-1}\le r_1-r_n\le r_1/3$. Thus, $\Reg_\infty^\star\le r_1/3$ also when $k=\ell$.
}

\noindent\textbf{Tightness.} Take $n=3$ and
\[
(r_1,r_2,r_3)=(3,2,1).
\]
{%
For $n=3$, the pairs $(k,\ell)$ in~\eqref{eq:Ginf-explicit-intro} are $(1,2)$, $(1,3)$, and $(2,3)$. Their three-term maxima are, respectively,
\[
\max\{0,0,1\}
=\max\{0,1,0\}
=\max\{1,0,0\}
=1.
\]
Therefore, $\Reg_\infty^\star=1=r_1/3$.
}

\section{Proofs for Section~\ref{sec:limit}}\label{app:limit}

\subsection{Proof of Theorem~\ref{thm:vanish}}\label{app:limit:minimax}

{%
We establish the lower bound using two preference vectors with nearly identical experimental choice probabilities. We then bound the proposed policy's losses from rounding and estimation to obtain the upper bound.

\noindent\textbf{Lower bound.} Fix integers $m\ge d\ge1$, any sequence of experimental assortments $S_1,\ldots,S_m$ using at most $d$ distinct assortments, and any decision rule $\mathcal A$. We show that this policy has worst-case regret at least $\Reg_{\infty,d}^\star$. If $\Reg_{\infty,d}^\star=0$, the lower bound follows from nonnegativity of regret. Otherwise, fix any $0<\delta<\Reg_{\infty,d}^\star$. It suffices to show that this policy has worst-case regret at least $\delta$, since $\delta$ can approach $\Reg_{\infty,d}^\star$ arbitrarily closely.

We first show that we can select $d+2$ products whose revenues exceed $\delta$ and are pairwise more than $\delta$ apart. To do so, we group the revenues, including zero, by repeatedly taking the smallest unassigned value and grouping it with all unassigned values at most $\delta$ above it. Each group is consecutive and has diameter at most $\delta$. After a group is formed, every remaining unassigned revenue exceeds that group's minimum by more than $\delta$. If there were at most $d+2$ groups, \eqref{eq:distinct:partition} would imply $\Reg_{\infty,d}^\star\le\delta$, contradicting the choice of $\delta$. The constructed partition therefore has at least $d+3$ consecutive revenue groups. The first group has minimum zero, so at least $d+2$ groups have minima greater than $\delta$. Select $d+2$ of these groups, and let $p_1,\ldots,p_{d+2}$ denote the products corresponding to their minimum revenues, ordered by decreasing revenue. Setting $\rho_j:=r_{p_j}$, we obtain
\[
\rho_j-\rho_{j+1}>\delta\quad\forall\,j\in[d+1],
\qquad \rho_{d+2}>\delta.
\]

For each $t\in\{1,\ldots,m\}$, let $j_{S_t}$ record the position in the list $p_1,\ldots,p_{d+2}$ of the lowest-revenue selected product offered in $S_t$:
\[
j_{S_t}:=\max\{j:p_j\in S_t\},
\]
with $j_{S_t}:=0$ if $S_t$ contains none of the selected products. The at most $d$ distinct experimental assortments yield at most $d$ distinct values of $j_{S_t}$. Among the $d+1$ indices $1,\ldots,d+1$, there is therefore an $i$ such that $j_{S_t}\ne i$ for every $t$. Whenever $p_i\in S_t$, the definition gives $j_{S_t}\ge i$. Since equality is excluded, $j_{S_t}>i$, so $S_t$ also contains $p_{j_{S_t}}$, whose revenue is lower than that of $p_i$.

We construct two preference vectors that differ only in the weight of $p_i$. For $H>1$, set
\[
v_{p_i}^{(a)}:=0,\qquad v_{p_i}^{(b)}:=H,\qquad
v_{p_j}^{(a)}=v_{p_j}^{(b)}:=H^{j-i+1}\quad\forall\,j>i,
\]
and set all other weights to zero. If $j_{S_t}<i$, every offered product has weight zero, so the observation is a no-purchase under both vectors. If $j_{S_t}>i$, the weight of $p_{j_{S_t}}$ grows faster than every other offered weight, so this product is chosen with probability converging to one under both vectors as $H\to\infty$. Since $j_{S_t}\ne i$, these cases specify an observation vector $o^\star$, independent of $H$: its $t$th entry is $0$ when $j_{S_t}<i$ and $p_{j_{S_t}}$ when $j_{S_t}>i$. Since $m$ is fixed, all $m$ observations take these values simultaneously with probability converging to one under both preference vectors as $H\to\infty$. On observing $o^\star$, the decision rule commits to the fixed assortment $C:=\mathcal A(o^\star)$.

If $C$ contains none of $p_{i+1},\ldots,p_{d+2}$, every product in $C$ has weight zero under $v^{(a)}$, so its expected revenue is zero. Offering $p_{i+1}$ alone gives expected revenue converging to $\rho_{i+1}$ as $H\to\infty$. Thus, the limiting loss under $v^{(a)}$ is at least $\rho_{i+1}>\delta$. Otherwise, let $j^\star:=\max\{j>i:p_j\in C\}$. Under $v^{(b)}$, the weight of $p_{j^\star}$ grows faster than every other weight in $C$, so the expected revenue of $C$ converges to $\rho_{j^\star}\le\rho_{i+1}$ as $H\to\infty$. Offering $p_i$ alone gives expected revenue converging to $\rho_i$. Thus, the limiting loss under $v^{(b)}$ is at least $\rho_i-\rho_{i+1}>\delta$.

In either case, the limiting loss from committing to $C$ is greater than $\delta$ under one of the two preference vectors. Under that vector, $O=o^\star$ has probability converging to one as $H\to\infty$, and the policy commits to $C$ on this event. Since loss is nonnegative on all other observations, taking $H\to\infty$ shows that the policy has worst-case expected regret at least $\delta$. Letting $\delta$ increase to $\Reg_{\infty,d}^\star$ and minimizing over the sequence of experimental assortments and the decision rule gives
\[
\Reg_{m,d}^\star\ge\Reg_{\infty,d}^\star.
\]
\noindent\textbf{Upper bound.}
We use a minimizing partition in~\eqref{eq:distinct:partition} and the corresponding policy defined in~\eqref{eq:distinct:policy}. Let $\widetilde r_i$ be the revenue of product $i$ rounded down to its group's minimum revenue. For $v\in\mathbb R_+^n$ and $C\subseteq[n]$, define the expected rounded revenue of $C$ and the optimal expected rounded revenue by
\[
\widetilde R_v(C):=\frac{\sum_{i\in C}\widetilde r_i v_i}{1+V(C)},
\qquad
\widetilde R_v^\star:=\max_{C\subseteq[n]}\widetilde R_v(C).
\]
The chosen partition has largest group diameter $\Reg_{\infty,d}^\star$, so rounding gives $0\le r_i-\widetilde r_i\le\Reg_{\infty,d}^\star$ for every product $i$. For every assortment $C$, the corresponding decrease in expected revenue satisfies
\[
0\le R_v(C)-\widetilde R_v(C)
=\frac{\sum_{i\in C}(r_i-\widetilde r_i)v_i}{1+V(C)}
\le \frac{\Reg_{\infty,d}^\star\,V(C)}{1+V(C)}
\le \Reg_{\infty,d}^\star.
\]
Taking the maximum over assortments gives
\[
R_v^\star
=\max_{C\subseteq[n]}R_v(C)
\le \max_{C\subseteq[n]}\bigl(\widetilde R_v(C)+\Reg_{\infty,d}^\star\bigr)
=\widetilde R_v^\star+\Reg_{\infty,d}^\star.
\]
Combining this bound with $R_v(C)\ge\widetilde R_v(C)$, we obtain, for every assortment $C$,
\begin{equation}\label{eq:distinct:rounding}
R_v^\star-R_v(C)
\le \Reg_{\infty,d}^\star+\widetilde R_v^\star-R_v(C)
\le \Reg_{\infty,d}^\star+\widetilde R_v^\star-\widetilde R_v(C).
\end{equation}

Let $h\le d+1$ be the number of groups with positive rounded revenue, ordered by decreasing rounded revenue. For $j=1,\ldots,h$, let $k_j$ be the largest product index in group $j$. Under the rounded revenues, an optimal assortment can be chosen as a revenue prefix containing whole groups. Thus, if $h\ge1$, one of $[k_1],\dots,[k_h]$ is optimal. For $h\le1$, the policy defined immediately after~\eqref{eq:distinct:policy} commits to all products with positive rounded revenue, choosing the empty assortment when $h=0$. This assortment is optimal under the rounded revenues, so~\eqref{eq:distinct:rounding} bounds the policy's worst-case regret under the original revenues by $\Reg_{\infty,d}^\star$.

{%
We therefore assume $h\ge2$ and fix $v\in\mathbb R_+^n$. Write ${\widetilde R_j}$ for the expected rounded revenue of $[k_j]$. Among the experimental assortments $[k_1],\ldots,[k_{h-1}]$, let ${\widetilde R^{\max}}$ be the largest expected rounded revenue:
\[
\begin{gathered}
{\widetilde R_j}:=\widetilde R_v([k_j]),\qquad \forall\,j\in[h],\\
{\widetilde R^{\max}}:=\max_{j<h}{\widetilde R_j}.
\end{gathered}
\]
}
Every product in the last positive-revenue group, $[k_h]\setminus[k_{h-1}]$, has rounded revenue $\widetilde r_{k_h}$, so
\[
{\widetilde R_h}=
\frac{(1+V([k_{h-1}])){\widetilde R_{h-1}}
+V([k_h]\setminus[k_{h-1}])\widetilde r_{k_h}}
{1+V([k_h])}.
\]
Since ${\widetilde R_h}$ is a weighted average of ${\widetilde R_{h-1}}$ and $\widetilde r_{k_h}$, and ${\widetilde R_{h-1}}\le{\widetilde R^{\max}}$, we have ${\widetilde R_h}\le\max\{{\widetilde R^{\max}},\widetilde r_{k_h}\}$.

To bound ${\widetilde R_h}$ from below, choose $j<h$ with ${\widetilde R_j}={\widetilde R^{\max}}$. Every product in $[k_h]\setminus[k_j]$ has rounded revenue at least $\widetilde r_{k_h}$. The expected rounded revenue of $[k_h]$ is a weighted average of ${\widetilde R_j}={\widetilde R^{\max}}$ and the rounded revenues of these added products, so it is at least $\min\{{\widetilde R^{\max}},\widetilde r_{k_h}\}$. Combining the two bounds, the expected rounded revenue of the last prefix $[k_h]$ satisfies
{%
\begin{equation}\label{eq:distinct:last-prefix-bounds}
\min\{{\widetilde R^{\max}},\widetilde r_{k_h}\}
\le{\widetilde R_h}\le
\max\{{\widetilde R^{\max}},\widetilde r_{k_h}\}.
\end{equation}
}

{%
For each $j<h$, let $\epsilon_j:=\widehat R_j-{\widetilde R_j}$, and let $\epsilon:=\max\{0,\epsilon_1,\ldots,\epsilon_{h-1}\}$ be the largest positive estimation error. Fix an index $j^\star<h$ with ${\widetilde R_{j^\star}}={\widetilde R^{\max}}$, independently of the observations. Write $C:=\mathcal A(O)$ for the final assortment and $L:=\widetilde R_v^\star-\widetilde R_v(C)$ for its loss under the rounded revenues. If ${\widetilde R^{\max}}\ge\widetilde r_{k_h}$, then $\min\{{\widetilde R^{\max}},\widetilde r_{k_h}\}=\widetilde r_{k_h}$ and $\max\{{\widetilde R^{\max}},\widetilde r_{k_h}\}={\widetilde R^{\max}}$. Substituting these values into the bounds on ${\widetilde R_h}$ in~\eqref{eq:distinct:last-prefix-bounds} gives
\[
\widetilde r_{k_h}\le{\widetilde R_h}\le{\widetilde R^{\max}}.
\]
Since one of $[k_1],\ldots,[k_h]$ is optimal under the rounded revenues, we have
\[
\widetilde R_v^\star
=\max_{j\in[h]}{\widetilde R_j}
=\max\{{\widetilde R^{\max}},{\widetilde R_h}\}
={\widetilde R^{\max}}.
\] If the policy commits to $[k_{\widehat j}]$, then $\widehat R_{\widehat j}\ge\widehat R_{j^\star}$, so
\[
L={\widetilde R^{\max}}-{\widetilde R_{\widehat j}}
\le\epsilon_{\widehat j}-\epsilon_{j^\star}
\le\epsilon-\epsilon_{j^\star}.
\]
If the policy commits to $[k_h]$, then $\widehat R_{j^\star}\le\widehat R_{\widehat j}<\widetilde r_{k_h}$, so
\[
L={\widetilde R^{\max}}-{\widetilde R_h}
\le{\widetilde R^{\max}}-\widetilde r_{k_h}
<-\epsilon_{j^\star}
\le\epsilon-\epsilon_{j^\star}.
\]
If instead ${\widetilde R^{\max}}<\widetilde r_{k_h}$, the bounds in~\eqref{eq:distinct:last-prefix-bounds} give
\[
{\widetilde R^{\max}}\le{\widetilde R_h}\le\widetilde r_{k_h},
\qquad
\widetilde R_v^\star
=\max\{{\widetilde R^{\max}},{\widetilde R_h}\}
={\widetilde R_h}.
\] Committing to $[k_h]$ gives $L=0$. Committing to $[k_{\widehat j}]$ implies $\widehat R_{\widehat j}\ge\widetilde r_{k_h}$ and gives
\[
L={\widetilde R_h}-{\widetilde R_{\widehat j}}
\le\widetilde r_{k_h}-{\widetilde R_{\widehat j}}
\le\epsilon_{\widehat j}
\le\epsilon.
\]
The comparison between ${\widetilde R^{\max}}$ and $\widetilde r_{k_h}$ depends only on the fixed preference vector $v$. Since the estimates are unbiased, $\mathbb E_v[\epsilon_{j^\star}]=0$, so both cases give $\mathbb E_v[L]\le\mathbb E_v[\epsilon]$. Combining this bound with~\eqref{eq:distinct:rounding} yields
\begin{equation}\label{eq:distinct:estimation}
\mathbb E_v[R_v^\star-R_v(C)]
\le\Reg_{\infty,d}^\star+\mathbb E_v[\epsilon].
\end{equation}

For each $j<h$, $\widehat R_j$ averages $\tau=\lfloor m/d\rfloor$ independent rounded revenues $X_{j,t}:=\widetilde r_{O_{j,t}}$, $t=1,\ldots,\tau$. Each $X_{j,t}$ lies in $[0,r_1]$ and has mean ${\widetilde R_j}$. For one such observation,
\[
\operatorname{Var}_v(X_{j,t})
=\mathbb E_v[(X_{j,t}-r_1/2)^2]-({\widetilde R_j}-r_1/2)^2
\le\frac{r_1^2}{4}.
\]
The equality uses $\mathbb E_v[X_{j,t}]={\widetilde R_j}$. The inequality holds because $|X_{j,t}-r_1/2|\le r_1/2$. Consequently,
\[
\mathbb E_v[\epsilon_j^2]
=\operatorname{Var}_v(\widehat R_j)
=\frac{1}{\tau^2}\sum_{t=1}^{\tau}\operatorname{Var}_v(X_{j,t})
\le\frac{r_1^2}{4\tau},
\qquad\forall\,j<h.
\]
The first equality uses the unbiasedness of $\widehat R_j$, the second uses independence across repetitions, and the inequality applies the single-observation variance bound to all $\tau$ terms.

For the largest positive error $\epsilon$, we obtain
\begin{equation}\label{eq:distinct:second-moment-error}
\mathbb E_v[\epsilon]
\le\sqrt{\mathbb E_v[\epsilon^2]}
\le\sqrt{\sum_{j<h}\mathbb E_v[\epsilon_j^2]}
\le r_1\sqrt{\frac{h-1}{4\tau}}.
\end{equation}
The first inequality is Cauchy--Schwarz. The second uses $\epsilon^2\le\sum_{j<h}\epsilon_j^2$ and linearity of expectation. The third applies $\mathbb E_v[\epsilon_j^2]\le r_1^2/(4\tau)$ to each of the $h-1$ terms.

Alternatively, Hoeffding's lemma gives
\[
\mathbb E_v[\exp(\lambda\epsilon_j)]
\le\exp\!\left(\frac{\lambda^2r_1^2}{8\tau}\right),
\qquad \forall\,j<h,\quad\forall\,\lambda>0.
\]
Using $\exp(\lambda\epsilon)\le1+\sum_{j<h}\exp(\lambda\epsilon_j)$ and Jensen's inequality, we obtain
\begin{equation}\label{eq:distinct:hoeffding-error}
\begin{aligned}
\mathbb E_v[\epsilon]
&\le\frac{1}{\lambda}\ln\mathbb E_v[\exp(\lambda\epsilon)]
\le\frac{\ln h}{\lambda}+\frac{\lambda r_1^2}{8\tau},
\qquad \forall\,\lambda>0,\\
\mathbb E_v[\epsilon]
&\le\inf_{\lambda>0}\left\{\frac{\ln h}{\lambda}+\frac{\lambda r_1^2}{8\tau}\right\}
=r_1\sqrt{\frac{\ln h}{2\tau}}.
\end{aligned}
\end{equation}
The second line follows by minimizing over $\lambda>0$. Combining bounds~\eqref{eq:distinct:second-moment-error} and~\eqref{eq:distinct:hoeffding-error} with~\eqref{eq:distinct:estimation} and $h-1\le d$ yields
\[
\mathbb E_v[R_v^\star-R_v(C)]
\le\Reg_{\infty,d}^\star+
r_1\min\!\left\{
\frac{\sqrt d}{2\sqrt{\tau}},
\sqrt{\frac{\ln(d+1)}{2\tau}}
\right\}.
\]
The right-hand side is independent of $v$, so it bounds the policy's worst-case regret and gives the claimed upper bound on $\Reg_{m,d}^\star$. For $d=1$, $\tau=m$ and the first term in the minimum gives the error $r_1/(2\sqrt m)$ from Theorem~\ref{thm:rep:main}. For every fixed $d$, the error converges to zero as $m\to\infty$. Together with the minimax lower bound proved above, this establishes the stated limit.
}

}

{%
\subsection{Proof of Corollary~\ref{cor:distinct:universal}}\label{app:limit:corollary}

Fix integers $m\ge d\ge1$. Partition $[0,r_1]$ into $d+2$ consecutive intervals of length $r_1/(d+2)$, assigning each shared endpoint to only one interval. Group the revenue points, including zero, according to these intervals and omit empty groups. The resulting partition has at most $d+2$ nonempty consecutive groups, each with diameter at most $r_1/(d+2)$. By the definition of $\Reg_{\infty,d}^{\star}$,
\[
\Reg_{\infty,d}^{\star}\le\frac{r_1}{d+2}.
\]
Substituting this inequality into Theorem~\ref{thm:vanish} gives the stated bound.

For the observation-budget statement, take
\[
m=d\left\lceil\frac{(d+2)^2}{2}\ln(d+1)\right\rceil.
\]
Then $\lfloor m/d\rfloor\ge (d+2)^2\ln(d+1)/2$, so the logarithmic bound on the estimation term gives
\[
\frac{r_1\min\{\sqrt d,\sqrt{2\ln(d+1)}\}}
     {2\sqrt{\lfloor m/d\rfloor}}
\le r_1\sqrt{\frac{\ln(d+1)}{2\lfloor m/d\rfloor}}
\le\frac{r_1}{d+2}.
\]
Consequently, $\Reg_{m,d}^{\star}\le 2r_1/(d+2)$. This choice of $m$ is $O(d^3\ln(d+1))$, proving the observation-budget statement.
}

\subsection{Proof of Theorem~\ref{thm:distinct:vanishing}}\label{app:limit:vanishing}

{%
Fix a set $\mathcal F$ of $d$ distinct experimental assortments. For each integer $N\ge1$, it suffices to analyze $N$ repetitions of every assortment in $\mathcal F$, giving $dN$ observations. This sequence respects the limit of $N$ repetitions per assortment, and any policy permitted by the theorem can be reproduced by retaining its required observations in its experimental order. We prove necessity and sufficiency in two steps. First, we show that if the condition in the theorem fails, every policy has worst-case regret bounded below by a positive constant independent of $N$. Second, when the condition holds, we construct a policy whose worst-case regret converges to zero as $N\to\infty$.

\noindent\textbf{Necessity.} Suppose, for contradiction, that the minimum worst-case regret converges to zero as $N\to\infty$, but the condition in the theorem fails. Fix a product $i<n$ for which it fails. Either product $i$ is never offered, or every experiment containing it also contains a product with lower revenue. Fix $N$ and an arbitrary decision rule $\mathcal A_{d\times N}$. For $H>1$, define two preference vectors by
\[
v_i^{(a)}=0,\qquad v_i^{(b)}=H^i,
\]
\[
v_j^{(a)}=v_j^{(b)}=0\quad\forall\,j<i,
\qquad
v_j^{(a)}=v_j^{(b)}=H^j\quad\forall\,j>i.
\]
The vectors differ only at product $i$. Here $j$ ranges over the other products, and the exponent is simply the product's index: lower revenues correspond to larger powers of $H$.

Consider an experimental assortment $S$. If it contains a product with revenue below $r_i$, its lowest-revenue product has a weight that grows faster than every other offered weight and the no-purchase weight $1$. It is therefore selected with probability tending to one under both vectors. If $S$ contains no such product, the failed condition forces $i\notin S$. All offered weights are then zero, so the observation is a no-purchase. Let $o^\star$ be the sequence of these common limiting outcomes over the $dN$ rounds. This sequence does not depend on $H$, and for fixed $N$ its probability tends to one under both vectors.

On observing $o^\star$, the decision rule commits to the same fixed final assortment
\[
C:=\mathcal A_{d\times N}(o^\star).
\]
If $C$ contains no product with revenue below $r_i$, all its weights under $v^{(a)}$ are zero, so $R_{v^{(a)}}(C)=0$. Offering the next lower-revenue product, $i+1$, alone gives expected revenue tending to $r_{i+1}$.

Otherwise, under $v^{(b)}$, the lowest-revenue product in $C$ is selected with probability tending to one. Its revenue is at most $r_{i+1}$, whereas offering product $i$ alone gives expected revenue tending to $r_i$. The limiting loss is therefore at least $r_i-r_{i+1}$.

In either case, the loss on $O=o^\star$ has a positive limit under one of the two vectors. The probability of this event tends to one, and loss on other observations is nonnegative. Taking $H\to\infty$ therefore yields
\[
\sup_{v\ge0}\mathbb E_v\!\left[
R_v^\star-R_v\bigl(\mathcal A_{d\times N}(O)\bigr)
\right]
\ge\min\{r_{i+1},\,r_i-r_{i+1}\}>0.
\]
The right-hand side depends only on the revenues, so minimizing over decision rules preserves this positive bound for every $N$. This contradicts the assumed convergence to zero and proves necessity.

\noindent\textbf{Sufficiency.} Assume that the set $\mathcal F$ of $d$ distinct experimental assortments satisfies the condition in the theorem. We construct a policy whose experimental sequence consists of $N$ repetitions of each assortment in $\mathcal F$, and show that its worst-case regret converges to zero as $N\to\infty$.

Without loss of generality, let $S_j$ be an assortment in $\mathcal F$ satisfying $j\in S_j\subseteq[j]$, for each $j<n$. Denote the observed no-purchase and product-$j$ fractions among its $N$ observations by $\widehat\phi(0,S_j)$ and $\widehat\phi(j,S_j)$, respectively. The MNL identity $v_j=\phi_v(j,S_j)/\phi_v(0,S_j)$ lets us estimate $v_j$ when the observed no-purchase fraction is sufficiently large.

Fix $0<\delta\le1/2$, to be chosen as a function of $N$ at the end of the proof. Let ${i_\delta}$ be the first index $j<n$ whose observed no-purchase fraction satisfies $\widehat\phi(0,S_j)<{\delta}$; set ${i_\delta}=n$ if there is none. For $j<{i_\delta}$, these fractions are at least ${\delta}$, so define the estimated weight vector $\widehat v=(\widehat v_j)_{j<i_\delta}$ and its corresponding expected revenue function by
\[
\widehat v_j:=\frac{\widehat\phi(j,S_j)}{\widehat\phi(0,S_j)},
\qquad
{R_{\widehat v}}(C):=\frac{\sum_{j\in C}r_j\widehat v_j}{1+\sum_{j\in C}\widehat v_j}
\qquad\forall\,C\subseteq[{i_\delta}-1].
\]
With $[0]=\varnothing$, choose the earlier prefix with the largest estimated expected revenue:
\[
\widehat k\in\arg\max_{0\le k<{i_\delta}}{R_{\widehat v}}([k]),
\qquad
{R_{\widehat v}^{\max}}:={\max_{C\subseteq[i_\delta-1]}R_{\widehat v}(C)=R_{\widehat v}([\widehat k])},
\]
breaking ties deterministically. The decision rule $\mathcal A_{d\times N}$ commits to $[\widehat k]$ if ${R_{\widehat v}^{\max}}\ge r_{{i_\delta}}$, and to $[{i_\delta}]$ otherwise.

We now show that the final assortment $\mathcal A_{d\times N}(O)$ selected by the proposed decision rule maximizes expected revenue in the problem with estimated earlier weights. Here the product universe is $[{i_\delta}]$, with weights $\widehat v_j$ for $j<{i_\delta}$ and any nonnegative weight for product ${i_\delta}$. By revenue ordering, ${R_{\widehat v}^{\max}}$ is the best expected revenue on $[{i_\delta}-1]$. Adding product ${i_\delta}$ to any assortment in $[{i_\delta}-1]$ gives a convex combination of that assortment's expected revenue and $r_{{i_\delta}}$. When ${R_{\widehat v}^{\max}}\ge r_{{i_\delta}}$, both quantities are at most ${R_{\widehat v}^{\max}}$, so $[\widehat k]$ is optimal. If ${R_{\widehat v}^{\max}}<r_{{i_\delta}}$, every assortment's expected revenue is at most $r_{{i_\delta}}$, so adding any missing product in $[{i_\delta}]$ cannot decrease revenue. Hence, $[{i_\delta}]$ is optimal in this case.

We now bound the errors in the observed fractions and their contribution to regret under the true weights. Let ${\mathcal E_\delta}$ be the event that, for every $j<n$, both observed fractions $\widehat\phi(0,S_j)$ and $\widehat\phi(j,S_j)$ differ from their corresponding true probabilities by at most ${\delta^3}$. Each fraction is an average of $N$ independent indicators. For every $v\ge0$, $j<n$, and outcome $a\in\{0,j\}$, Hoeffding's inequality gives
\[
\mathbb P_v\!\left[
\left|\widehat\phi(a,S_j)-\phi_v(a,S_j)\right|>{\delta^3}
\right]
\le2e^{-2N{\delta^6}}.
\]
Let ${\mathcal E_\delta}^c$ denote the complement of $\mathcal E_\delta$, namely the event that at least one of these $2(n-1)$ errors exceeds ${\delta^3}$. The union bound gives
\[
\begin{aligned}
\mathbb P_v[{\mathcal E_\delta}^c]
&\le
\sum_{j=1}^{n-1}\sum_{a\in\{0,j\}}
\mathbb P_v\!\left[
\left|\widehat\phi(a,S_j)-\phi_v(a,S_j)\right|>{\delta^3}
\right]\\
&\le 2(n-1)\cdot2e^{-2N{\delta^6}}=4(n-1){e^{-2N\delta^6}}.
\end{aligned}
\]
The bound is independent of $v$, so it also holds after taking $\sup_{v\ge0}$. On ${\mathcal E_\delta}$, consider any $j<{i_\delta}$. By construction, $\widehat\phi(0,S_j)\ge{\delta}$, and $\delta\le1/2$ gives $\phi_v(0,S_j)\ge\widehat\phi(0,S_j)-{\delta^3}\ge{\delta}/2$.
Comparing the estimated and true probability ratios gives
\[
\begin{aligned}
|\widehat v_j-v_j|
&\le\frac{|\widehat\phi(j,S_j)-\phi_v(j,S_j)|}{\widehat\phi(0,S_j)}
+\frac{\phi_v(j,S_j)|\widehat\phi(0,S_j)-\phi_v(0,S_j)|}
{\widehat\phi(0,S_j)\phi_v(0,S_j)}\le{\delta^2}+2{\delta}\le3{\delta}.
\end{aligned}
\]
We next bound the change in expected revenue caused by replacing the true weights $v_j$ for $j<{i_\delta}$ with their estimates $\widehat v_j$. For the comparison, define the preference vector $v'$ by
{%
\[
v'_j:=
\begin{cases}
\widehat v_j,&j<i_\delta,\\
v_j,&j\ge i_\delta,
\end{cases}
\qquad\forall\,j\in[n].
\]
}
For any $C\subseteq[i_\delta]$, subtracting the two revenue expressions and using $\sum_{j\in C}r_jv_j=R_v(C)(1+\sum_{j\in C}v_j)$ gives
{%
\[
R_{v'}(C)-R_v(C)
=\frac{\sum_{j\in C}(r_j-R_v(C))(v'_j-v_j)}{1+\sum_{j\in C}v'_j}.
\]
}
The denominator is at least one, $|r_j-R_v(C)|\le r_1$, and $v'_j-v_j=0$ for $j\ge i_\delta$. Taking absolute values and applying the triangle inequality therefore gives, on $\mathcal E_\delta$,
\[
|{R_{v'}}(C)-R_v(C)|
\le r_1\sum_{\substack{j\in C\\j<{i_\delta}}}|\widehat v_j-v_j|
\le3(n-1)r_1{\delta}=:{\Delta_\delta},
\qquad\forall\,C\subseteq[{i_\delta}].
\]
The last inequality uses the fact that at most $n-1$ weights change, each with error at most $3\delta$.

We already showed that the final assortment $\mathcal A_{d\times N}(O)$ selected by the proposed decision rule maximizes expected revenue on $[i_\delta]$ for the estimated earlier weights and any nonnegative weight of product $i_\delta$. This applies to $v'$, giving the equality below. On $\mathcal E_\delta$, the first and last inequalities apply $|R_{v'}(C)-R_v(C)|\le\Delta_\delta$ to an assortment maximizing $R_v$ over subsets of $[i_\delta]$ and to the policy's final assortment $\mathcal A_{d\times N}(O)$, respectively:
\[
\max_{C\subseteq[{i_\delta}]}R_v(C)
\le\max_{C\subseteq[{i_\delta}]}{R_{v'}}(C)+{\Delta_\delta}
={R_{v'}}\bigl(\mathcal A_{d\times N}(O)\bigr)+{\Delta_\delta}
\le R_v\bigl(\mathcal A_{d\times N}(O)\bigr)+2{\Delta_\delta}.
\]
Thus, on ${\mathcal E_\delta}$, the final assortment loses at most $2{\Delta_\delta}$ relative to the best assortment in $[{i_\delta}]$.

It remains to compare this restricted optimum with $R_v^\star$. If ${i_\delta}=n$, they are equal. Suppose ${i_\delta}<n$. On ${\mathcal E_\delta}$, the defining inequality $\widehat\phi(0,S_{{i_\delta}})<{\delta}$ gives $\phi_v(0,S_{{i_\delta}})\le{\delta}+{\delta^3}\le2{\delta}$. Since $S_{{i_\delta}}\subseteq[{i_\delta}]$, the purchase probability from $[{i_\delta}]$ is at least $1-\phi_v(0,S_{{i_\delta}})$. Every product in $[{i_\delta}]$ has revenue at least $r_{{i_\delta}}$, so
\[
R_v([{i_\delta}])
\ge r_{{i_\delta}}\bigl(1-\phi_v(0,S_{{i_\delta}})\bigr)
\ge r_{{i_\delta}}-2r_1{\delta}.
\]
For any assortment in $[n]$, its expected revenue is a convex combination of the expected revenue of its intersection with $[{i_\delta}]$ and the revenues of its products outside $[{i_\delta}]$. Those outside revenues are below $r_{{i_\delta}}$, and the restricted optimum is at least $R_v([{i_\delta}])$. Therefore,
\[
R_v^\star
\le\max\!\left\{\max_{C\subseteq[{i_\delta}]}R_v(C),\,r_{{i_\delta}}\right\}
\le\max_{C\subseteq[{i_\delta}]}R_v(C)+2r_1{\delta}.
\]
Hence restricting the final assortment to $[{i_\delta}]$ loses at most $2r_1{\delta}$.

On ${\mathcal E_\delta}$, the final assortment's loss relative to $R_v^\star$ is therefore at most $2{\Delta_\delta}+2r_1{\delta}$. On ${\mathcal E_\delta}^c$, the loss is at most $r_1$. Taking expectations and using the bound on $\mathbb P_v[{\mathcal E_\delta}^c]$ gives
\[
\begin{aligned}
\sup_{v\ge0}\mathbb E_v\!\left[
R_v^\star-R_v\bigl(\mathcal A_{d\times N}(O)\bigr)
\right]
&\le2{\Delta_\delta}+2r_1{\delta}+r_1\sup_{v\ge0}\mathbb P_v[{\mathcal E_\delta}^c]\\
&\le(6n-4)r_1{\delta}+4r_1(n-1){e^{-2N\delta^6}}.
\end{aligned}
\]
For sufficiently large $N$, choose $\delta=(\ln N/(12N))^{1/6}$. Then $\delta\le1/2$ and $e^{-2N\delta^6}=N^{-1/6}$, so the proposed policy's worst-case regret is at most
{%
\[
r_1\left[(6n-4)\left(\frac{\ln N}{12N}\right)^{1/6}+4(n-1)N^{-1/6}\right].
\]
}%
For fixed $n$, this is $O\!\left(r_1(\ln N/N)^{1/6}\right)$, which converges to zero and proves sufficiency.
}

\section{Quantitative Bounds under Bounded Revenue Steps and Running Examples}

\subsection{Quantitative Bounds under Bounded Revenue Steps}\label{app:steps}

Here, we quantify how closely $\Reg_0^\star$ and $\Reg_1^\star$ approach their respective upper bounds $r_1/2$ and $8r_1/23$ when the revenue steps are small. Lemmas~\ref{lem:eps:G0}
and~\ref{lem:eps:G1} bound $\Reg_0^{{\star}}$ and $\Reg_1^{{\star}}$ using only the revenue steps $r_i-r_{i+1}$ at an optimal index or
triple, where $r_{n+1}=0$. These steps are at most $\epsilon$ whenever $r_i-r_{i+1}\le\epsilon$ for every $i\in[n]$.
The last step $r_n-r_{n+1}=r_n$ cannot be dropped from this condition: for $r=(1,0.9)$ and $\epsilon=0.1$, the only
other step is $r_1-r_2=0.1\le\epsilon$, but $\Reg_0^{{\star}}=\Reg_1^{{\star}}=0.1$ lies below both $\frac{r_1-\epsilon}{2}$
and $\frac{8r_1-15\epsilon}{23}$.

\begin{lemma}\label{lem:eps:G0}
Let ${\epsilon_{k_0}}:=r_{k_0^{{\star}}}-r_{k_0^{{\star}}+1}$ be the revenue step at the index $k_0^{{\star}}$ of Corollary~\ref{cor:G0half}. Then
\[
\frac{r_1-{\epsilon_{k_0}}}{2}\ \le\ \Reg_0^{{\star}}\ \le\ \frac{r_1}{2}.
\]
In particular, if ${\epsilon_{k_0}}\le\epsilon$, which holds whenever $r_i-r_{i+1}\le\epsilon$ for every $i\in[n]$, then
$\Reg_0^{{\star}}\ge\frac{r_1-\epsilon}{2}$.
\end{lemma}

\begin{myproof}[Proof of Lemma~\ref{lem:eps:G0}]
We use Corollary~\ref{cor:G0half} for the upper bound and obtain the lower bound by averaging the two terms in Theorem~\ref{thm:G0}'s regret formula. At the optimal index,
$\Reg_0^{{\star}}=\max\{r_{k_0^{{\star}}+1},\,r_1-r_{k_0^{{\star}}}\}$, so
\[
\Reg_0^{{\star}}\ge\tfrac12\big(r_{k_0^{{\star}}+1}+r_1-r_{k_0^{{\star}}}\big)=\frac{r_1-{\epsilon_{k_0}}}{2}.
\]
\end{myproof}

\begin{lemma}\label{lem:eps:G1}
Let $({k^{{\star}}},{p^{{\star}}},{\ell^{{\star}}})$ attain the minimum in Theorem~\ref{thm:one:opt}, and let
\[
{\epsilon_\ell}:=r_{{\ell^{{\star}}}}-r_{{\ell^{{\star}}}+1},
\qquad
{\epsilon_p}:=r_{{p^{{\star}}}}-r_{{p^{{\star}}}+1},
\qquad
{\epsilon_k}:=r_{{k^{{\star}}}}-r_{{k^{{\star}}}+1}
\]
be the revenue steps at these indices. Then
\[
\frac{16r_1-(13{\epsilon_\ell}+16{\epsilon_p}+{\epsilon_k})}{46}\ \le\ \Reg_1^{{\star}}\ \le\ \frac{8r_1}{23}.
\]
In particular, if ${\epsilon_\ell},{\epsilon_p},{\epsilon_k}\le\epsilon$, which holds whenever $r_i-r_{i+1}\le\epsilon$ for every $i\in[n]$, then
$\Reg_1^{{\star}}\ge\frac{8r_1-15\epsilon}{23}$. For every $0<\epsilon\le r_1/22$, there exist an integer $n$ and a revenue vector
$r_1>\cdots>r_n>0$ with $r_i-r_{i+1}\le\epsilon$ for every $i\in[n]$ and $\Reg_1^{{\star}}=\frac{8r_1-15\epsilon}{23}$. If
$\epsilon>r_1/22$, no such revenue vector exists.
\end{lemma}

\begin{myproof}[Proof of Lemma~\ref{lem:eps:G1}]
We apply Lemma~\ref{lem:one:823} for the upper bound. We then combine four regret terms for the lower bound and determine when it is attained. Write $k,p,\ell$ for ${k^{{\star}}},{p^{{\star}}},{\ell^{{\star}}}$, and let
\[
M_1:=r_{\ell+1},\quad M_2:=r_1-r_p,\quad M_3:=r_{p+1}-r_\ell+\tfrac{r_\ell^2}{4r_k},\quad M_4:=r_{k+1}-r_\ell,\quad
M_5:=\tfrac{(r_1-r_\ell)^2}{4r_1}
\]
be the five terms of $M(k,p,\ell)$ in~\eqref{eq:M}.

\noindent\textbf{Lower bound.} By Theorem~\ref{thm:one:opt}, $\Reg_1^{{\star}}=M(k,p,\ell)=\max_iM_i$. Since $\Reg_1^{{\star}}\ge M_i$ for every $i$
and $13+16+16+1=46$,
\[
46\,\Reg_1^{{\star}}\ge13M_1+16M_2+16M_3+M_4.
\]
Expanding the sum, then substituting $r_{\ell+1}=r_\ell-{\epsilon_\ell}$, $r_{p+1}=r_p-{\epsilon_p}$ and $r_{k+1}=r_k-{\epsilon_k}$, and
completing the square give
\begin{align*}
13M_1+16M_2+16M_3+M_4&=16r_1+13r_{\ell+1}-17r_\ell-16r_p+16r_{p+1}+r_{k+1}+\tfrac{4r_\ell^2}{r_k}\\
&=16r_1-4r_\ell+r_k+\tfrac{4r_\ell^2}{r_k}-(13{\epsilon_\ell}+16{\epsilon_p}+{\epsilon_k})\\
&=16r_1+\tfrac{(r_k-2r_\ell)^2}{r_k}-(13{\epsilon_\ell}+16{\epsilon_p}+{\epsilon_k}).
\end{align*}
The square is nonnegative, so
\[
\Reg_1^{{\star}}\ge\frac{16r_1-(13{\epsilon_\ell}+16{\epsilon_p}+{\epsilon_k})}{46}.
\]
If ${\epsilon_\ell},{\epsilon_p},{\epsilon_k}\le\epsilon$, this bound is at least $\frac{16r_1-30\epsilon}{46}=\frac{8r_1-15\epsilon}{23}$.

\noindent\textbf{Tightness.} Let $L:=\frac{8r_1-15\epsilon}{23}$, and suppose that $r_i-r_{i+1}\le\epsilon$ for every $i\in[n]$
and $\Reg_1^{{\star}}=L$. Then every inequality used to obtain the lower bound $L$ holds with equality. Because the coefficients
$13,16,16,1$ are positive, this gives ${\epsilon_\ell}={\epsilon_p}={\epsilon_k}=\epsilon$, $r_k=2r_\ell$, and $M_1=M_2=M_3=M_4=L$. Hence
\begin{equation}\label{eq:eps:tight-values}
r_{\ell+1}=L,\quad r_\ell=L+\epsilon,\quad r_k=2L+2\epsilon,\quad r_{k+1}=2L+\epsilon,\quad r_p=r_1-L,\quad
r_{p+1}=r_1-L-\epsilon.
\end{equation}
If $p=k$, then $r_p=r_k$ gives $r_1=3L+2\epsilon=\frac{24r_1+\epsilon}{23}$, which is impossible. Hence $k<p$, so
$r_{k+1}\ge r_p$, that is,
\[
r_{k+1}-r_p=3L+\epsilon-r_1=\frac{r_1-22\epsilon}{23}\ge0.
\]
Thus $\epsilon\le r_1/22$, so no such revenue vector exists if $\epsilon>r_1/22$.

Conversely, let $0<\epsilon\le r_1/22$. Choose indices $1<k<p<\ell$ and revenues
satisfying~\eqref{eq:eps:tight-values}, with $p=k+1$ if $\epsilon=r_1/22$. These revenues satisfy
$
r_1>r_k>r_{k+1}\ge r_p>r_{p+1}>r_\ell>r_{\ell+1}>0,
$
since $r_1-r_k=r_{p+1}-r_\ell=\frac{7r_1-16\epsilon}{23}>0$ and $r_{\ell+1}=L>0$. With enough products, the intervals
between consecutive values, and the interval from $r_{\ell+1}$ down to $0$, can be filled with revenues whose steps
are positive and at most $\epsilon$. The resulting revenue vector has $r_i-r_{i+1}\le\epsilon$ for every $i\in[n]$, so
the lower bound gives $\Reg_1^{{\star}}\ge L$. Substituting~\eqref{eq:eps:tight-values} gives $M_1=M_2=M_3=M_4=L$, and
\[
M_5=\frac{(r_1-r_\ell)^2}{4r_1}=\frac{(15r_1-8\epsilon)^2}{2116\,r_1}\le L,
\]
because this inequality is equivalent to $\epsilon\le7r_1/16$. Theorem~\ref{thm:one:opt} then gives
$\Reg_1^{{\star}}\le M(k,p,\ell)=L$, so $\Reg_1^{{\star}}=L$.
\end{myproof}

\subsection{Derivations for the Running Examples}\label{app:running-examples}

{%
Here, we derive the values, bounds, and policies in Tables~\ref{tab:examples-G0}--\ref{tab:examples-Ginf-finite} for examples (a)--(d). Throughout, $n\ge20$ is even and $\epsilon_n=R/n\le R/20$. The values of $\Reg_0^\star$ and $\Reg_1^\star$ in Tables~\ref{tab:examples-G1-finite} and~\ref{tab:examples-Ginf-finite} are exact, except for the bounds on $\Reg_1^\star$ in example (a). Theorems~\ref{thm:G0} and~\ref{thm:one:opt} give $\Reg_0^\star$ and $\Reg_1^\star$, respectively, and definition~\eqref{eq:Ginf-explicit-intro} gives $\Reg_\infty^\star$. Theorem~\ref{thm:rep:main} then gives the finite-$m$ regret bounds. Letting $n\to\infty$ with $R$ fixed gives the reported limits and, for each fixed $m$, the lower bounds on relative regret reduction.
}

\begin{table}[h]
\centering
\renewcommand{\arraystretch}{2}
\setlength{\aboverulesep}{0pt}
\setlength{\belowrulesep}{0pt}
\begin{tabular}{clc|ccc}
\hline
& Example & $\Reg_0^{{\star}}$ & $(k^{{\star}},p^{{\star}},\ell^{{\star}})$ & $\Reg_1^{{\star}}$ & $\lim\limits_{n\to\infty}\frac{\Reg_0^{{\star}}-\Reg_1^{{\star}}}{\Reg_0^{{\star}}}$ \\
\hline
(a) & Uniformly spaced revenues & $\frac{R}{2}$ & -- &
$\in\left[\frac{8R-15\epsilon_n}{23},\frac{8R}{23}\right]$ & $\frac{7}{23}$ \\
(b) & One product at $R$; others near $R/2$ & $\frac{R}{2}$ & $(1,1,n)$ &
$\epsilon_n+\frac{(\frac{R}{2}-\epsilon_n)^2}{4R}$ & $\frac{7}{8}$ \\
(c) & Clusters near $R$ and $R/3$ & $\frac{R}{3}+\frac{\epsilon_n}{2}$ & $(\frac{n}{2},\frac{n}{2},n)$ &
$\frac{(\frac{2R}{3}+\frac{\epsilon_n}{2})^2}{4R}$ & $\frac{2}{3}$ \\
(d) & Clusters near $R$ and $R/8$ & $\frac{R}{8}+\frac{\epsilon_n}{2}$ & $(\frac{n}{2},\frac{n}{2},\frac{n}{2})$
& $\frac{R}{8}+\frac{\epsilon_n}{2}$ & $0$ \\
\hline
\end{tabular}
{%
\caption{Finite-$n$ policies and regrets with one observation for examples (a)--(d), with even $n\ge20$ and $\epsilon_n=R/n$. The dash in example (a) indicates that a finite-$n$ minimizing triple is not reported. The last column gives limiting regret reduction relative to the no-observation minimax regret $\Reg_0^\star$.}
\label{tab:examples-G1-finite}
}
\end{table}

\begin{table}[h]
\centering
\renewcommand{\arraystretch}{2}
\setlength{\aboverulesep}{0pt}
\setlength{\belowrulesep}{0pt}
\resizebox{\linewidth}{!}{%
\begin{tabular}{clcc|ccc}
\hline
& Example & $\Reg_0^{{\star}}$ & $\Reg_1^{{\star}}$ & ${k^\star}$
& $\Reg_m^{{\star}}$
& \renewcommand{\arraystretch}{1}\begin{tabular}{@{}c@{}}$\liminf\limits_{n\to\infty}$\\[1mm] $\frac{\Reg_1^{{\star}}-\Reg_m^{{\star}}}{\Reg_1^{{\star}}}$\end{tabular} \\
\hline
(a) & Uniformly spaced revenues & $\frac{R}{2}$ & $\in\left[\frac{8R-15\epsilon_n}{23},\frac{8R}{23}\right]$
& $\left\lceil\frac{n+1}{3}\right\rceil$
& $\le \frac{R}{n}\left\lceil\frac{n-2}{3}\right\rceil+\frac{R}{2\sqrt m}$ & $\ge \frac{1}{24}-\frac{23}{16\sqrt m}$ \\
(b) & One product at $R$; others near $R/2$ & $\frac{R}{2}$
& $\epsilon_n+\frac{(\frac{R}{2}-\epsilon_n)^2}{4R}$ & $1$
& $\le \epsilon_n+\frac{R}{2\sqrt m}$ & $\ge 1-\frac{8}{\sqrt m}$ \\
(c) & Clusters near $R$ and $R/3$ & $\frac{R}{3}+\frac{\epsilon_n}{2}$
& $\frac{(\frac{2R}{3}+\frac{\epsilon_n}{2})^2}{4R}$ & $\frac{n}{2}$
& $\le \epsilon_n+\frac{R}{2\sqrt m}$ & $\ge 1-\frac{9}{2\sqrt m}$ \\
(d) & Clusters near $R$ and $R/8$ & $\frac{R}{8}+\frac{\epsilon_n}{2}$
& $\frac{R}{8}+\frac{\epsilon_n}{2}$ & $\frac{n}{2}$
& $\le \epsilon_n+\frac{R}{2\sqrt m}$ & $\ge 1-\frac{4}{\sqrt m}$ \\
\hline
\end{tabular}}
{%
\caption{Experimental cutoffs and regret bounds with $m$ observations for examples (a)--(d), with even $n\ge20$ and $\epsilon_n=R/n$. The last column gives asymptotic lower bounds on regret reduction relative to the one-observation minimax regret $\Reg_1^\star$.}
\label{tab:examples-Ginf-finite}
}
\end{table}

\subsubsection*{Example \textup{(a)}}

Here $r_{j+1}=(n-j)\epsilon_n$ and $R-r_j=(j-1)\epsilon_n$ for $j\in[n]$. By Theorem~\ref{thm:G0},
\[
\Reg_0^{{\star}}
=\epsilon_n\min_{j\in[n]}\max\{n-j,j-1\}
=\frac{n}{2}\,\epsilon_n
=\frac{R}{2}.
\]
For $j\le n/2$ the first term is at least $n/2$, and for $j>n/2$ the second term is at least $n/2$. Since $j=n/2$
attains $n/2$, the prefix $[n/2]$ is an optimal policy.

Every revenue step, including $r_n-r_{n+1}$, equals $\epsilon_n$. Lemma~\ref{lem:eps:G1}, applied with
$\epsilon=\epsilon_n$, therefore gives
\[
\frac{8R-15\epsilon_n}{23}\le \Reg_1^{{\star}}\le\frac{8R}{23}.
\]

{%
These bounds also determine the limiting cutoffs in Table~\ref{tab:examples-policies}. For each $n$, choose any minimizing triple $(k,p,\ell)$ in Theorem~\ref{thm:one:opt}. Let $M_1,\ldots,M_4$ denote the first four terms of $M(k,p,\ell)$ in~\eqref{eq:M}. The proof of Lemma~\ref{lem:eps:G1} and the inequalities $M_i\le\Reg_1^\star\le8R/23$ give
\[
16R-30\epsilon_n
\le13M_1+16M_2+16M_3+M_4
\le16R.
\]
The four differences $8R/23-M_i$ are nonnegative, and their weighted sum is at most $30\epsilon_n$. Each coefficient is at least one, so
\[
0\le\frac{8R}{23}-M_i\le30\epsilon_n,
\qquad {\forall\,i\in[4]}.
\]
Thus, $M_i\to8R/23$ for $i=1,\ldots,4$, since $\epsilon_n\to0$.

Using
\[
M_1=R\left(1-\frac{\ell}{n}\right),\qquad
M_2=R\frac{p-1}{n},\qquad
M_4=R\frac{\ell-k-1}{n},
\]
we obtain
\[
\frac{k}{n}\to\frac7{23},\qquad
\frac{p}{n}\to\frac8{23},\qquad
\frac{\ell}{n}\to\frac{15}{23}.
\]
}

For $\Reg_\infty^\star$, dividing the three terms in the definition~\eqref{eq:Ginf-explicit-intro} by $\epsilon_n$ gives the nonnegative integers $k-1$, $\ell-k-1$, and $n-\ell$. Their sum is $n-2$, so their maximum is at least $\lceil(n-2)/3\rceil$. Take
$k=\lceil(n+1)/3\rceil$ and $\ell=n-\lceil(n-2)/3\rceil$, which satisfy $k<\ell$ because $\ell-k=n-1-2\lceil(n-2)/3\rceil\ge n/3-1>0$ for $n\ge20$. Then $k-1=n-\ell=\lceil(n-2)/3\rceil$, so
${\ell-k-1}\mathrel{{=}}n-2-2\lceil(n-2)/3\rceil$, which is at most $\lceil(n-2)/3\rceil$. Hence
\[
\Reg_\infty^{{\star}}=\frac{R}{n}\left\lceil\frac{n-2}{3}\right\rceil.
\]

\subsubsection*{Example \textup{(b)}}

Here $r_2=R/2$ and $r_n=R/2-\epsilon_n$. With no observation, committing to $[1]$ gives worst-case regret
$\max\{r_2,R-r_1\}=R/2$ by Theorem~\ref{thm:G0}, while every longer prefix gives at least $R-r_2=R/2$. Hence
$\Reg_0^{{\star}}=R/2$, and the prefix $[1]$ is an optimal policy.

With one observation, consider the triple $(1,1,n)$. Its five terms in~\eqref{eq:M} are $0$, $0$,
$\epsilon_n+(R/2-\epsilon_n)^2/(4R)$, $\epsilon_n$, and $(R/2+\epsilon_n)^2/(4R)$. The third term is the largest: it
is at least $\epsilon_n$, and it exceeds the fifth term by $\epsilon_n/2$. Its value is at most $R/20+R/16=9R/80$.
Conversely, if $p\ge2$, the second term of $M(k,p,\ell)$ is at least $R-r_2=R/2$. If $p=1$, then $k=1$; if moreover
$\ell<n$, the first term is at least $r_n=R/2-\epsilon_n\ge9R/20$. Both bounds exceed $9R/80$, so $(1,1,n)$ is the
only minimizing triple, and
\[
\Reg_1^{{\star}}=\epsilon_n+\frac{(R/2-\epsilon_n)^2}{4R}.
\]

{%
For $\Reg_\infty^\star$, the pair $(k,\ell)=(1,n)$ separates $R$, the cluster near $R/2$, and zero into three groups with diameters $0,\epsilon_n,0$. The four revenues $R,R/2,R/2-\epsilon_n,0$ are pairwise at least $\epsilon_n$ apart, since $\epsilon_n\le R/20$. Any partition into three groups puts two of these four revenues in the same group, whose diameter is therefore at least $\epsilon_n$. Hence $\Reg_\infty^\star=\epsilon_n$.
}

\subsubsection*{Example \textup{(c)}}

Here $r_{n/2}=R-\epsilon_n$, $r_{n/2+1}=R/3+\epsilon_n/2$, and $r_n=R/3-\epsilon_n/2$. With no observation,
committing to $[n/2]$ gives worst-case regret $\max\{r_{n/2+1},R-r_{n/2}\}=R/3+\epsilon_n/2$ by
Theorem~\ref{thm:G0}. Every shorter prefix gives at least $r_{n/2}=R-\epsilon_n$, and every longer prefix gives at
least $R-r_{n/2+1}=2R/3-\epsilon_n/2$. Both bounds are at least $R/3+\epsilon_n/2$ because $\epsilon_n\le R/20$.
Hence $\Reg_0^{{\star}}=R/3+\epsilon_n/2$, and the prefix $[n/2]$ is an optimal policy.

With one observation, consider the triple $(n/2,n/2,n)$. Its five terms in~\eqref{eq:M} are $0$, $\epsilon_n$,
$\epsilon_n+(R/3-\epsilon_n/2)^2/(4(R-\epsilon_n))$, $\epsilon_n$, and $(2R/3+\epsilon_n/2)^2/(4R)$. The last term is
the largest: it is at least $R/9$, whereas $\epsilon_n\le R/20$ and
\[
\epsilon_n+\frac{(R/3-\epsilon_n/2)^2}{4(R-\epsilon_n)}
\le\frac{R}{20}+\frac{5R}{171}<\frac{R}{9}.
\]
Conversely, every triple with $\ell=n$ has the same fifth term, $(2R/3+\epsilon_n/2)^2/(4R)$. Every triple with $\ell<n$
has first term at least $r_n=R/3-\epsilon_n/2>R/4$, which exceeds this fifth term because
$(2R/3+\epsilon_n/2)^2/(4R)<(3R/4)^2/(4R)<R/4$. Hence
\[
\Reg_1^{{\star}}=\frac{(2R/3+\epsilon_n/2)^2}{4R}.
\]

{%
For $\Reg_\infty^\star$, the pair $(k,\ell)=(n/2,n)$ separates the two revenue clusters and zero into three groups with diameters $\epsilon_n,\epsilon_n,0$. Conversely, the four endpoint revenues of the clusters are pairwise at least $\epsilon_n$ apart: the endpoints of each cluster are $\epsilon_n$ apart, and the gap between the clusters is $2R/3-3\epsilon_n/2>\epsilon_n$, since $\epsilon_n\le R/20$. Any partition into three groups puts two of these four revenues in the same group, whose diameter is therefore at least $\epsilon_n$. Hence $\Reg_\infty^\star=\epsilon_n$.
}

\subsubsection*{Example \textup{(d)}}

Here $r_{n/2}=R-\epsilon_n$, $r_{n/2+1}=R/8+\epsilon_n/2$, and $r_n=R/8-\epsilon_n/2$. With no observation,
committing to $[n/2]$ gives worst-case regret $\max\{r_{n/2+1},R-r_{n/2}\}=R/8+\epsilon_n/2$ by
Theorem~\ref{thm:G0}. Every shorter prefix gives at least $r_{n/2}=R-\epsilon_n$, and every longer prefix gives at
least $R-r_{n/2+1}=7R/8-\epsilon_n/2$. Both bounds exceed $R/8+\epsilon_n/2$ because $\epsilon_n\le R/20$. Hence
$\Reg_0^{{\star}}=R/8+\epsilon_n/2$, and the prefix $[n/2]$ is an optimal policy.

{%
With one observation, the policy with $(k,p,\ell)=(n/2,n/2,n/2)$ ignores the observation and always commits to $[n/2]$. Its worst-case regret therefore equals $\Reg_0^\star=R/8+\epsilon_n/2$.
}
Conversely, if $\ell\le n/2$, the first term of $M(k,p,\ell)$ is
at least $r_{n/2+1}$. If $\ell>n/2$, then $r_\ell\le r_{n/2+1}\le3R/20$, so the fifth term satisfies
\[
\frac{(R-r_\ell)^2}{4R}\ge\frac{289R}{1600}>\frac{3R}{20}\ge r_{n/2+1}.
\]
Hence $\Reg_1^{{\star}}=R/8+\epsilon_n/2=\Reg_0^{{\star}}$.

{%
For $\Reg_\infty^\star$, the pair $(k,\ell)=(n/2,n)$ gives three group diameters $\epsilon_n,\epsilon_n,0$. Both clusters have diameter $\epsilon_n$, and the gap between them is $7R/8-3\epsilon_n/2>\epsilon_n$, since $\epsilon_n\le R/20$. Thus, the argument using the four endpoint revenues in example (c) gives $\Reg_\infty^\star\ge\epsilon_n$. Hence $\Reg_\infty^\star=\epsilon_n$.
}

{%
For each example, Theorem~\ref{thm:rep:main} and the option to ignore observations give
\[
\Reg_m^\star
\le\min\!\left\{\Reg_1^\star,\,\Reg_\infty^\star+\frac{R}{2\sqrt m}\right\}.
\]
Substituting the preceding values gives the finite-$m$ bounds in Table~\ref{tab:examples-Ginf-finite}. For the reduction bounds, the benchmarks are $\Reg_0^\star$ in Table~\ref{tab:examples-Ginf} and $\Reg_1^\star$ in Table~\ref{tab:examples-Ginf-finite}. We divide by the corresponding benchmark, subtract from one, and take the limit inferior as $n\to\infty$ through even integers, with $R$ and $m$ fixed.
}

\end{APPENDICES}
\end{document}